# A mixed Petrov–Galerkin Cosserat rod finite element formulation

Marco Herrmann[1], Domenico Castello[2], Jonas Breuling[3], Idoia Cortes Garcia[1], Leopoldo Greco[2], and Simon R. Eugster[1]

[1] Eindhoven University of Technology, Department of Mechanical Engineering-Dynamics and Control, Eindhoven, The Netherlands
[2] University of Catania, Department of Civil Engineering and Architecture (DICAr), Catania, Italy
[3] University of Stuttgart, Institute for Nonlinear Mechanics, Stuttgart, Germany

This paper presents a total Lagrangian mixed Petrov–Galerkin finite element formulation that provides a computationally efficient approach for analysing Cosserat rods that is free of singularities and locking. To achieve a singularity-free orientation parametrisation of the rod, the nodal kinematical unknowns are defined as the nodal centreline positions and unit quaternions. We apply Lagrange interpolation to all nodal kinematic coordinates, and in combination with a projection of non-unit quaternions, this leads to an interpolation with orthonormal cross-section-fixed bases. To eliminate locking effects such as shear locking, the variational Hellinger–Reissner principle is applied, resulting in a mixed approach with additional fields composed of resultant contact forces and moments. Since the mixed formulation contains the constitutive law in compliance form, it naturally incorporates constrained theories, such as the Kirchhoff–Love theory. This study specifically examines the influence of the additional internal force fields on the numerical performance, including locking mitigation and robustness. Using well-established benchmark examples, the method demonstrates enhanced computational robustness and efficiency, as evidenced by the reduction in required load steps and iterations when applying the standard Newton–Raphson method.

## 1 Introduction

The theoretical foundation of rod mechanics dates back to the early works of Euler (1744), who introduced the elastica to describe planar, shear-rigid, and inextensible rods. Later, Kirchhoff (1859) extended this theory to spatial rods, while Love (1892) generalised it to extensible cases, leading to the term Kirchhoff–Love rod for spatial, nonlinear, shear-rigid structures. Further developments in shear-deformable rod formulations were made by E. Cosserat and F. Cosserat (1909), Timoshenko (1921), Reissner (1981), and Simo and Vu-Quoc (1986). Depending on the chosen literature, these rods are referred to as special Cosserat rods (see Antman 1995), Simo–Reissner beams (see C. Meier et al. 2019), spatial Timoshenko beams (see Eugster 2015), or geometrically exact beams (see Betsch and Steinmann 2002). In this work, we use the term Cosserat rod to describe shear-deformable rods where the centerline and the cross-section orientation are the independent fields of interest.

In addition, numerous rod finite element formulations have been proposed in the literature. A comprehensive discussion of their historical evolution and subtle differences is beyond the scope of this work. For an extensive review of both shear-deformable and shear-rigid rod formulations, we refer the interested reader to the survey by C. Meier et al. (2019). However, to provide context for our approach, we highlight key developments that represent significant milestones in the progression toward the present formulation. In particular, the results of this work build upon and extend the findings of Harsch et al. (2023) and Eugster and Harsch (2023), further contributing to the ongoing development of advanced rod finite element methods.

One of the fundamental challenges in rod formulations is the description of the cross-section orientation when it comes to large rotations. As the cross-section orientation is represented by an orthonormal transformation matrix, which is an element of the special orthogonal group $SO(3)$ and can be parametrised by at least 3 real numbers, some parametrisations suffer from singularities and non-uniqueness. Examples are Euler-angles, Tait–Bryan angles or Rodrigues parameters. To avoid large rotations within one load increment or time step, early formulations employed updated Lagrangian schemes (see Cardona and Géradin 1988; Ibrahimbegović 1997; Simo and Vu-Quoc 1986). While these approaches avoid coming closer to singularities in the parametrisation, they require specific update procedures to be handled. Furthermore, they might lead to path-dependent solutions, as noticed by M. A. Crisfield and Jelenić (1999). In contrast, we pursue a total Lagrangian formulation, where all quantities are consistently defined with respect to an inertial basis. For a total Lagrangian formulation, widely known rotation parametrisations used in rod finite element formulations are rotation vectors (Cardona and Géradin 1988; Harsch et al. 2023), quaternions (Ghosh and Roy 2008; Zupan et al. 2009; Harsch and Eugster 2023; Wasmer and Betsch 2024) and the 9-parameter method (Betsch and Steinmann 2002; I. Romero and Armero 2002; Harsch et al. 2021). Each parametrisation poses additional requirement, such as a complementary update for rotation vectors, a constraint for the length of the quaternion and orthonormality constraints for the 9-parameter method.

Furthermore, an additional challenge arises in the requirement for an objective interpolation of the kinematics. M. A. Crisfield and Jelenić (1999) discovered that the direct interpolation of the nodal (total) rotation vectors leads to a non-objective formulation. Interpolating relative rotation vectors that describe the rotation relative to the first node solves this issue but introduces the additional computational step of evaluating these relative rotation vectors (see Jelenić and M. Crisfield 1999).

Another key challenge of rod finite element formulations is locking. As discussed in Balobanov and Niiranen (2018), shear and membrane locking in rods can occur if Kirchhoff (shear-rigidity) and inextensibility constraints follow in the limit case of a parameter tending to zero. This appears for instance for very slender rods if the stiffness parameters are computed in the sense of Saint–Venant by using the material's Young's and shear moduli, respectively, together with the cross-section geometry. Finite elements that are prone to locking cannot fulfill these constraints exactly over the entire element and introduce parasitic dilatation and shear strains. A common approach to mitigate this issue is reduced integration, which decreases deformation errors but fluctuations in resultant contact forces and moments are still present (see for instance C. Meier et al. 2015). A reinterpolation of the resultant contact forces and moments is possible to decrease the fluctuations. Other strategies to completely eliminate locking are B-bar methods (see for instance Greco et al. 2017), mixed formulations and three-field formulations based on the Hellinger–Reissner or the Hu–Washizu principle, respectively. Mixed and three-field formulations can be found in plate and shell finite element formulations (see Atluri 1984; Betsch and Janz 2016), Kirchhoff–Love rods (see C. Meier et al. 2014; Greco et al. 2024), or Cosserat rods (see J. G. Kim and Y. Y. Kim 1998; H. Santos et al. 2010; H. A. F. A. Santos et al. 2011). In a planar study of arches, J. G. Kim and Y. Y. Kim (1998) use a cubic interpolation of the kinematic fields and quadratic interpolation for the resultant stresses together with the Bubnov–Galerkin approach. In (H. Santos et al. 2010; H. A. F. A. Santos et al. 2011), a so-called hybrid-mixed finite element formulation was derived, resulting in different polynomial degrees for the interpolation of orientation, displacement, force and moment. Mixed formulations can also be found in isogeometric analysis of Cosserat rods, see Weeger et al. (2017); Marino (2017), where the resultant contact forces and moments are interpolated with components in the inertial fixed basis and the collocation method is used. Also in rod formulations with enhanced kinematics, where warping or extension of the cross-section is taken into account, mixed formulations are applied, see for example Choi et al. (2025); Choi et al. (2023).

Mixed and three-field methods do not only solve the problem of locking, they can also significantly improve computational performance. This has been noticed by various authors and with different applications (Magisano et al. 2017b; Magisano et al. 2017a; Pfefferkorn et al. 2021; Marino 2017; Banovec 1981). From a purely numerical perspective for solving systems of nonlinear equations, Albersmeyer and Diehl (2010) introduced the so-called lifted Newton method as a strategy to improve the convergence of Newton's method, when applied to nonlinear root finding problems. Motivated by the multiple shooting method, they highlight that a better rate of convergence can be achieved if intermediate variables are introduced as additional degrees of freedom, constrained to their forward computed value. In a mixed formulation, the independent resultant contact forces and moments take the role of these intermediate variables. Since the equivalence of the mixed and displacement-based formulations can be shown under certain circumstances, as demonstrated in (Malkus and Hughes 1978; Noor and Peters 1981), the lifted Newton method arises naturally from the mixed formulation.

Another effective way to remove locking is to employ intrinsically locking-free formulations, such as the $SE(3)$ interpolation which has piece-wise constant strains. However, the equations are highly nonlinear, involve trigonometric functions for interpolation, couple position and rotation interpolations, and have computationally expensive evaluations of the force residuals and their derivatives. Additionally, these approaches can impose restrictions on scalability, since a refinement with a higher polynomial degree is not straightforwardly possible. While the authors of this paper used the $SE(3)$ interpolation in a total Lagrangian framework (see Harsch et al. 2023), the interpolation can be done in the $SE(3)$ Lie group setting (see Sonneville et al. 2014), requiring specialised Lie group solvers. Santana et al. (2022) also presented a total Lagrangian $SE(3)$ formulation. In contrast to the work of Harsch, where a Petrov–Galerkin projection was applied on the weak form of the Cosserat rod's equilibrium equations, Santana derived a so-called equilibrium-based formulation, where the resultant contact forces and moments are obtained from the strong form of the rod's differential equations. Furthermore, strain-based formulations should be mentioned, as they also can represent vanishing strains. They are parametrised directly by the strain measures instead of positions and orientations. A well-known approach is the piece-wise constant strain formulation (see for instance Renda et al. 2018), which is widely used in the field of soft robotics. However, a major drawback of strain-based formulations is that positions and orientations can only be obtained by recursive algorithms, leading to asymmetric formulations.

Beyond these numerical challenges, there is also the question of if and how to incorporate constrained rod theories within the general Cosserat rod framework. One approach, as explored by Harsch et al. (2021), introduces explicit constraints to enforce the Kirchhoff–Love assumptions. Of course, strain-based rod formulations can fulfil the constraints without any further effort, but the previously mentioned drawback remains. The common alternative is to intrinsically satisfy the Kirchhoff–Love assumption by advanced kinematic interpolation itself. Several formulations follow this path, yet they face similar challenges as Cosserat rods do, regarding objectivity, locking, and equation complexity. Moreover, intrinsically Kirchhoff–Love rods are often formulated within an updated Lagrangian framework, inheriting the associated numerical difficulties (see for instance C. A. Meier 2016).

In this work, we address these challenges by proposing a rod finite element formulation with the following key contributions:

- We present a total Lagrangian (*i.e.* path-independent) and objective rod finite element formulation parameterised by nodal centerline positions and unit quaternions, guaranteeing a singularity-free description of finite rotations. A polynomial interpolation of these nodal quantities together with a proper mapping for non-unit quaternions leads to an objective interpolation with orthonormal cross-section-fixed bases.
- We apply a Petrov–Galerkin projection method, where the virtual centerline displacements are represented in the inertial basis and the virtual rotations are represented in the cross-section-fixed bases.
- The internal virtual work is derived using the Hellinger–Reissner principle, leading to the mixed formulation. An interpolation of the resultant contact forces and moments is applied, resulting in additional degrees of freedom of the system.
- The resulting mixed rod finite element formulation is locking-free, meaning that both displacement errors and errors in the resultant contact forces and moments are eliminated. In particular, the internal force fluctuations commonly observed in displacement-based formulations are resolved entirely.

- Constrained rod theories, such as the Kirchhoff–Love assumption, are embedded in the proposed mixed formulation. By setting parameters of the compliance to zero, all 6 strain measures can be constrained independently without the requirement of additional concepts.
- The computational benefits of the mixed formulation are demonstrated through numerical evidence from a study of the Newton–Raphson method. For all the chosen examples, solving for static equilibrium using the mixed formulation requires fewer load increments and shows better local convergence behaviour. Therefore, it requires fewer iterations per load increment in comparison to the displacement-based formulation.

The remainder of this paper is organised as follows. In Section 2, the Cosserat rod theory is briefly recapitulated. It includes the spatially continuous formulations of the internal virtual work and the external virtual work. Two formulations for the internal virtual work are introduced. The first formulation uses only kinematics quantities, while the second is derived using the variational Hellinger–Reissner principle resulting in a formulation that incorporates kinematics and independent resultant contact forces and moments. In Section 3, the Petrov–Galerkin rod finite element formulation of the Cosserat rod is derived. Therein, the ansatz functions for kinematics, the resultant contact forces and moments, and the corresponding test functions are discretised. The final rod finite element formulation is then obtained by inserting the discretised fields into the continuous formulation from Section 2, leading to the equations of the generalised force equilibrium. Section 4 is dedicated to comparing the obtained rod finite element formulations. The comparison includes the displacement-based $SE(3)$ interpolation outlined in Harsch et al. (2023), and the proposed mixed quaternion formulation of this paper with different polynomial degrees. For completeness, we also consider a mixed $SE(3)$ interpolation and displacement-based quaternion interpolations. Further, we investigate the influence of different spatial integration schemes. Well-studied benchmark experiments are conducted to show the validity and correctness of the proposed formulations. Finally, Section 5 presents the concluding remarks. In Appendix A, the boundary value problem of the Cosserat rod is formulated. Appendix B gives further details on the finite element formulation. The detailed derivation of the tangent operator to the quaternion function for non-unit quaternions is given in Appendix C. Information on the numerical convergence analyses are given in Appendix D.

## 2 Cosserat rod theory

We introduce the Euclidean 3-space $\mathbb{E}^3$ as an abstract 3-dimensional real inner-product space (see Antman 1995). A basis for $\mathbb{E}^3$ is a linearly independent set of three vectors $\boldsymbol{e}_x, \boldsymbol{e}_y, \boldsymbol{e}_z \in \mathbb{E}^3$. The basis is said to be right-handed if $\boldsymbol{e}_x \cdot (\boldsymbol{e}_y \times \boldsymbol{e}_z) > 0$ and orthonormal if their base vectors are mutually orthogonal and have unit length. In this paper, only right-handed orthonormal bases are considered. For a given basis $B = \{\boldsymbol{e}_x^B, \boldsymbol{e}_y^B, \boldsymbol{e}_z^B\}$, the respective components $a_i^B$, $i \in \{x, y, z\}$ of a vector $\boldsymbol{a} = a_x^B \boldsymbol{e}_x^B + a_y^B \boldsymbol{e}_y^B + a_z^B \boldsymbol{e}_z^B \in \mathbb{E}^3$ can be collected in the triple ${}_B\boldsymbol{a} = (a_x^B, a_y^B, a_z^B) \in \mathbb{R}^3$. Thus, we carefully distinguish $\mathbb{R}^3$ from the Euclidean 3-space $\mathbb{E}^3$. The notation $\boldsymbol{0} \in \mathbb{R}^3$ is used for the zero-triple. In contrast, the zero-tuple in the $n$-dimensional vector space $\mathbb{R}^n$ is indicated as $\boldsymbol{0}_n \in \mathbb{R}^n$, and the zero-matrix is indicated by $\boldsymbol{0}_{n\times m} \in \mathbb{R}^{n\times m}$. The rotation between two bases $B_0$ and $B_1$ is captured by the proper orthogonal transformation matrix $\boldsymbol{A}_{B_0 B_1} \in SO(3)$, which relates the coordinate representations ${}_{B_0}\boldsymbol{a}$ and ${}_{B_1}\boldsymbol{a}$ in accordance with ${}_{B_0}\boldsymbol{a} = \boldsymbol{A}_{B_0 B_1}\, {}_{B_1}\boldsymbol{a}$.

Let $\xi \in \mathcal{J} = [0, 1] \subset \mathbb{R}$ be the centerline parameter. The motion of a Cosserat rod is captured by a centerline curve represented in an inertial $I$-basis ${}_I\boldsymbol{r}_{OC} = {}_I\boldsymbol{r}_{OC}(\xi) \in \mathbb{R}^3$ augmented by the cross-section orientations $\boldsymbol{A}_{IB} = \boldsymbol{A}_{IB}(\xi) \in SO(3) = \{\boldsymbol{A} \in \mathbb{R}^{3\times3} | \boldsymbol{A}^\top \boldsymbol{A} = \boldsymbol{I} \wedge \det(\boldsymbol{A}) = 1\}$, where $\boldsymbol{I}$ refers to the identity matrix of $\mathbb{R}^3$. The subscripts $O$ and $C$ in the centerline curve refer to the origin and the centerline point, respectively. The cross-section orientation $\boldsymbol{A}_{IB}$ can also be interpreted as a transformation matrix that relates the representation of a vector in the cross-section-fixed $B$-basis to its representation in the inertial $I$-basis.

The derivatives with respect to the centerline parameter $\xi$ are denoted by $(\bullet)_{,\xi}$. The variation of a function is indicated by $\delta(\bullet)$. With this, we can introduce the virtual displacement ${}_I\delta\boldsymbol{r}_C = \delta\left({}_I\boldsymbol{r}_{OC}\right)$.

The virtual rotation of the cross-section-fixed $B$-basis relative to the inertial $I$-basis, in components with respect to the $B$-basis, is defined by ${}_B\delta\boldsymbol{\phi}_{IB} := j^{-1}\big(\boldsymbol{A}_{IB}^\top \delta\left(\boldsymbol{A}_{IB}\right)\big)$, where $j\colon \mathbb{R}^3 \to \mathfrak{so}(3) = \{\boldsymbol{B} \in \mathbb{R}^{3\times 3} \mid \boldsymbol{B}^\top = -\boldsymbol{B}\}$ is the linear and bijective map such that $\tilde{\boldsymbol{a}}\boldsymbol{b} = j(\boldsymbol{a})\boldsymbol{b} = \boldsymbol{a}\times\boldsymbol{b}$ for all $\boldsymbol{a}, \boldsymbol{b} \in \mathbb{R}^3$. Analogously, the scaled curvature is defined as ${}_B\bar{\boldsymbol{\kappa}}_{IB} := j^{-1}\big(\boldsymbol{A}_{IB}^\top(\boldsymbol{A}_{IB})_{,\xi}\big)$.

For the reference centerline curve ${}_I\boldsymbol{r}_{OC}^0$, the length of the rod's tangent vector is $J = \|{}_I\boldsymbol{r}_{OC,\xi}^0\|$. Thus, for a given centerline parameter $\xi$, the reference arc length increment is $\mathrm{d}s = J\mathrm{d}\xi$. The derivative with respect to the reference arc length $s$ of a function $\boldsymbol{f} = \boldsymbol{f}(\xi) \in \mathbb{R}^3$ can then be defined as $\boldsymbol{f}_{,s}(\xi) := \boldsymbol{f}_{,\xi}(\xi)/J(\xi)$.

The curvature is given by ${}_B\boldsymbol{\kappa}_{IB} = {}_B\bar{\boldsymbol{\kappa}}_{IB}/J$ and the dilatation and shear strains are given by ${}_B\boldsymbol{\gamma} = {}_B\bar{\boldsymbol{\gamma}}/J$ determined by ${}_B\bar{\boldsymbol{\gamma}} := \boldsymbol{A}_{IB}^\top\, {}_I\boldsymbol{r}_{OC,\xi}$. The strain measures of a Cosserat rod are the changes of these strains in comparison to the rod's reference configuration. The strain measures are then given by

$$\boldsymbol{\varepsilon}_{\boldsymbol{\gamma}} = \bar{\boldsymbol{\varepsilon}}_{\boldsymbol{\gamma}}/J\,, \quad \bar{\boldsymbol{\varepsilon}}_{\boldsymbol{\gamma}} = {}_B\bar{\boldsymbol{\gamma}} - {}_B\bar{\boldsymbol{\gamma}}^0 \quad \text{and} \quad \boldsymbol{\varepsilon}_{\boldsymbol{\kappa}} = \bar{\boldsymbol{\varepsilon}}_{\boldsymbol{\kappa}}/J\,, \quad \bar{\boldsymbol{\varepsilon}}_{\boldsymbol{\kappa}} = {}_B\bar{\boldsymbol{\kappa}}_{IB} - {}_B\bar{\boldsymbol{\kappa}}_{IB}^0\,, \tag{1}$$

where the superscript 0 refers to the evaluation in the rod's reference configuration.

We restrict ourselves to hyperelastic material models where the strain energy density with respect to the reference arc length $W = W(\boldsymbol{\varepsilon}_{\boldsymbol{\gamma}}, \boldsymbol{\varepsilon}_{\boldsymbol{\kappa}}; \xi)$ depends on the strain measures (1) and possibly explicitly on the centerline parameter $\xi$. By that, the internal virtual work functional is defined as

$$\delta W^{\mathrm{int}} := -\int_{\mathcal{J}} \delta W J \mathrm{d}\xi = -\int_{\mathcal{J}} \left\{\delta\boldsymbol{\varepsilon}_{\boldsymbol{\gamma}}^\top\, {}_B\boldsymbol{n} + \delta\boldsymbol{\varepsilon}_{\boldsymbol{\kappa}}^\top\, {}_B\boldsymbol{m}\right\} J \mathrm{d}\xi\,, \tag{2}$$

where we have recognised the constitutive equations

$${}_B\boldsymbol{n} := \left(\frac{\partial W}{\partial \boldsymbol{\varepsilon}_{\boldsymbol{\gamma}}}\right)^\top, \quad {}_B\boldsymbol{m} := \left(\frac{\partial W}{\partial \boldsymbol{\varepsilon}_{\boldsymbol{\kappa}}}\right)^\top. \tag{3}$$

Using $\delta\boldsymbol{\varepsilon}_{\boldsymbol{\gamma}} = \delta({}_B\bar{\boldsymbol{\gamma}})/J$ and $\delta\boldsymbol{\varepsilon}_{\boldsymbol{\kappa}} = \delta({}_B\bar{\boldsymbol{\kappa}}_{IB})/J$ with $\delta({}_B\bar{\boldsymbol{\gamma}}) = \boldsymbol{A}_{IB}^\top({}_I\delta\boldsymbol{r}_C)_{,\xi} - {}_B\delta\boldsymbol{\phi}_{IB} \times {}_B\bar{\boldsymbol{\gamma}}$ and $\delta({}_B\bar{\boldsymbol{\kappa}}_{IB}) = ({}_B\delta\boldsymbol{\phi}_{IB})_{,\xi} - {}_B\delta\boldsymbol{\phi}_{IB} \times {}_B\bar{\boldsymbol{\kappa}}_{IB}$ (for a detailed derivation, see Harsch et al. 2023, Appendix B), the internal virtual work (2) takes the form

$$\delta W^{\mathrm{int}} = -\int_{\mathcal{J}} \left\{({}_I\delta\boldsymbol{r}_C)_{,\xi}^\top \boldsymbol{A}_{IB}\, {}_B\boldsymbol{n} + ({}_B\delta\boldsymbol{\phi}_{IB})_{,\xi}^\top\, {}_B\boldsymbol{m} - {}_B\delta\boldsymbol{\phi}_{IB}^\top \left[{}_B\bar{\boldsymbol{\gamma}} \times {}_B\boldsymbol{n} + {}_B\bar{\boldsymbol{\kappa}}_{IB} \times {}_B\boldsymbol{m}\right]\right\} \mathrm{d}\xi\,. \tag{4}$$

As in Simo and Vu-Quoc (1986, Equation 2.10), we introduce the diagonal elasticity matrices $\boldsymbol{C}_{\boldsymbol{\gamma}} = \mathrm{diag}(k_{\mathrm{e}}, k_{\mathrm{s}_y}, k_{\mathrm{s}_z})$ and $\boldsymbol{C}_{\boldsymbol{\kappa}} = \mathrm{diag}(k_{\mathrm{t}}, k_{\mathrm{b}_y}, k_{\mathrm{b}_z})$ with constant coefficients. In the following, the simple quadratic strain energy density

$$W(\boldsymbol{\varepsilon}_{\boldsymbol{\gamma}}, \boldsymbol{\varepsilon}_{\boldsymbol{\kappa}}) = \frac{1}{2}\boldsymbol{\varepsilon}_{\boldsymbol{\gamma}}^\top \boldsymbol{C}_{\boldsymbol{\gamma}}\, \boldsymbol{\varepsilon}_{\boldsymbol{\gamma}} + \frac{1}{2}\boldsymbol{\varepsilon}_{\boldsymbol{\kappa}}^\top \boldsymbol{C}_{\boldsymbol{\kappa}}\, \boldsymbol{\varepsilon}_{\boldsymbol{\kappa}} \tag{5}$$

is used. The complementary quadratic strain energy density $W^\star$ (considering a stress-free reference configuration, *i.e.*, reference contact forces and moments ${}_B\boldsymbol{n}^0 = {}_B\boldsymbol{m}^0 = \boldsymbol{0}$) is given by

$$W^\star({}_B\boldsymbol{n}, {}_B\boldsymbol{m}) = \frac{1}{2}{}_B\boldsymbol{n}^\top \boldsymbol{C}_{\boldsymbol{\gamma}}^{-1}\, {}_B\boldsymbol{n} + \frac{1}{2}{}_B\boldsymbol{m}^\top \boldsymbol{C}_{\boldsymbol{\kappa}}^{-1}\, {}_B\boldsymbol{m}\,. \tag{6}$$

Between (5) and (6), there exists the following classic relation of Fenchel's equation:

$$W(\boldsymbol{\varepsilon}_{\boldsymbol{\gamma}}, \boldsymbol{\varepsilon}_{\boldsymbol{\kappa}}) + W^\star({}_B\boldsymbol{n}, {}_B\boldsymbol{m}) = {}_B\boldsymbol{n}^\top\boldsymbol{\varepsilon}_{\boldsymbol{\gamma}} + {}_B\boldsymbol{m}^\top\boldsymbol{\varepsilon}_{\boldsymbol{\kappa}}\,. \tag{7}$$

From that relation, the strain energy density considered in the two-field Hellinger–Reissner functional follows as

$$W = W^{\mathrm{HR}}\left(\boldsymbol{\varepsilon}_{\boldsymbol{\gamma}}, \boldsymbol{\varepsilon}_{\boldsymbol{\kappa}}, {}_B\boldsymbol{n}, {}_B\boldsymbol{m}\right) = -W^\star({}_B\boldsymbol{n}, {}_B\boldsymbol{m}) + {}_B\boldsymbol{n}^\top\boldsymbol{\varepsilon}_{\boldsymbol{\gamma}} + {}_B\boldsymbol{m}^\top\boldsymbol{\varepsilon}_{\boldsymbol{\kappa}}\,. \tag{8}$$

The internal virtual work functional becomes

$$\begin{aligned}\delta W^{\mathrm{int,HR}} := -\int_{\mathcal{J}} \delta W^{\mathrm{HR}} J \mathrm{d}\xi = \int_{\mathcal{J}} &\left\{\delta({}_B\boldsymbol{n})^\top \left(\frac{\partial W^\star}{\partial\, {}_B\boldsymbol{n}}\right)^\top + \delta({}_B\boldsymbol{m})^\top \left(\frac{\partial W^\star}{\partial\, {}_B\boldsymbol{m}}\right)^\top\right\} J \mathrm{d}\xi \\ &- \int_{\mathcal{J}} \left\{\delta({}_B\boldsymbol{n})^\top\boldsymbol{\varepsilon}_{\boldsymbol{\gamma}} + {}_B\boldsymbol{n}^\top\, \delta\boldsymbol{\varepsilon}_{\boldsymbol{\gamma}} + \delta({}_B\boldsymbol{m})^\top\boldsymbol{\varepsilon}_{\boldsymbol{\kappa}} + {}_B\boldsymbol{m}^\top\, \delta\boldsymbol{\varepsilon}_{\boldsymbol{\kappa}}\right\} J \mathrm{d}\xi\,,\end{aligned} \tag{9}$$

where it is notable that ${}_B\delta\boldsymbol{n} \neq \delta({}_B\boldsymbol{n}) = \delta(\boldsymbol{A}_{IB}^\top {}_I\boldsymbol{n}) = {}_B\delta\boldsymbol{n} - {}_B\delta\boldsymbol{\phi}_{IB} \times {}_B\boldsymbol{n}$ and similarly for the variation of ${}_B\boldsymbol{m}$. Ordering the terms in the integral by the variations and using $\boldsymbol{\varepsilon}_{\boldsymbol{\gamma}} = \bar{\boldsymbol{\varepsilon}}_{\boldsymbol{\gamma}}/J$ and $\boldsymbol{\varepsilon}_{\boldsymbol{\kappa}} = \bar{\boldsymbol{\varepsilon}}_{\boldsymbol{\kappa}}/J$ leads to

$$\begin{aligned}\delta W^{\mathrm{int,HR}} = -\int_{\mathcal{J}} \bigg\{ \delta({}_B\boldsymbol{n})^\top \bigg[ \bar{\boldsymbol{\varepsilon}}_{\boldsymbol{\gamma}} - J\Big(\frac{\partial W^\star}{\partial {}_B\boldsymbol{n}}\Big)^\top \bigg] + \delta({}_B\boldsymbol{m})^\top \bigg[ \bar{\boldsymbol{\varepsilon}}_{\boldsymbol{\kappa}} - J\Big(\frac{\partial W^\star}{\partial {}_B\boldsymbol{m}}\Big)^\top \bigg] \bigg\} \mathrm{d}\xi \\ - \int_{\mathcal{J}} \Big\{ \delta\boldsymbol{\varepsilon}_{\boldsymbol{\gamma}}^\top {}_B\boldsymbol{n} + \delta\boldsymbol{\varepsilon}_{\boldsymbol{\kappa}}^\top {}_B\boldsymbol{m} \Big\} J \mathrm{d}\xi \,. \end{aligned} \tag{10}$$

Note that the second line reads equivalently to (2) with the only difference that the resultant contact forces and moments in (2) are obtained by the constitutive equations (3). Here, the resultant contact forces and moments are independent fields. Using the same variations as introduced after (3) the second line can be transformed identically. For the first line, using the derivatives of (6) w.r.t. ${}_B\boldsymbol{n}$ and ${}_B\boldsymbol{m}$ results in

$$\begin{aligned}\delta W^{\mathrm{int,HR}} = -\int_{\mathcal{J}} \Big\{ \delta({}_B\boldsymbol{n})^\top \Big[ \bar{\boldsymbol{\varepsilon}}_{\boldsymbol{\gamma}} - J\,\boldsymbol{C}_{\boldsymbol{\gamma}}^{-1} {}_B\boldsymbol{n} \Big] + \delta({}_B\boldsymbol{m})^\top \Big[ \bar{\boldsymbol{\varepsilon}}_{\boldsymbol{\kappa}} - J\,\boldsymbol{C}_{\boldsymbol{\kappa}}^{-1} {}_B\boldsymbol{m} \Big] \Big\} \mathrm{d}\xi \\ - \int_{\mathcal{J}} \Big\{ ({}_I\delta\boldsymbol{r}_C)_{,\xi}^\top \boldsymbol{A}_{IB}\, {}_B\boldsymbol{n} + ({}_B\delta\boldsymbol{\phi}_{IB})_{,\xi}^\top {}_B\boldsymbol{m} - {}_B\delta\boldsymbol{\phi}_{IB}^\top \big[ {}_B\bar{\boldsymbol{\gamma}} \times {}_B\boldsymbol{n} + {}_B\bar{\boldsymbol{\kappa}}_{IB} \times {}_B\boldsymbol{m} \big] \Big\} \mathrm{d}\xi \,. \end{aligned} \tag{11}$$

Note that even in the inextensible case, where no strain energy density $W$ with finite stiffness coefficients is available, the complementary strain energy density $W^\star$ still has finite coefficients. This can be seen as the limit of stiffness components approaching infinity, which results in zeros on the diagonals of $\boldsymbol{C}_{\boldsymbol{\gamma}}^{-1}$ and $\boldsymbol{C}_{\boldsymbol{\kappa}}^{-1}$. Shear-rigidity can therefore easily be introduced by using $\boldsymbol{C}_{\boldsymbol{\gamma}}^{-1} = \mathrm{diag}(k_{\mathrm{e}}^{-1}, 0, 0)$. Additionally, adding the inextensibility constraint leads to $\boldsymbol{C}_{\boldsymbol{\gamma}}^{-1} = \boldsymbol{0}_{3\times3}$. The same results were obtained by Eugster and Harsch (2020), where still the internal virtual work from (4) with the constitutive equations (3) was used for the unconstrained deformations.

Assume the line distributed external forces ${}_I\boldsymbol{b} = {}_I\boldsymbol{b}(\xi) \in \mathbb{R}^3$ and moments ${}_B\boldsymbol{c} = {}_B\boldsymbol{c}(\xi) \in \mathbb{R}^3$ to be given as densities with respect to the reference arc length. Moreover, for $i \in \{0, 1\}$, point forces ${}_I\boldsymbol{b}_i \in \mathbb{R}^3$ and point moments ${}_B\boldsymbol{c}_i \in \mathbb{R}^3$ can be applied to the rod's boundaries at $\xi_0 = 0$ and $\xi_1 = 1$. The corresponding external virtual work functional is defined as

$$\delta W^{\mathrm{ext}} := \int_{\mathcal{J}} \big\{ {}_I\delta\boldsymbol{r}_C^\top {}_I\boldsymbol{b} + {}_B\delta\boldsymbol{\phi}_{IB}^\top {}_B\boldsymbol{c} \big\} J \mathrm{d}\xi + \sum_{i=0}^{1} \big[ {}_I\delta\boldsymbol{r}_C^\top {}_I\boldsymbol{b}_i + {}_B\delta\boldsymbol{\phi}_{IB}^\top {}_B\boldsymbol{c}_i \big]_{\xi_i} \,. \tag{12}$$

From these virtual work contributions one can either derive the boundary value problem, see Appendix A, or derive finite element formulations, which is done in the following section.

# 3 Petrov–Galerkin finite element formulation

This section introduces the discretisation of the Cosserat rod. The kinematic fields ${}_I\boldsymbol{r}_{OC}$ and $\boldsymbol{A}_{IB}$, and the virtual displacement ${}_I\delta\boldsymbol{r}_C$ and rotation ${}_B\delta\boldsymbol{\phi}_{IB}$ are discretised with continuous Lagrange polynomials of degree $p$. The resultant contact force ${}_B\boldsymbol{n}$ and moment ${}_B\delta\boldsymbol{m}$ and their consistent variation are discretised with discontinuous Lagrange polynomials of degree $p-1$. Thus, they are discontinuous at the element boundaries. The discretisation results in discrete virtual work functionals from which the finite dimensional equations of the static equilibrium are derived.

## 3.1 Quaternion interpolation

For the discretisation, the rod's parameter space $\mathcal{J}$ is divided into $n_{\mathrm{el}}$ linearly spaced element intervals $\mathcal{J}^e = [\xi^e, \xi^{e+1})$ via $\mathcal{J} = \bigcup_{e=0}^{n_{\mathrm{el}}-1} \mathcal{J}^e$ with $\xi^e = e/n_{\mathrm{el}}$. To obtain a $p$-th order finite element, we introduce for each element $p+1$ linearly spaced nodes $\xi_i^e$, with $\xi_0^e = \xi^e$ and $\xi_p^e = \xi^{e+1}$, leading to $N = (pn_{\mathrm{el}} + 1)$ distinct nodes $\xi_k$ with $k = ep + i$. For the interpolation of the centerline curve ${}_I\boldsymbol{r}_{OC}$ and the cross-section orientations $\boldsymbol{A}_{IB}$, we introduce nodal centerline positions ${}_I\boldsymbol{r}_{OC_k} \in \mathbb{R}^3$ and nodal quaternions $\boldsymbol{P}_{IB_k} \in \mathbb{R}^4$ at these $N$ nodes. The nodal quaternions are related to the nodal transformation matrices by $\boldsymbol{A}_{IB_k} = \boldsymbol{A}(\boldsymbol{P}_{IB_k})$, where the mapping $\boldsymbol{A}(\boldsymbol{P}) \in SO(3)$ is defined for a quaternion $\boldsymbol{P} = (p_0, \boldsymbol{p})$ with the scalar part $p_0 \in \mathbb{R}$ and the vectorial part $\boldsymbol{p} \in \mathbb{R}^3$ by

$$\boldsymbol{A}(\boldsymbol{P}) = \boldsymbol{I} + 2\big(p_0 \tilde{\boldsymbol{p}} + \tilde{\boldsymbol{p}}^2\big)/\|\boldsymbol{P}\|^2 , \tag{13}$$

(see I. Romero 2004; Rucker 2018; Harsch and Eugster 2023)[1]. The generalised coordinates of element $e$ are collected in the tuple $\boldsymbol{q}^e \in \mathbb{R}^{7(p+1)}$, containing the nodal centerline points ${}_I\boldsymbol{r}_{OC_k}$ and nodal quaternions $\boldsymbol{P}_{IB_k}$ with indices $k \in \{pe, \dots pe + p\}$.

The generalised coordinates of the whole rod are collected in the tuple $\boldsymbol{q} \in \mathbb{R}^{7N}$, where the Boolean connectivity matrix $\boldsymbol{C}_{\boldsymbol{q},e} \in \mathbb{R}^{7(p+1)\times 7N}$ allows to extract the element tuple via $\boldsymbol{q}^e = \boldsymbol{C}_{\boldsymbol{q},e}\, \boldsymbol{q}$. Within each element, the interpolation of the centerline and the quaternion are given by

$${}_I\boldsymbol{r}^e_{OC}(\xi, \boldsymbol{q}^e) = \boldsymbol{N}^e_{\boldsymbol{r}}(\xi)\, \boldsymbol{q}^e , \quad \boldsymbol{P}^e_{IB}(\xi, \boldsymbol{q}^e) = \boldsymbol{N}^e_{\boldsymbol{P}}(\xi)\, \boldsymbol{q}^e , \tag{14}$$

where the matrices of Lagrange shape functions $\boldsymbol{N}^e_{\boldsymbol{r}} \in \mathbb{R}^{3\times 7(p+1)}$ and $\boldsymbol{N}^e_{\boldsymbol{P}} \in \mathbb{R}^{4\times 7(p+1)}$ with polynomial degree $p$ were used. A more detailed explanation of these matrices together with the ordering of the nodal quantities in the tuples of generalised coordinates can be found in Appendix B. For the sake of completeness, a global function for the centerline can be assembled by

$${}_I\boldsymbol{r}_{OC}(\xi, \boldsymbol{q}) = \sum_{e=0}^{n_{\mathrm{el}}-1} \chi_{\mathcal{J}^e}(\xi)\, {}_I\boldsymbol{r}^e_{OC}(\xi, \boldsymbol{C}_{\boldsymbol{q},e}\, \boldsymbol{q}) , \tag{15}$$

where $\chi$ is the characteristic function $\chi_{\mathcal{J}^e} : \mathcal{J} \to \{0, 1\}$, which is one for $\xi \in \mathcal{J}^e = [\xi^e, \xi^{e+1})$ and zero elsewhere. Similarly, a global function can also be assembled for the quaternion $\boldsymbol{P}_{IB}(\xi, \boldsymbol{q})$ based on the element functions $\boldsymbol{P}^e_{IB}(\xi, \boldsymbol{q}^e)$. The interpolated transformation matrix and the interpolated scaled curvature are then given by

$$\boldsymbol{A}_{IB} = \boldsymbol{A}(\boldsymbol{P}_{IB}) , \quad {}_B\bar{\boldsymbol{\kappa}}_{IB} = j^{-1}_{SO(3)}\big(\boldsymbol{A}^\top_{IB} \boldsymbol{A}_{IB,\xi}\big) = \boldsymbol{T}(\boldsymbol{P}_{IB})\, \boldsymbol{P}_{IB,\xi} , \tag{16}$$

where the tangent operator for the quaternion map $\boldsymbol{T}(\boldsymbol{P}) \in \mathbb{R}^{3\times 4}$ is

$$\boldsymbol{T}(\boldsymbol{P}) = \frac{2}{\|\boldsymbol{P}\|^2} \big( -\boldsymbol{p} \quad p_0 \boldsymbol{I} - \tilde{\boldsymbol{p}} \big) . \tag{17}$$

Its derivation is shown in Appendix C. This kinematic interpolation was already proposed by I. Romero (2004), who also provided a proof of objectivity. In contrast to the $\mathbb{R}^{12}$ interpolation introduced by Betsch and Steinmann (2002) and I. Romero and Armero (2002), the quaternion map in (13) ensures that $\boldsymbol{A}_{IB}$ is always an orthonormal transformation matrix.

At the same nodes $\xi_k$, we introduce the nodal virtual centerline displacement ${}_I\delta\boldsymbol{r}_{C_k} \in \mathbb{R}^3$ and the nodal virtual rotation ${}_B\delta\boldsymbol{\phi}_{IB_k} \in \mathbb{R}^3$. Similarly to the generalised coordinates, the generalised variations are collected element-wise in the tuple $\delta\boldsymbol{s}^e \in \mathbb{R}^{6(p+1)}$ and for the whole rod in $\delta\boldsymbol{s} \in \mathbb{R}^{6N}$. They are related by the Boolean connectivity matrix $\boldsymbol{C}_{\delta s,e} \in \mathbb{R}^{6(p+1)\times 6N}$ via $\delta\boldsymbol{s}^e = \boldsymbol{C}_{\delta s,e}\, \delta\boldsymbol{s}$. Within each element, the interpolation of the virtual displacement and the virtual rotation is given by

$${}_I\delta\boldsymbol{r}^e_C(\xi, \delta\boldsymbol{s}^e) = \boldsymbol{N}^e_{\delta\boldsymbol{r}}(\xi)\, \delta\boldsymbol{s}^e , \quad {}_B\delta\boldsymbol{\phi}^e_{IB}(\xi, \delta\boldsymbol{s}^e) = \boldsymbol{N}^e_{\delta\boldsymbol{\phi}}(\xi)\, \delta\boldsymbol{s}^e , \tag{18}$$

where the matrices of Lagrange shape functions $\boldsymbol{N}^e_{\delta\boldsymbol{r}} \in \mathbb{R}^{3\times 6(p+1)}$ and $\boldsymbol{N}^e_{\delta\boldsymbol{\phi}} \in \mathbb{R}^{3\times 6(p+1)}$ with polynomial degree $p$ were used. Similarly to (15), one can also build global functions based on these element functions.

For the interpolation of the resultant contact force and moment, we allow for discontinuities at the element boundaries. This is motivated due to the fact, that the strain measures are also discontinuous at the element boundaries after discretisation. Furthermore, we use $p-1$ as polynomial degree for the interpolation of the resultant contact forces and moments. This choice leads to $p$ nodes $\xi^e_{\sigma,i}$ per element for the interpolation of the resultant contact forces and moments.

[1] Note that this mapping also returns orthonormal matrices for non-unit quaternions, as the division by $\|\boldsymbol{P}\|^2$ can be seen as normalisation of the quaternion being relevant when interpolating the quaternions.

Considering the whole rod, there are $pn_{\mathrm{el}}$ nodes, where independent nodal resultant contact forces ${}_B\boldsymbol{n}_i^e \in \mathbb{R}^3$ and moments ${}_B\boldsymbol{m}_i^e \in \mathbb{R}^3$ are introduced. The nodal resultant contact forces and moments are collected element-wise in the tuple of resultant contact forces and moments $\boldsymbol{\lambda}_c^e \in \mathbb{R}^{6p}$ and for the whole rod in $\boldsymbol{\lambda}_c \in \mathbb{R}^{6pn_{\mathrm{el}}}$. They are related by the Boolean connectivity matrix $\boldsymbol{C}_{\boldsymbol{\lambda}_c,e} \in \mathbb{R}^{6p\times 6pn_{\mathrm{el}}}$ via $\boldsymbol{\lambda}_c^e = \boldsymbol{C}_{\boldsymbol{\lambda}_c,e}\,\boldsymbol{\lambda}_c$.

Within each element, the contact forces and moments are interpolated by

$$
{}_B\boldsymbol{n}^e(\xi,\boldsymbol{\lambda}_c^e) = \boldsymbol{N}_{\boldsymbol{n}}^e(\xi)\,\boldsymbol{\lambda}_c^e\,,\quad {}_B\boldsymbol{m}^e(\xi,\boldsymbol{\lambda}_c^e) = \boldsymbol{N}_{\boldsymbol{m}}^e(\xi)\,\boldsymbol{\lambda}_c^e\,, \tag{19}
$$

where the matrices of Lagrange shape functions $\boldsymbol{N}_{\boldsymbol{n}}^e \in \mathbb{R}^{3\times 6p}$ and $\boldsymbol{N}_{\boldsymbol{m}}^e \in \mathbb{R}^{3\times 6p}$ with polynomial degree $p-1$ were used. Similarly to (15), one can also build global functions based on these element functions. Since the variations of the resultant contact force and moment follow from the variation of (19), their interpolation uses the same matrices of shape functions $\boldsymbol{N}_{\boldsymbol{n}}^e$ and $\boldsymbol{N}_{\boldsymbol{m}}^e$ and Boolean connectivity matrices $\boldsymbol{C}_{\boldsymbol{\lambda}_c,e}$ together with element variations $\delta\boldsymbol{\lambda}_c^e \in \mathbb{R}^{6p}$ and global variations $\delta\boldsymbol{\lambda}_c \in \mathbb{R}^{6pn_{\mathrm{el}}}$.

An example on the assembly of the tuples $\boldsymbol{q}$, $\boldsymbol{q}^e$, $\boldsymbol{\lambda}_c$ and $\boldsymbol{\lambda}_c^e$ and the matrices $\boldsymbol{C}_{\boldsymbol{q},e}$, $\boldsymbol{N}_{\boldsymbol{r}}^e$, $\boldsymbol{N}_{\boldsymbol{P}}^e$, $\boldsymbol{C}_{\boldsymbol{\lambda}_c,e}$, $\boldsymbol{N}_{\boldsymbol{n}}^e$ and $\boldsymbol{N}_{\boldsymbol{m}}^e$ is given in Appendix B.2 for $n_{\mathrm{el}} = 3$ and $p = 2$.

## 3.2 Discrete virtual work functionals

With the introduced interpolation strategy for ansatz and test functions from (14), (18) and (19), we can now discretise the virtual work functionals. Since the integrals of the virtual work functional can be subdivided into the integrals of the element intervals due to the element-wise interpolation introduced in (15), the element contributions dependent only on the element functions. However, to increase readability, the superscript denoting the element number in the element functions will be suppressed. Additionally, we partly suppress the function arguments, which should be clear from the context.

To obtain the discretised form of the displacement-based virtual work, we insert the approximation of the virtual centerline displacement and the virtual cross-section rotation together with the corresponding approximations for cross-section orientations and strain measures into (4). The displacement-based internal virtual work is given by

$$
\begin{aligned}
&\delta W^{\mathrm{int},DB}(\boldsymbol{q};\delta\boldsymbol{s}) = \delta\boldsymbol{s}^\top \boldsymbol{f}^{\mathrm{int}}(\boldsymbol{q})\,,\quad \boldsymbol{f}^{\mathrm{int}}(\boldsymbol{q}) = \sum_{e=0}^{n_{\mathrm{el}}-1} \boldsymbol{C}_{\delta s,e}^\top \boldsymbol{f}_e^{\mathrm{int}}(\boldsymbol{C}_{\boldsymbol{q},e}\,\boldsymbol{q})\,,\\
&\boldsymbol{f}_e^{\mathrm{int}}(\boldsymbol{q}^e) = -\int_{\mathcal{J}^e} \Big\{\big(\boldsymbol{N}_{\delta\boldsymbol{r},\xi}^e\big)^\top \boldsymbol{A}_{IB}\,{}_B\boldsymbol{n} + \big(\boldsymbol{N}_{\delta\boldsymbol{\phi},\xi}^e\big)^\top {}_B\boldsymbol{m} - \big(\boldsymbol{N}_{\delta\boldsymbol{\phi}}^e\big)^\top \big({}_B\bar{\boldsymbol{\gamma}} \times {}_B\boldsymbol{n} + {}_B\bar{\boldsymbol{\kappa}}_{IB} \times {}_B\boldsymbol{m}\big)\Big\}\mathrm{d}\xi\,,
\end{aligned} \tag{20}
$$

where we have introduced the internal forces $\boldsymbol{f}^{\mathrm{int}}$ and their element contribution $\boldsymbol{f}_e^{\mathrm{int}}$.

To discretise the internal virtual work obtained by the Hellinger–Reissner principle, we insert the interpolation of the resultant contact forces and moments and their variations together with the corresponding approximations for virtual displacement, virtual rotation, cross-section orientations and strain measures into (11). The mixed internal virtual work is given by

$$
\begin{aligned}
\delta W^{\mathrm{int},MX}(\boldsymbol{q},\boldsymbol{\lambda}_c;\delta\boldsymbol{s},\delta\boldsymbol{\lambda}_c) &= \delta\boldsymbol{s}^\top \boldsymbol{W}_c(\boldsymbol{q})\boldsymbol{\lambda}_c + \delta\boldsymbol{\lambda}_c^\top\big(\boldsymbol{K}_c^{-1}\,\boldsymbol{\lambda}_c - \boldsymbol{l}_c(\boldsymbol{q})\big)\,,\\
\boldsymbol{W}_c(\boldsymbol{q}) &= \sum_{e=0}^{n_{\mathrm{el}}-1} \boldsymbol{C}_{\delta s,e}^\top \boldsymbol{W}_{c,e}(\boldsymbol{C}_{\boldsymbol{q},e}\,\boldsymbol{q})\boldsymbol{C}_{\boldsymbol{\lambda}_c,e}\,,\\
\boldsymbol{K}_c^{-1} &= \sum_{e=0}^{n_{\mathrm{el}}-1} \boldsymbol{C}_{\boldsymbol{\lambda}_c,e}^\top \boldsymbol{K}_{c,e}^{-1}\,\boldsymbol{C}_{\boldsymbol{\lambda}_c,e}\,,\\
\boldsymbol{l}_c(\boldsymbol{q}) &= \sum_{e=0}^{n_{\mathrm{el}}-1} \boldsymbol{C}_{\boldsymbol{\lambda}_c,e}^\top \boldsymbol{l}_{c,e}(\boldsymbol{C}_{\boldsymbol{q},e}\,\boldsymbol{q})\,,
\end{aligned} \tag{21}
$$

where we have introduced the matrix of compliance force directions $\boldsymbol{W}_c$, the constant compliance matrix $\boldsymbol{K}_c^{-1}$ and the compliance length $\boldsymbol{l}_c$, whose element contributions are given by

$$\begin{aligned}
\boldsymbol{W}_{c,e}(\boldsymbol{q}^e) &= -\int_{\mathcal{J}^e} \Big\{ \big(\mathbf{N}^e_{\delta\boldsymbol{r},\xi}\big)^\top \boldsymbol{A}_{IB}\, \mathbf{N}^e_{\boldsymbol{n}} + \big(\mathbf{N}^e_{\delta\boldsymbol{\phi},\xi}\big)^\top \mathbf{N}^e_{\boldsymbol{m}} - \big(\mathbf{N}^e_{\delta\boldsymbol{\phi}}\big)^\top \big({}_B\tilde{\bar{\boldsymbol{\gamma}}}\, \mathbf{N}^e_{\boldsymbol{n}} + {}_B\tilde{\bar{\boldsymbol{\kappa}}}_{IB}\, \mathbf{N}^e_{\boldsymbol{m}}\big) \Big\} \mathrm{d}\xi\,, \\
\boldsymbol{K}^{-1}_{c,e} &= \int_{\mathcal{J}^e} \Big\{ \big(\mathbf{N}^e_{\boldsymbol{n}}\big)^\top \boldsymbol{C}_{\boldsymbol{\gamma}}^{-1}\, \mathbf{N}^e_{\boldsymbol{n}} + \big(\mathbf{N}^e_{\boldsymbol{m}}\big)^\top \boldsymbol{C}_{\boldsymbol{\kappa}}^{-1}\, \mathbf{N}^e_{\boldsymbol{m}} \Big\} J \mathrm{d}\xi\,, \\
\boldsymbol{l}_{c,e}(\boldsymbol{q}^e) &= \int_{\mathcal{J}^e} \Big\{ \big(\mathbf{N}^e_{\boldsymbol{n}}\big)^\top \bar{\boldsymbol{\varepsilon}}_{\boldsymbol{\gamma}} + \big(\mathbf{N}^e_{\boldsymbol{m}}\big)^\top \bar{\boldsymbol{\varepsilon}}_{\boldsymbol{\kappa}} \Big\} \mathrm{d}\xi\,.
\end{aligned} \tag{22}$$

Similarly, the external virtual work is discretised by inserting approximations for virtual displacement, virtual rotation, into (12). The discretised external virtual work is given by

$$\begin{aligned}
\delta W^{\mathrm{ext}}(\boldsymbol{q};\delta\boldsymbol{s}) &= \delta\boldsymbol{s}^\top \boldsymbol{f}^{\mathrm{ext}}(\boldsymbol{q})\,, \\
\boldsymbol{f}^{\mathrm{ext}}(\boldsymbol{q}) &= \sum_{e=0}^{n_{\mathrm{el}}-1} \boldsymbol{C}^\top_{\delta s,e}\, \boldsymbol{f}^{\mathrm{ext}}_e(\boldsymbol{C}_{\boldsymbol{q},e}\boldsymbol{q}) + \sum_{i=0}^{1} \boldsymbol{C}^\top_{\delta s,e_i} \Big[ \big(\mathbf{N}^{e_i}_{\delta\boldsymbol{r}}\big)^\top {}_I\boldsymbol{b}_i + \big(\mathbf{N}^{e_i}_{\delta\boldsymbol{\phi}}\big)^\top {}_B\boldsymbol{c}_i \Big]_{\xi_i}\,, \\
\boldsymbol{f}^{\mathrm{ext}}_e(\boldsymbol{q}^e) &= \int_{\mathcal{J}^e} \Big\{ \big(\mathbf{N}^e_{\delta\boldsymbol{r}}\big)^\top {}_I\boldsymbol{b} + \big(\mathbf{N}^e_{\delta\boldsymbol{\phi},i}\big)^\top {}_B\boldsymbol{c} \Big\} J \mathrm{d}\xi\,,
\end{aligned} \tag{23}$$

where we have introduced the external forces $\boldsymbol{f}^{\mathrm{ext}}$ with their element contributions $\boldsymbol{f}^{\mathrm{ext}}_e$ and with $e_0 = 0$, $e_1 = n_{\mathrm{el}} - 1$, $\xi_0 = 0$ and $\xi_1 = 1$ for the point forces and moments at the rod's boundaries. The used discretisation of ansatz and test functions allows further to apply additional external point forces and moments at the internal element boundaries. Therefore, external loads can not only be applied at the two end points of the rod, but at all $n_{\mathrm{el}} + 1$ element boundaries $\xi^e$.

While the dependency of the internal virtual work on the generalised position coordinates $\boldsymbol{q}$ is clearly obvious through the occurrence of $\boldsymbol{A}_{IB}$, ${}_B\bar{\boldsymbol{\gamma}}$ and ${}_B\bar{\boldsymbol{\kappa}}_{IB}$ inside the integrals, the dependency of the external virtual work on the generalised position coordinates is only present, when the external force (${}_I\boldsymbol{b}$, ${}_I\boldsymbol{b}_0$ or ${}_I\boldsymbol{b}_1$) is not constant w.r.t. the inertial fixed basis $I$, or the external moment (${}_B\boldsymbol{c}$, ${}_B\boldsymbol{c}_0$ or ${}_B\boldsymbol{c}_1$) is not constant w.r.t. the cross-section-fixed basis $B$. This can happen, *e.g.*, when the external moment is constant w.r.t. the inertial fixed $I$-basis, which results in ${}_B\boldsymbol{c}_0 = \boldsymbol{A}^\top_{IB_0}\, {}_I\boldsymbol{c}_0$, as applied in Section 4.3.

Element integrals of the form $\int_{\mathcal{J}^e} f(\xi)\mathrm{d}\xi$ arising in the discretised internal forces as well as in the matrix of compliance force directions, in the compliance matrix and in the compliance length are subsequently computed using an $m$-point Gauss–Legendre quadrature rule. Because of the rational interpolation of $\boldsymbol{A}_{IB}$, the integrals in $\boldsymbol{f}^{\mathrm{int}}_e$, $\boldsymbol{W}_{c,e}$ and $\boldsymbol{l}_{c,e}$ cannot be computed exactly with this numerical quadrature rule. For polynomial degree $p = 1$, we refer to full integration when $m_{\mathrm{full}} = 2$ quadrature points are used. Similarly, for $p = 2$, full integration is obtained using $m_{\mathrm{full}} = 5$ quadrature points. As it will turn out in Section 4.1, this choice leads to the well-known locking behaviour for displacement-based formulations. Thus, we further introduce a reduced number of integration points $m_{\mathrm{red}} = p$ to evaluate the internal virtual work contributions. Further remarks on the reduced integration scheme are given in Section 3.4. For the evaluation of the integral in the external virtual work, we recommend to use $m_{\mathrm{ext}} \geqslant (p + p_{\mathrm{ext}} + 1)/2$ quadrature points, based on the polynomial degree $p_{\mathrm{ext}}$ of the distributed external force ${}_I\boldsymbol{b}$ and moment ${}_B\boldsymbol{c}$.

## 3.3 Equations of static equilibrium

The principle of virtual work states that the sum of all virtual work functional has to vanish for arbitrary virtual displacements (see dell'Isola and Steigmann 2020, Chapter 8), that is,

$$\begin{aligned}
DB: &\quad \delta W^{\mathrm{tot},DB} = \delta W^{\mathrm{int},DB} + \delta W^{\mathrm{ext}} \overset{!}{=} 0 \quad \forall \delta\boldsymbol{s} \quad \text{or} \\
MX: &\quad \delta W^{\mathrm{tot},MX} = \delta W^{\mathrm{int},MX} + \delta W^{\mathrm{ext}} \overset{!}{=} 0 \quad \forall \delta\boldsymbol{s}, \forall \delta\boldsymbol{\lambda}_c\,.
\end{aligned} \tag{24}$$

Together with the quaternion constraints to ensure a well-posed formulation

$$\boldsymbol{g}_S = (g_{S,0}, \dots g_{S,N-1}) = \boldsymbol{0}_N \quad \text{with} \quad g_{S,k} = \|\boldsymbol{P}_{IB_k}\|^2 - 1 = 0\,, \tag{25}$$

which constrains each nodal quaternion to unit length, the nonlinear generalised force equilibria

$$DB: \quad \begin{aligned} \boldsymbol{f}^{\mathrm{int}}(\boldsymbol{q}) + \boldsymbol{f}^{\mathrm{ext}}(\boldsymbol{q}) &= \mathbf{0}_{6N}\,, \\ \boldsymbol{g}_S(\boldsymbol{q}) &= \mathbf{0}_N\,, \end{aligned} \quad \text{or} \quad MX: \quad \begin{aligned} \boldsymbol{W}_c(\boldsymbol{q})\boldsymbol{\lambda}_c + \boldsymbol{f}^{\mathrm{ext}}(\boldsymbol{q}) &= \mathbf{0}_{6N}\,, \\ \boldsymbol{K}_c^{-1}\boldsymbol{\lambda}_c - \boldsymbol{l}_c(\boldsymbol{q}) &= \mathbf{0}_{6pn_{\mathrm{el}}}\,, \\ \boldsymbol{g}_S(\boldsymbol{q}) &= \mathbf{0}_N\,, \end{aligned} \tag{26}$$

are obtained. Note that the number of equations equals in both formulations the number of unknowns, since $\boldsymbol{q} \in \mathbb{R}^{7N}$ and $\boldsymbol{\lambda}_c \in \mathbb{R}^{6pn_{\mathrm{el}}}$ with the relation $N = pn_{\mathrm{el}} + 1$. Following Géradin and Cardona (2001), prescribed boundary conditions can be incorporated into the principle of virtual work using perfect bilateral constraints. In both formulations, the nonlinear equilibrium equations can be represented in the form $\boldsymbol{f}(\boldsymbol{x}) = \mathbf{0}_n$, where $\boldsymbol{x} \in \mathbb{R}^n$ contains all unknowns, *i.e.*, the generalised coordinates, possibly Lagrange multipliers for bilateral constraints and in the mixed formulation also the generalised forces. Such a nonlinear equation is solved by any root-finding algorithm, *e.g.*, Newton–Raphson, Riks. Note that a system of linear equations with a non-symmetric iteration matrix $\partial \boldsymbol{f}/\partial \boldsymbol{x} \in \mathbb{R}^{n\times n}$ must be solved in each iteration. In this work, the static equilibrium will be computed using Newton's method. As termination criteria we require the Euclidean norm $\sqrt{\boldsymbol{f}(\boldsymbol{x})^\top \boldsymbol{f}(\boldsymbol{x})}$ to be smaller than $\varepsilon\sqrt{n}$, where $\varepsilon$ is the tolerance.

Although we derived only the equations for static equilibrium in this work, for dynamics the equations of motion can also be derived. Since the same generalised coordinates are chosen as in Harsch and Eugster (2023) and also the interpolation of the variations is identical, the inertial virtual work and the kinematic differential equation can be derived in the same way, leading to equations of motion.

## 3.4 Constraint deformations and static condensation

Following the discussion after (11), the displacement-based formulation cannot numerically capture the limiting case in which stiffness parameters tend to infinity. In contrast, the mixed formulation covers this by setting $k_i = 0$. With that the matrix $\boldsymbol{K}_c^{-1}$ is not invertible anymore. To be more specific, let $n_g$ be the number of constrained deformations, e.g., $n_g = 2$ for shear-stiff rods, and $n_g = 3$ for shear-stiff and inextensible rods, then $\boldsymbol{K}_c^{-1}$ has a rank-deficiency of $n_g pn_{\mathrm{el}}$. However, it is easily possible to split the compliance equation into two parts by

$$\bar{\boldsymbol{K}}_c^{-1}\bar{\boldsymbol{\lambda}}_c - \bar{\boldsymbol{l}}_c(\boldsymbol{q}) = \mathbf{0}_{(6-n_{\mathrm{g}})pn_{\mathrm{el}}}\,, \quad \boldsymbol{g}_c(\boldsymbol{q}) = \mathbf{0}_{n_{\mathrm{g}}pn_{\mathrm{el}}}\,, \tag{27}$$

containing the same equations as $\boldsymbol{K}_c^{-1}\boldsymbol{\lambda}_c - \boldsymbol{l}_c(\boldsymbol{q}) = \mathbf{0}_{6pn_{\mathrm{el}}}$ but with the matrix $\bar{\boldsymbol{K}}_c^{-1}$ having full rank and thus being invertible. The nodal resultant contact forces and moments of $\boldsymbol{\lambda}_c$ that are not contained in $\bar{\boldsymbol{\lambda}}_c$ are then the Lagrange multipliers required to satisfy $\boldsymbol{g}_c(\boldsymbol{q}) = \mathbf{0}_{n_{\mathrm{g}}pn_{\mathrm{el}}}$. Thus, the constraints derived through this method are equivalent to constrained rod finite element models based on the Lagrange multiplier methods, see for example Harsch et al. (2021). Furthermore, we can also use only the constraints from the mixed formulation and keep the unconstrained deformations in the displacement-based formulation. This approach will be later used in the experiment presented in Section 4.6, where we show the case of a non-convex strain energy density function together with constraints.

For simplicity, assume $n_{\mathrm{g}} = 0$, such that $\boldsymbol{K}_c^{-1}$ has full rank and is invertible. That allows to obtain $\boldsymbol{\lambda}_c = \boldsymbol{K}_c\,\boldsymbol{l}_c(\boldsymbol{q})$, which can then be inserted into the force equation of the mixed generalised force equilibrium (26), resulting in the static condensed generalised force equilibrium

$$SC: \quad \begin{aligned} \boldsymbol{W}_c(\boldsymbol{q})\boldsymbol{K}_c\,\boldsymbol{l}_c(\boldsymbol{q}) + \boldsymbol{f}^{\mathrm{ext}}(\boldsymbol{q}) &= \mathbf{0}_{6N}\,, \\ \boldsymbol{g}_S(\boldsymbol{q}) &= \mathbf{0}_N\,. \end{aligned} \tag{28}$$

The strategy of static condensation can also be used for $n_{\mathrm{g}} \neq 0$ for the invertible part $\bar{\boldsymbol{K}}_c^{-1}$ of $\boldsymbol{K}_c^{-1}$, corresponding to $\bar{\boldsymbol{\lambda}}_c$ and $\bar{\boldsymbol{l}}_c$.

Following the argumentation of Noor and Peters (1981); Malkus and Hughes (1978), the solution of the generalised coordinates $\boldsymbol{q}$ using the displacement-based formulation is equivalent to the solution of the generalised coordinates $\boldsymbol{q}$ using the mixed formulation, when in both cases reduced integration is applied. Since the resultant contact forces and moments are interpolated with polynomial degree $p-1$, and therefore having $p$ unknowns per component, this matches the number of quadrature points used in case of reduced integration. In the mixed formulation, the resultant contact forces and moments therefore fulfill at these quadrature points the constitutive equations ${}_B\boldsymbol{n} = \boldsymbol{C}_\gamma\, \boldsymbol{\varepsilon}_\gamma$ and ${}_B\boldsymbol{m} = \boldsymbol{C}_\kappa\, \boldsymbol{\varepsilon}_\kappa$. In the displacement-based formulation, the constitutive equations are fulfilled at all points and especially also at the quadrature points. Consequently, also the static condensed formulation has the same solution $\boldsymbol{q}$, and it must hold therefore $\boldsymbol{f}^{\text{int}}(\boldsymbol{q}) = \boldsymbol{W}_c(\boldsymbol{q})\boldsymbol{K}_c\, \boldsymbol{l}_c(\boldsymbol{q})$ in case of reduced integration.

We want to highlight that, although the solution of the generalised coordinates $\boldsymbol{q}$ of the static condensed formulation and the mixed formulation are identical, the iterates taken by the solver are different. This is also in agreement with (Albersmeyer and Diehl 2010), where the mixed formulation naturally introduces the auxiliary variables leading to a much better convergence when the Newton–Raphson method is applied. Note also that, due to numerical computations, the iterations and solutions of the different formulations differ slightly.

### 3.5 *SE*3 interpolation

Additionally, the *SE*3 interpolation outlined in Harsch et al. (2023); Eugster and Harsch (2023) will be analysed. The formulation is adapted to also use a quaternion parametrisation for the nodal transformation matrices instead of rotation vectors. Therefore, the quaternion constraints $\boldsymbol{g}_S(\boldsymbol{q})$ are also introduced. Since the *SE*3 interpolation leads to element-wise constant strains, a mixed internal virtual work with element-wise constant resultant contact forces and moments is applied. The tuples $\boldsymbol{q}$, $\delta\boldsymbol{s}$, $\boldsymbol{\lambda}_c$ and $\delta\boldsymbol{\lambda}_c$ are therefore identical to the ones used in the quaternion interpolation with polynomial degree $p = 1$.

## 4 Numerical experiments

In this chapter several benchmark experiments will be shown. These experiments allow us to compare the different formulations and show their advantages. In total, we will compare 12 different formulations which are combinations of different interpolations, the different internal virtual work formulations and the different number of spatial integration points. The formulations and their shorthand abbreviations are given in Table 1. For the quantitative analyses of the formulations, we consider the spatial convergence and the convergence rate of the solver iterations. The detailed formulae for computing the spatial error and the convergence rate are provided in Appendix D. When solving the equations of static equilibrium, we always solve the full equations (26). The static condensed equations are only used when explicitly stated. For these, we use an abbreviation following the logic of Table 1.

### 4.1 45° bent experiment

The first benchmark experiment was used by several authors (see for instance Bathe and Bolourchi 1979; Lo 1992; Ibrahimbegović 1995; Jelenić and M. Crisfield 1999; Ghosh and Roy 2008; Bauer et al. 2016; C. A. Meier 2016; Greco et al. 2024). In this experiment, an initially precurved rod with quadratic cross-section of width $w$ is clamped at $\xi = 0$ such that $\boldsymbol{A}_{IB}(0) = \boldsymbol{I}$ and ${}_I\boldsymbol{r}_{OC}(0) = \boldsymbol{0}$. The rod's centerline in the undeformed configuration represents $1/8$ of a circle (45°) with radius $R = 100$, as shown in Figure 1(a). The strategy to obtain the initial set of generalised coordinates $\boldsymbol{q}$ is explained in Appendix B.3. The stiffnesses $k_\mathrm{e} = EA$, $k_{\mathrm{s}_y} = k_{\mathrm{s}_z} = GA$, $k_\mathrm{t} = 2GI$ and $k_{\mathrm{b}_y} = k_{\mathrm{b}_z} = EI$ are given in terms of Young's and shear moduli $E = 10^7$ and $G = E/2$, the cross-section area $A = w^2$ and the second moment of area $I = w^4/12$. An external force ${}_I\boldsymbol{b}_1 = (0, 0, F_z)$ acts at $\xi = 1$. The slenderness ratio is defined in terms of the radius $R$ and the width $w$. We considered the ratios $\rho = R/w \in \{10^1, 10^2, 10^3, 10^4\}$, resulting in $w \in \{10, 1, 0.1, 0.01\}$. The reference solution for the analysis of the spatial convergence behavior was obtained by using the $Q^2_{MX_\mathrm{full}}$ formulation with 256 elements, *i.e.*, $N = 513$ nodes. The solver tolerance $\varepsilon$ ($\varepsilon_{DB}$ for displacement-based formulations and $\varepsilon_{MX}$ for mixed formulations) and the force component $F_z$ for each slenderness ratio are given in Table 2. The external force was applied within 50 load increments, except for the reference rod, where the load was applied within 5 load increments. The deformed configurations for the slenderness ratio $\rho = 10^2$ of the reference rod are shown in Figure 1(b). The vertical force component of the shown deformed configurations is $F_z \in \{0, 120, 240, 360, 480, 600\}$. The components of the tip displacements for these loads are shown in Figure 1(c).

| | Internal virtual work / Spatial integration | displacement-based full | displacement-based reduced | mixed full | mixed reduced |
|---|---|---|---|---|---|
| Interpolation | Quaternion, $p = 1$ | $Q^1_{DB_\text{full}}$ | $Q^1_{DB_\text{red}}$ | $Q^1_{MX_\text{full}}$ | $Q^1_{MX_\text{red}}$ |
| | Quaternion, $p = 2$ | $Q^2_{DB_\text{full}}$ | $Q^2_{DB_\text{red}}$ | $Q^2_{MX_\text{full}}$ | $Q^2_{MX_\text{red}}$ |
| | *SE*3 | $SE3_{DB_\text{full}}$ | $SE3_{DB_\text{red}}$ | $SE3_{MX_\text{full}}$ | $SE3_{MX_\text{red}}$ |

**Table 1** Different formulations with their shorthand abbreviations used for the numerical analyses as combinations of interpolation, internal virtual work and spatial integration.

| $\rho$ | $10^1$ | $10^2$ | $10^3$ | $10^4$ |
|---|---|---|---|---|
| $F_z$ | $6 \cdot 10^6$ | $6 \cdot 10^2$ | $6 \cdot 10^{-2}$ | $6 \cdot 10^{-6}$ |
| $\varepsilon_{DB}$ | $10^{-2}$ | $10^{-6}$ | $10^{-8}$ | $10^{-10}$ |
| $\varepsilon_{MX}$ | $10^{-2}$ | $10^{-6}$ | $10^{-10}$ | $10^{-13}$ |

**Table 2** 45° bent experiment: Parameters for the different slendernesses.

Figure 2 shows the spatial convergence behaviour for the different formulations and different number of integration points as well as the different slenderness ratios. From top to bottom, the rows correspond to slenderness ratios of $10^1$, $10^2$, $10^3$, and $10^4$. While the first two columns show

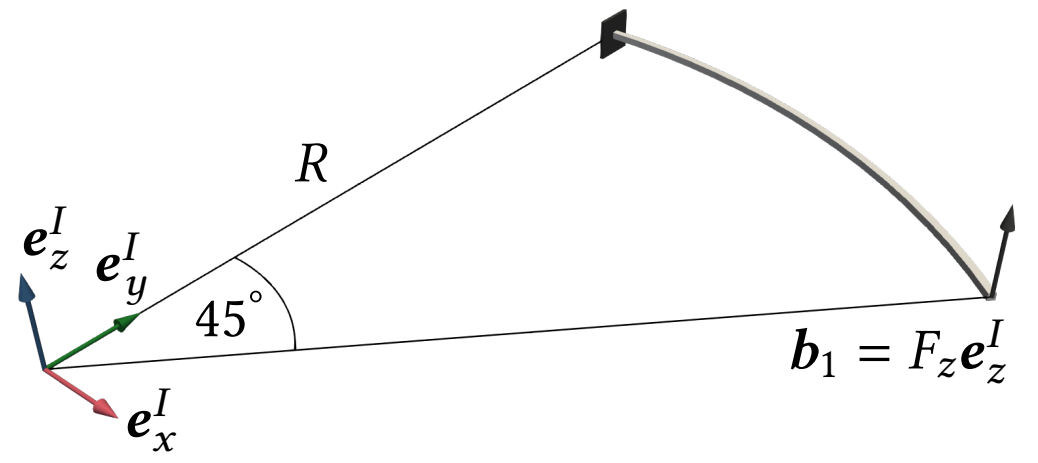

**(a)** Setup of the experiment

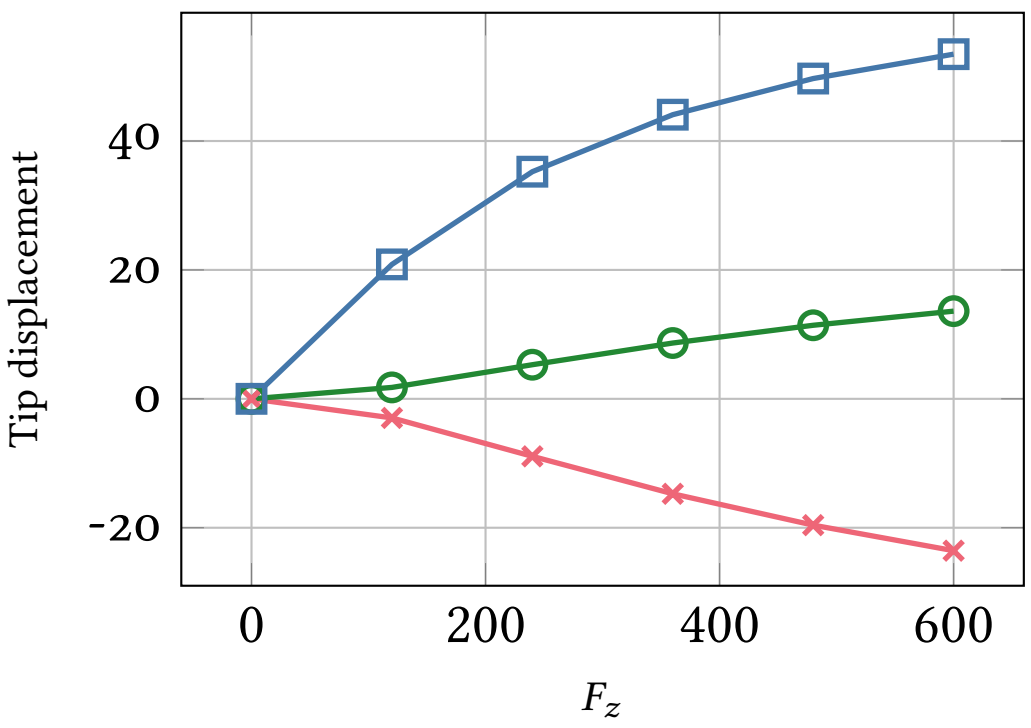

**(b)** Tip displacement

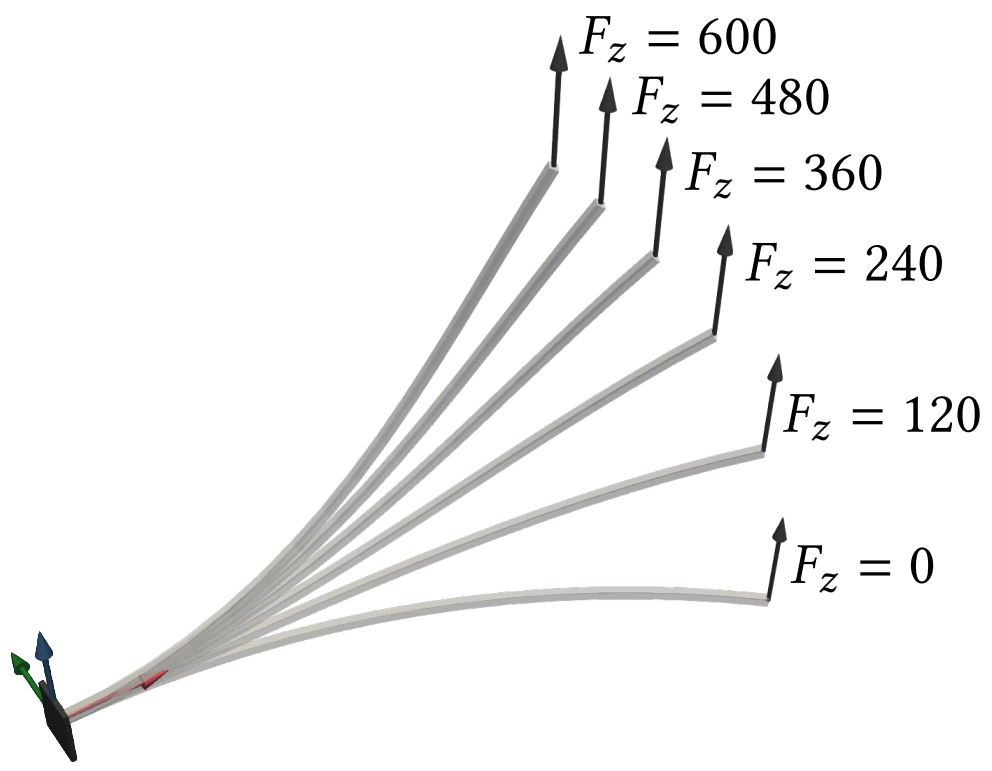

**(c)** Deformed configurations

**Figure 1** 45° bent experiment: Overview of the experiment. The shown configurations and tip displacements were obtained by using the reference rod ($Q^2_{MX_\mathrm{full}}$ formulation with $n_\mathrm{el} = 256$ elements) and 5 load increments for the slenderness ratio $\rho = 10^2$. The components of the tip displacement $(\Delta r_x, \Delta r_y, \Delta r_z) = {}_I\boldsymbol{r}_{OC}(1) - {}_I\boldsymbol{r}^0_{OC}(1)$ are indicated by (——), (——) and (——), respectively.

displacement-based formulations (full integration in the first columns and reduced integration in the second column), the latter two columns show mixed formulations (full integration in the third column and reduced integration in the fourth column). The spatial convergence reveals, that locking occurs in the displacement-based formulations of the quaternion interpolation when full integration is used. Increasing the polynomial degree can further alleviate locking, see the difference between $\boldsymbol{Q}^1$ and $\boldsymbol{Q}^2$ for the high slenderness ratios, but locking is still present. Applying the reduced integration method removes locking from the displacement-based formulations. Furthermore, there occurs no locking when the mixed formulation is used, no matter which integration method is applied. Since the $SE3$ interpolation is intrinsically locking free, the spatial convergence rate is for all four cases identical.

Figure 3 shows the resultant contact forces and moments for $\boldsymbol{Q}^1$, $\boldsymbol{Q}^2$ and $SE3$ formulations, respectively. In each figure, the top row shows the resultant contact forces ${}_B\boldsymbol{n}$ and the bottom row shows the resultant contact moments ${}_B\boldsymbol{m}$. The resultant contact forces and moments are obtained for the slenderness ratio $\rho = 10$, where each rod was discretised using $N = 9$ nodes. The displacement-based formulations of the quaternion interpolation lead to a high fluctuation of the resultant contact forces with similar results for the different integration methods. Note that the shown experiment with the slenderness ratio $\rho = 10$ shows no locking behaviour in the spatial convergence. Using a higher polynomial reduces the magnitude of fluctuations, but it is still large in comparison to the results of the mixed formulations. The resultant contact forces and moments of the $SE3$ interpolation are identical, no matter which formulation and integration method is used.

## 4.2 Helix experiment

This benchmark of the helix experiment was used in C. Meier et al. (2015), Harsch et al. (2021) and Harsch and Eugster (2023). In this work, it is used to demonstrate the computational efficiency of the mixed approach. A helical curve with $n = 2$ coils along the $\boldsymbol{e}_z^I$ axis with height $h = 50$ and radius $R_0 = 10$ is described with $\xi \in [0, 1]$ by

$$ {}_I\boldsymbol{r}_{OC}^*(\xi) = R_0\Big(\sin\big(\alpha(\xi)\big),\ -\cos\big(\alpha(\xi)\big),\ c\alpha(\xi)\Big), \quad \alpha(\xi) = 2\pi n\xi\,. \tag{29} $$

where $c = h/(2\pi R_0 n)$ is the pitch of the helix. The length of the curve is $L = 2\pi R_0 n\sqrt{1+c^2}$. The considered rod is initially straight with length $L$ and has a circular cross-section with radius $r = L/(2\rho)$, where $\rho \in \{10^1, 10^2, 10^3, 10^4\}$ is the slenderness ratio. The stiffnesses are given by $k_\mathrm{e} = EA$, $k_{\mathrm{s}_y} = k_{\mathrm{s}_z} = GA$, $k_\mathrm{t} = 2GI$ and $k_{\mathrm{b}_y} = k_{\mathrm{b}_z} = EI$ in terms of Young's and shear moduli $E = 1$ and $G = 1/2$, the cross-section area $A = \pi r^2$ and the second moment of area $I = \pi r^4/4$. Note that by this choice, it holds $k_\mathrm{t} = k_{\mathrm{b}_y} = k_{\mathrm{b}_z}$. The rod is clamped at $\xi = 0$, such that ${}_I\boldsymbol{r}_{OC}(0) = (0, -R_0, 0)$, ${}_I\boldsymbol{e}_x^B(0) = (1, 0, c)/\sqrt{1+c^2}$, ${}_I\boldsymbol{e}_y^B(0) = (0, 1, 0)$ and ${}_I\boldsymbol{e}_z^B(0) = {}_I\boldsymbol{e}_x^B(0) \times {}_I\boldsymbol{e}_y^B(0)$. To bend the rod into the shape from (29), a moment ${}_B\boldsymbol{c}_1 = (ck_\mathrm{t}, 0, k_{\mathrm{b}_z})/\big(R_0(1+c^2)\big)$ is applied at the free end at $\xi = 1$. The analytical solution has constant resultant contact forces ${}_B\boldsymbol{n}(\xi) = \boldsymbol{0}$ and moments ${}_B\boldsymbol{m}(\xi) = {}_B\boldsymbol{c}_1$, as discussed in Harsch et al. (2021). All helix experiments were performed using $N = 17$ nodes. In Figure 4, the final deformed configuration of the slenderness ratios $\rho \in \{10^1, 10^2\}$ are shown.

Table 3 shows the tolerances used for Newton's solver and the minimum number of load increments to get a converged solution. The number of load increments is determined by testing powers of two, starting with $2^0 = 1$ load increments, followed by $2^1 = 2$, $2^2 = 4$, $2^3 = 8$ and so on. Note that the table only shows $DB_\mathrm{red}$, $MX_\mathrm{red}$ and $MX_\mathrm{full}$ formulations, which are the only formulations that are locking free for all interpolations. The number of required load increments for the displacement-based formulations is significantly larger than the required number for the mixed formulations, which is almost always one.

Figure 5 shows the resultant contact forces and moments for the slenderness ratio $\rho = 10^1$, obtained with 128 load increments for displacement-based formulations and with 1 load increment for mixed formulations. The resultant contact forces and moments of both quaternion interpolations with $DB_\mathrm{full}$ indicate locking, which can be seen as the resultant contact moment is not "periodic" with respect to the centerline parameter. In fact, the deformed rod is not even close to the desired helical shape. The resultant contact forces and moments of both quaternion interpolations with $DB_{\text{red}}$ show the typical fluctuation around the exact value, as already noted in the previous experiment. The resultant contact forces and moments for $MX_{\text{full}}$ and $MX_{\text{red}}$ are almost identical and represent in each case the expected solution. All formulations of the $SE3$ rod are not only able to represent the correct resultant contact forces and moments, but also to correctly represent the helix shape, as discussed in Harsch et al. (2021).

In Figure 6(a), the average and the standard deviation of the required number of Newton iterations per load increment are shown. It is remarkable that all mixed formulations require a

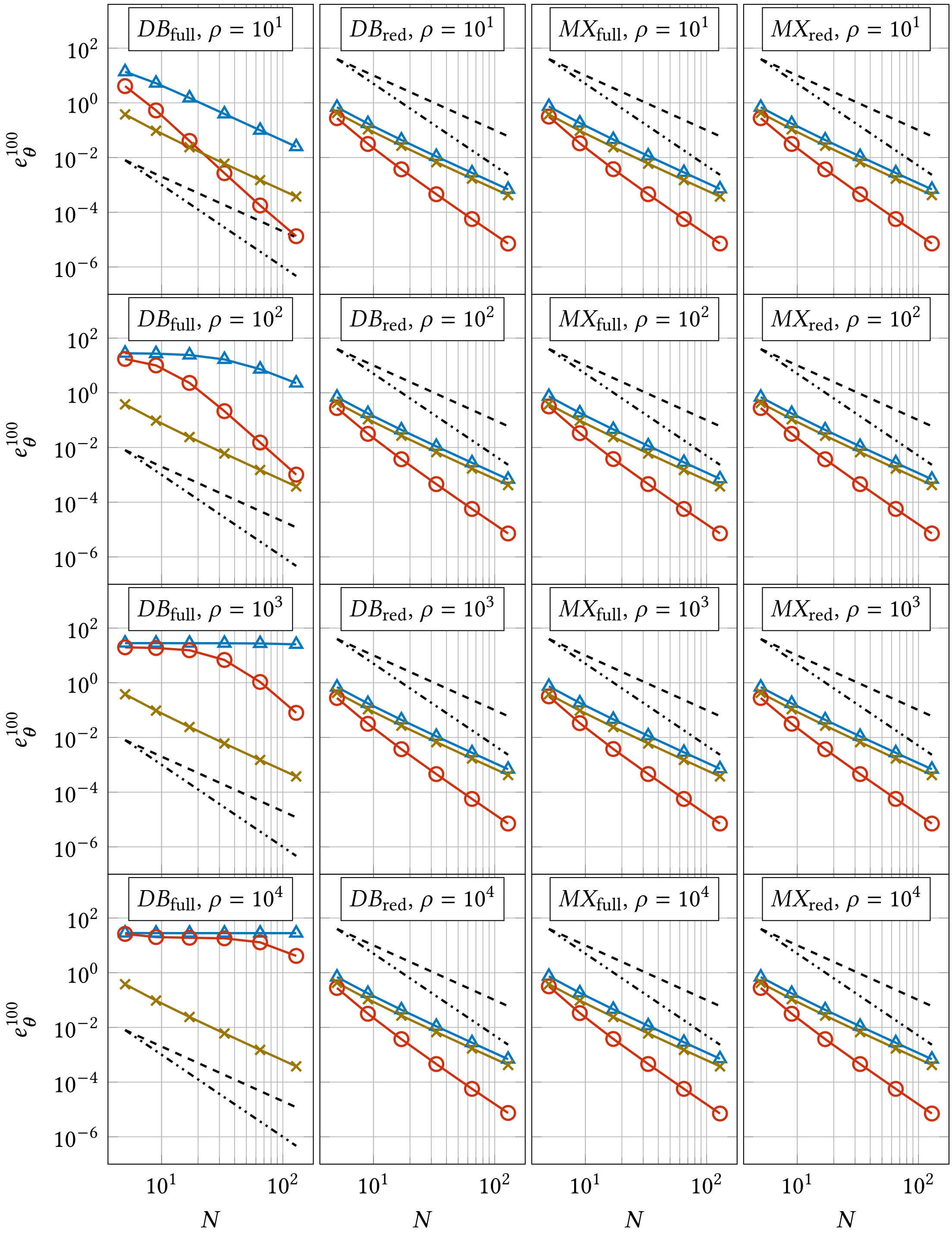

**Figure 2** 45° bent experiment: Spatial convergence rates for the different slenderness ratios and different formulations of the internal virtual work. For the rod, discretised with polynomial degree $p$ and $n_{\text{el}}$ elements, the number of nodes is $N = (pn_{\text{el}} + 1)$. The used error measure $e_\theta^{100}$ is a combination of position and orientation error at 100 points along the rod. The shown rates are for $Q^1$ (——), $Q^2$ (——) and $SE3$ (——) interpolations. The additional lines are proportional to $N^{-2}$ (- - -) and $N^{-3}$ (-···-).

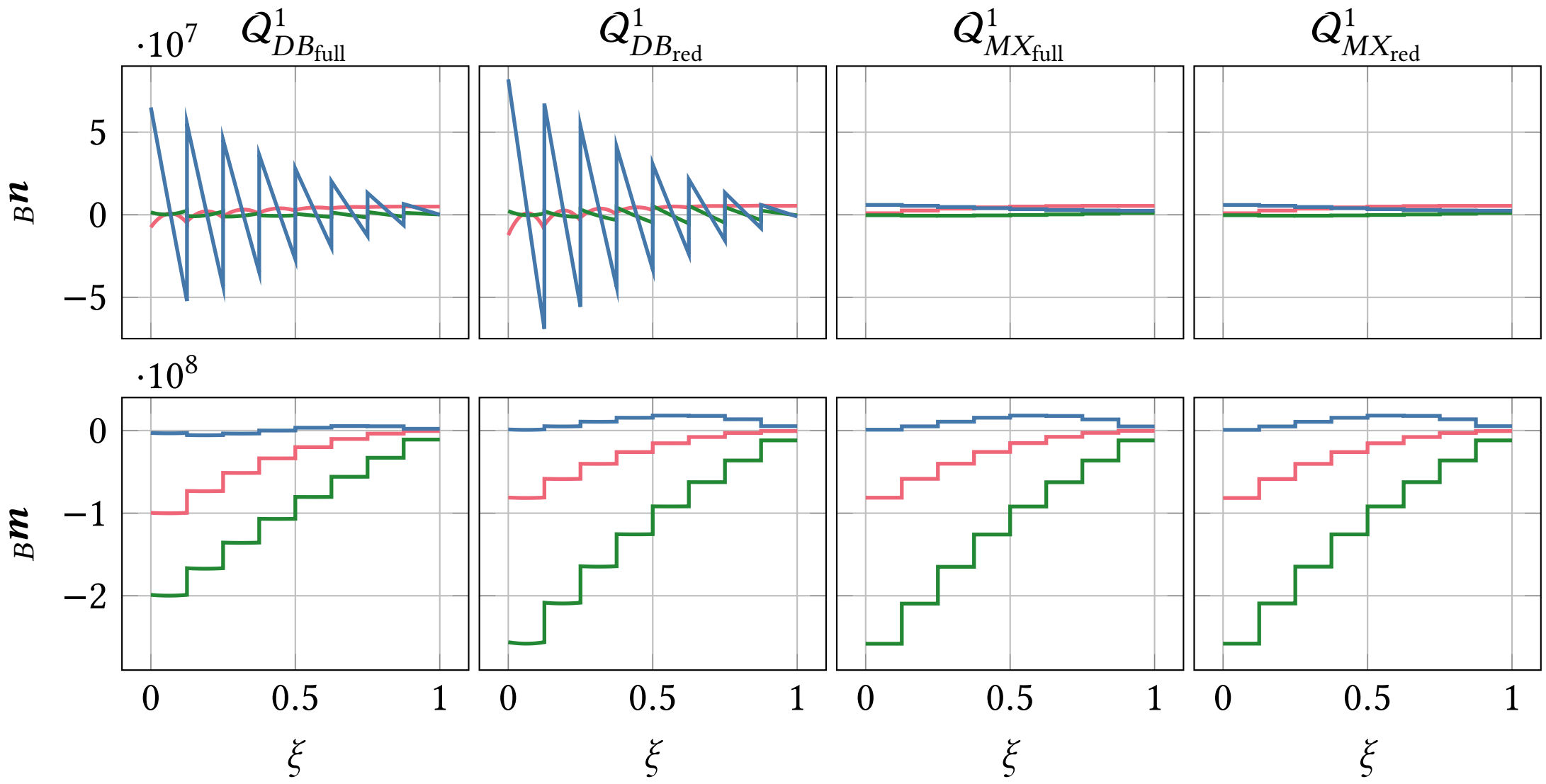

**(a)** Quaternion interpolation with polynomial degree $p = 1$ and $n_{\mathrm{el}} = 8$ elements ($N = 9$ nodes)

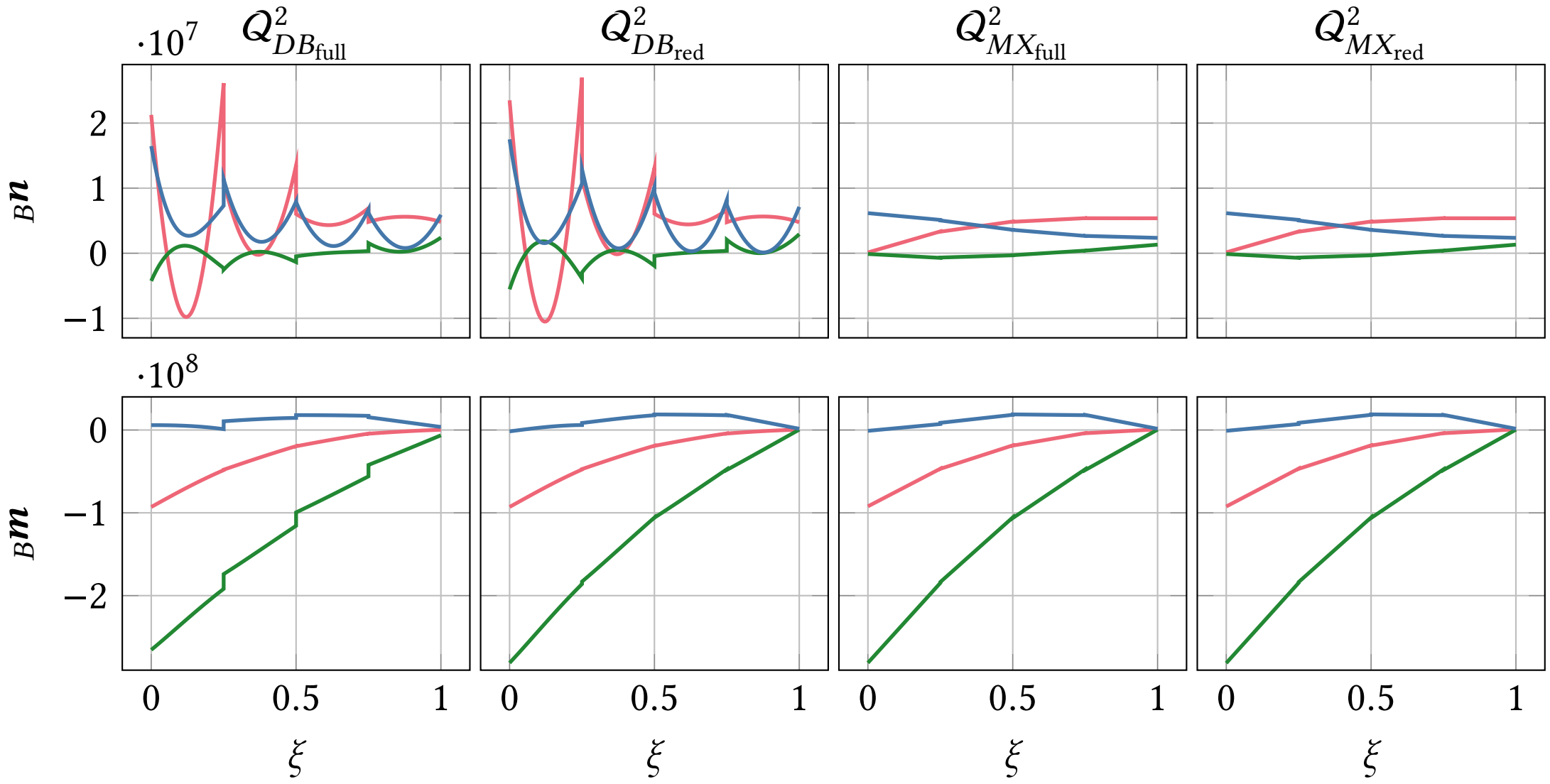

**(b)** Quaternion interpolation with polynomial degree $p = 2$ and $n_{\mathrm{el}} = 4$ elements ($N = 9$ nodes)

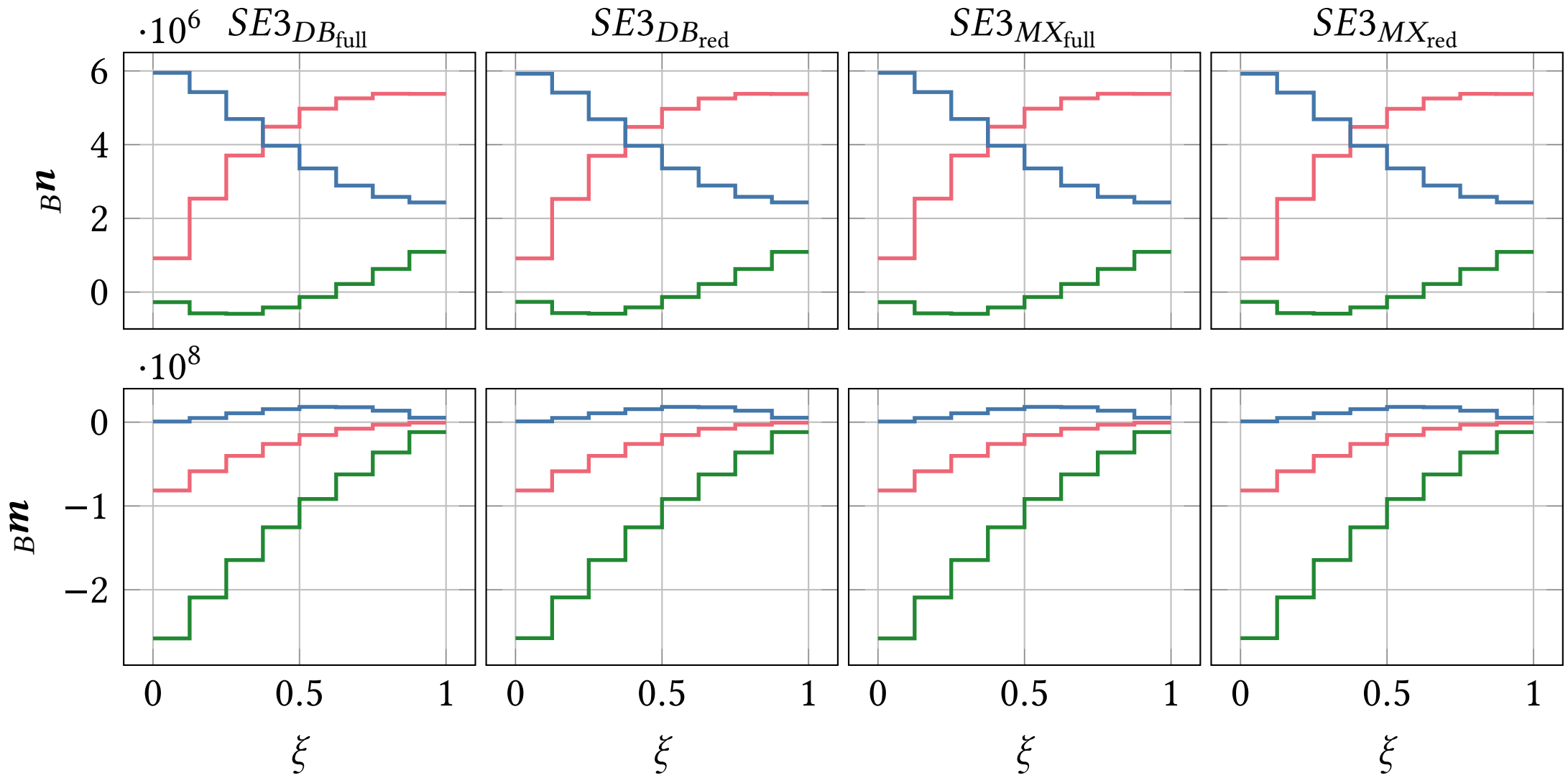

**(c)** $SE3$ interpolation with $n_{\mathrm{el}} = 8$ elements ($N = 9$ nodes)

**Figure 3** 45° bent experiment: Components of the resultant contact forces $_B\boldsymbol{n}$ (top row) and moments $_B\boldsymbol{m}$ (bottom row) of the different formulations for the slenderness $\rho = 10^1$. The components in $\boldsymbol{e}_x^B$, $\boldsymbol{e}_y^B$ and $\boldsymbol{e}_z^B$ direction are shown as (——), (——) and (——), respectively. Note, the scaling of the resultant contact forces differs between (a), (b) and (c).

much smaller number of Newton iterations than the displacement-based formulations. Also, the standard deviation of the mixed formulations is zero, *i.e.*, each load increment requires the same number of Newton iterations, while the displacement-based formulations show a large standard deviation. Together with the minimum number of load increments from Table 3, it is clearly evident how the computational efficiency is significantly increased by applying the mixed formulation. Note also that the computational efficiency increases when the mixed formulation is applied to the $SE3$-formulation. Although this formulation is intrinsically locking free with element-wise constant strains, the need for the application of the mixed method is still highly justified due to the computational improvements. Furthermore, in Figure 6(b), the average and the standard deviation of the local quadratic convergence rates, as computed in (D.3), are shown. The rates of the mixed formulations are much smaller than the rates of the displacement-based formulations, and are completely unaffected by the slenderness of the rod. This indicates not only better convergence, but also a computationally more robust formulation.

## 4.3 Rod bent to a helical form

The next benchmark was proposed the first time in Ibrahimbegović (1997) to simulate large inhomogeneous deformations. Here it is used to validate the shown computational improvements from the previous section on another problem. In this experiment, an initially straight cantilever with length $L = 10$ is considered. The stiffnesses are given by $k_\mathrm{e} = k_{\mathrm{s}_y} = k_{\mathrm{s}_z} = 10^4$ and $k_\mathrm{t} = k_{\mathrm{b}_y} = k_{\mathrm{b}_z} = 10^2$. On the cantilever's free tip an external moment ${}_B\boldsymbol{c}_1 = \boldsymbol{A}_{IB}^\top(0, 0, 20\pi k_{\mathrm{b}_y}/L)$ and an external force ${}_I\boldsymbol{b}_1 = (0, 0, 50)$ are applied, see Figure 7(a). Note that both the moment and the force are constant w.r.t. the inertial $I$-system, resulting in a generalised external force which depends on the generalised coordinates $\boldsymbol{q}$. As the stiffness parameters are not computed by geometric properties, the shown circular cross-section is only for visualisation. Further, note that without the perturbation force ${}_I\boldsymbol{b}_1$, the external moment would roll-up the cantilever to 10 coils, where the resultant contact forces and moments would be constant.

The tolerance of Newton's solver was set to $\epsilon = 10^{-8}$. We used for every rod $N = 61$ nodes, leading to $n_\mathrm{el} = 60$ elements for the $\boldsymbol{Q}^1$ and $SE3$ formulations and $n_\mathrm{el} = 30$ elements for the $\boldsymbol{Q}^2$

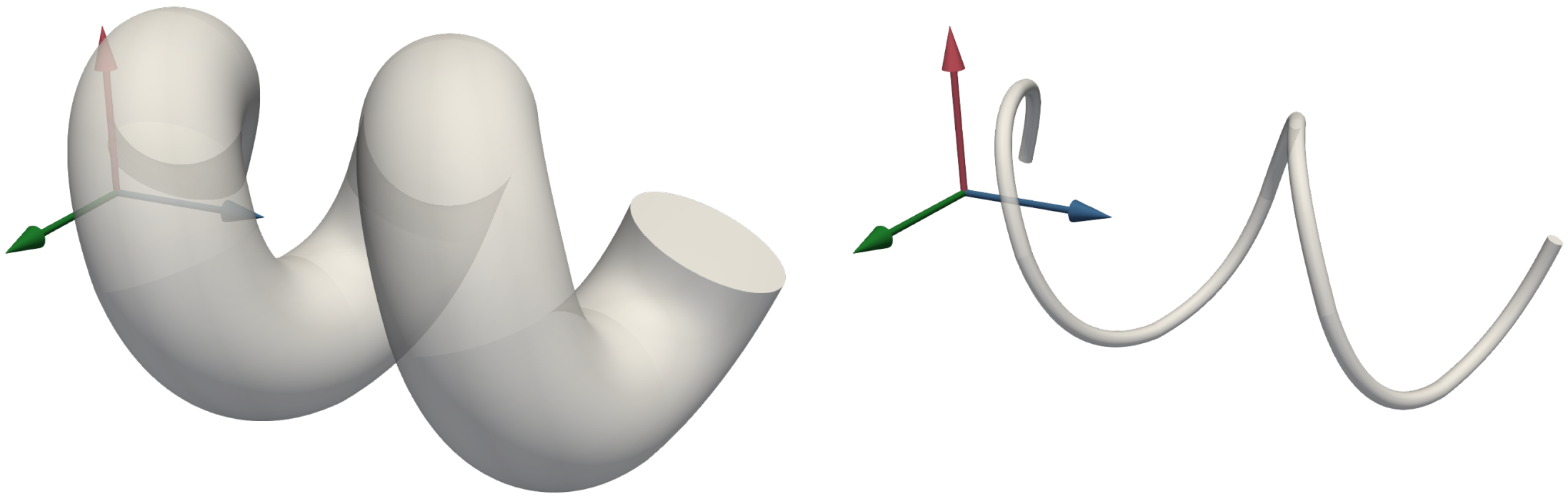

**Figure 4** Helix experiment: Final deformed configurations for the slendernesses $\rho \in \{10^1, 10^2\}$. The deformed configurations were obtained using the $\boldsymbol{Q}^2_{MX_\mathrm{full}}$ interpolation with $n_\mathrm{el} = 8$ elements, *i.e.*, $N = 17$ nodes. The red, green and blue arrows indicate the $\boldsymbol{e}_x^I$, $\boldsymbol{e}_y^I$ and $\boldsymbol{e}_z^I$ axes, respectively.

| | | $DB_\mathrm{red}$ / $SC_\mathrm{red}$ | | | $MX_\mathrm{red}$ | | | $MX_\mathrm{full}$ | | |
|---|---|---|---|---|---|---|---|---|---|---|
| Slenderness | Tolerance | $Q^1$ | $Q^2$ | $SE3$ | $Q^1$ | $Q^2$ | $SE3$ | $Q^1$ | $Q^2$ | $SE3$ |
| $10^1$ | $10^{-8}$ | 128 | 128 | 128 | 1 | 1 | 1 | 1 | 1 | 1 |
| $10^2$ | $10^{-10}$ | 64 | 64 | 64 | 1 | 1 | 1 | 1 | 1 | 1 |
| $10^3$ | $10^{-12}$ | 128 | 128 | 256 | 1 | 1 | 1 | 1 | 1 | 1 |
| $10^4$ | $10^{-14}$ | 1024 | 1024 | 512 / 1024 | 2 | 2 | 2 | 1 | 1 | 2 |

**Table 3** Helix experiment: Solver tolerances and minimum number of load increments needed to reach convergence for the different $DB_\mathrm{red}$ and $SC_\mathrm{red}$, $MX_\mathrm{red}$ and $MX_\mathrm{full}$ formulations. The number of load increments was determined by testing only powers of two. All formulations were discretised with $N = 17$ nodes.

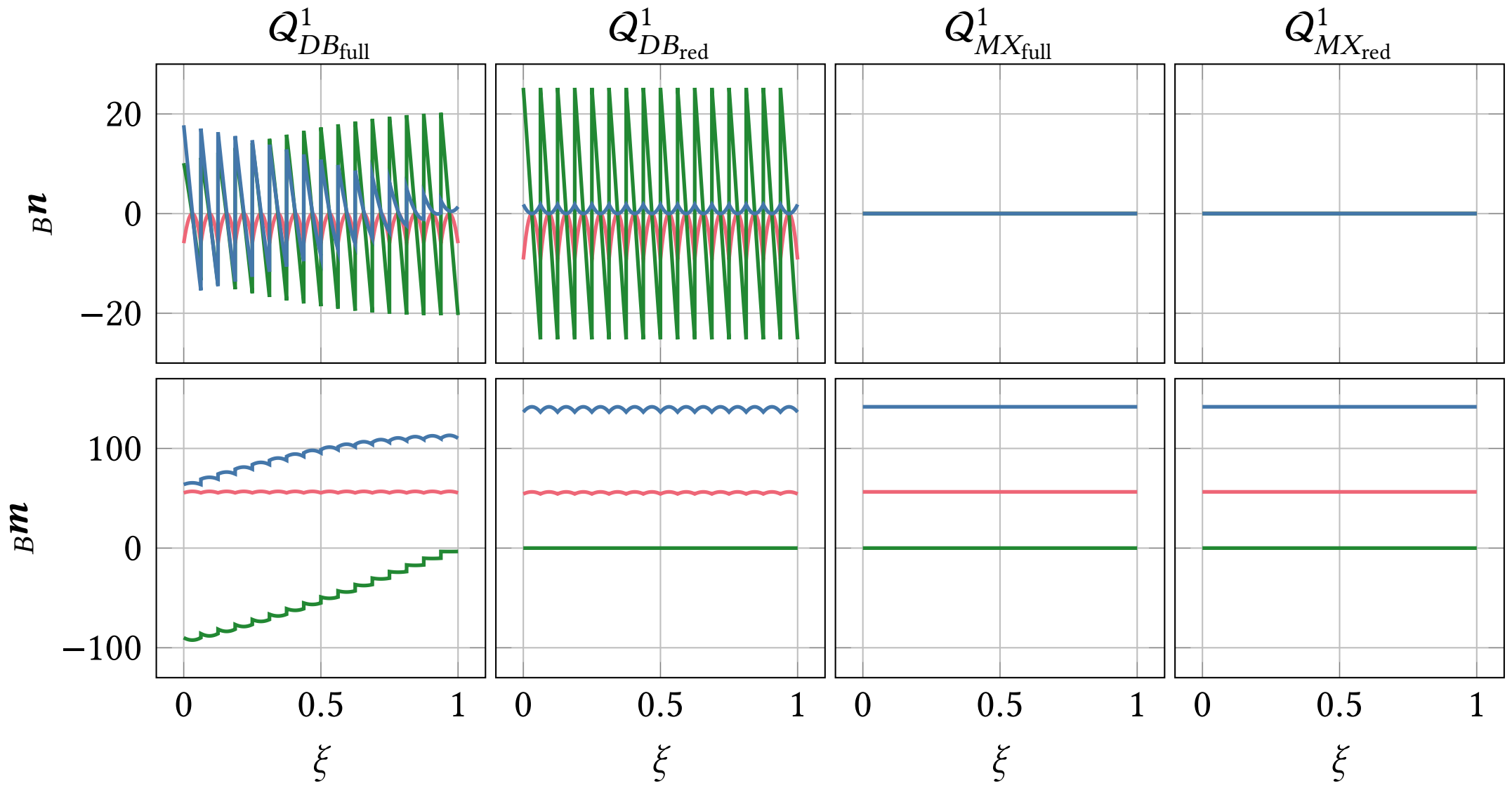

**(a)** Quaternion interpolation with polynomial degree $p = 1$ and $n_{el} = 16$ elements ($N = 17$ nodes)

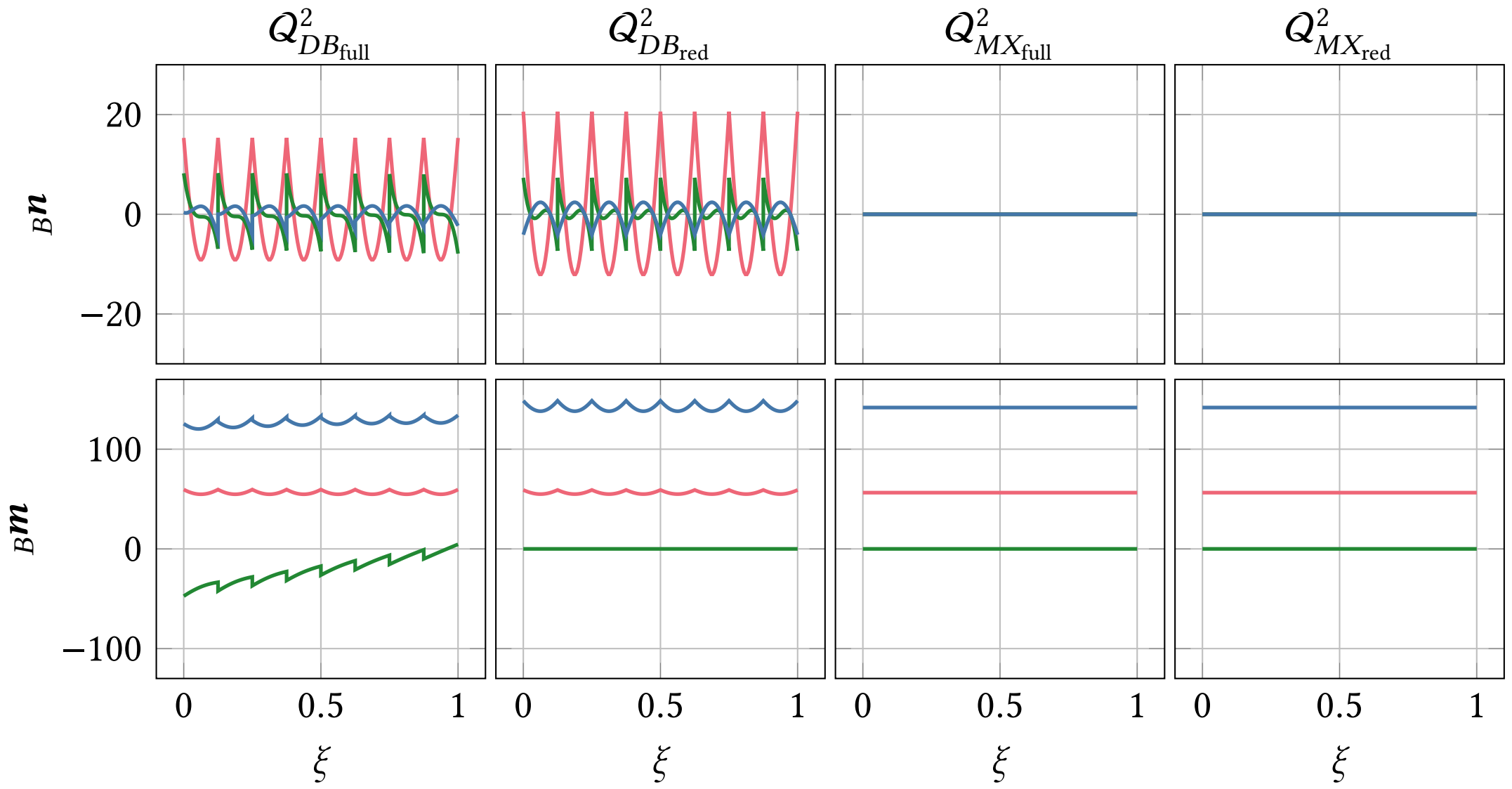

**(b)** Quaternion interpolation with polynomial degree $p = 2$ and $n_{el} = 8$ elements ($N = 17$ nodes)

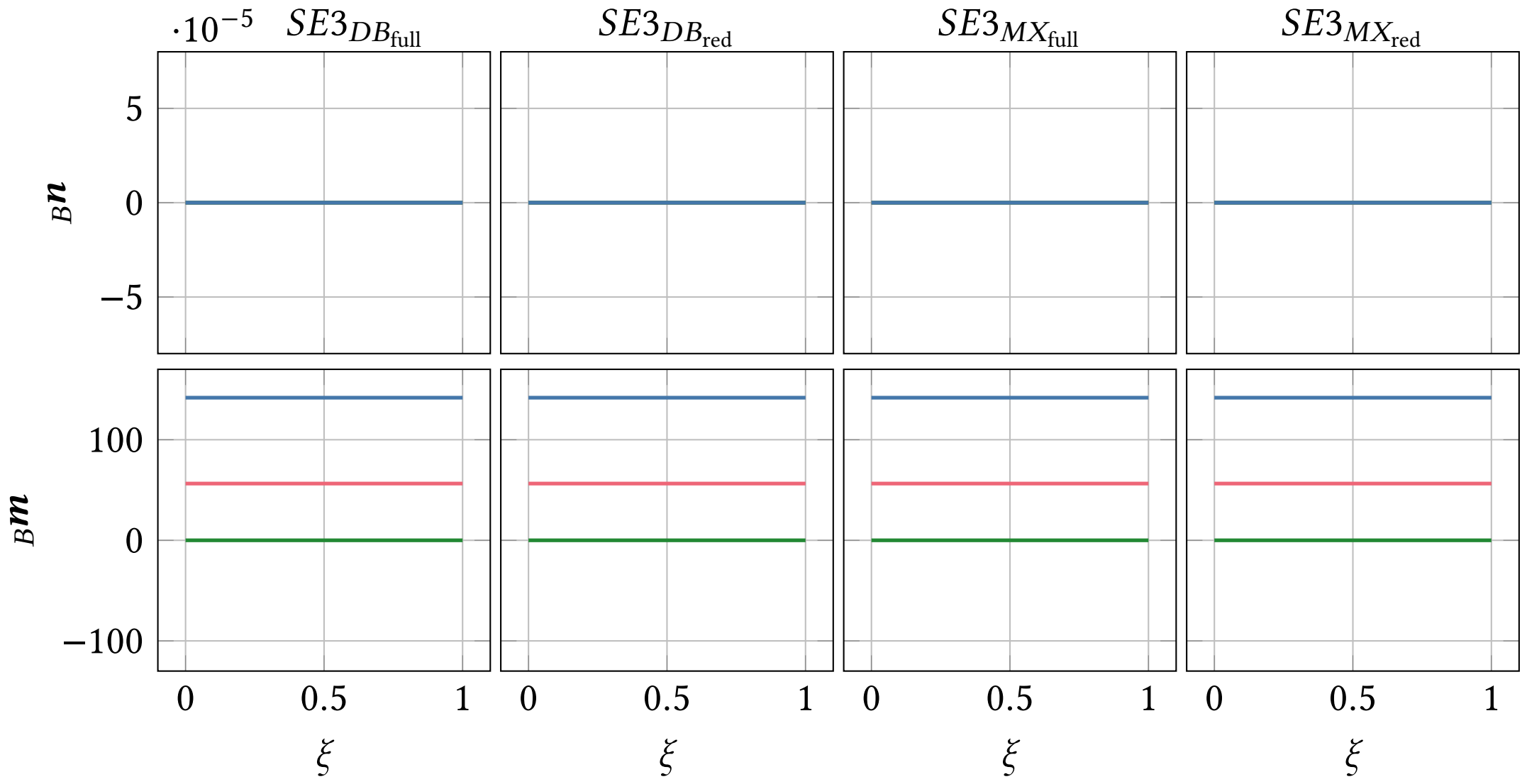

**(c)** $SE3$ interpolation with $n_{el} = 16$ elements ($N = 17$ nodes)

**Figure 5** Helix experiment: Components of the resultant contact forces $_B\boldsymbol{n}$ (top row) and moments $_B\boldsymbol{m}$ (bottom row) of the different formulations for the slenderness $\rho = 10^1$. The components in $\boldsymbol{e}^B_x$, $\boldsymbol{e}^B_y$ and $\boldsymbol{e}^B_z$ direction are shown as (——), (——) and (——), respectively.

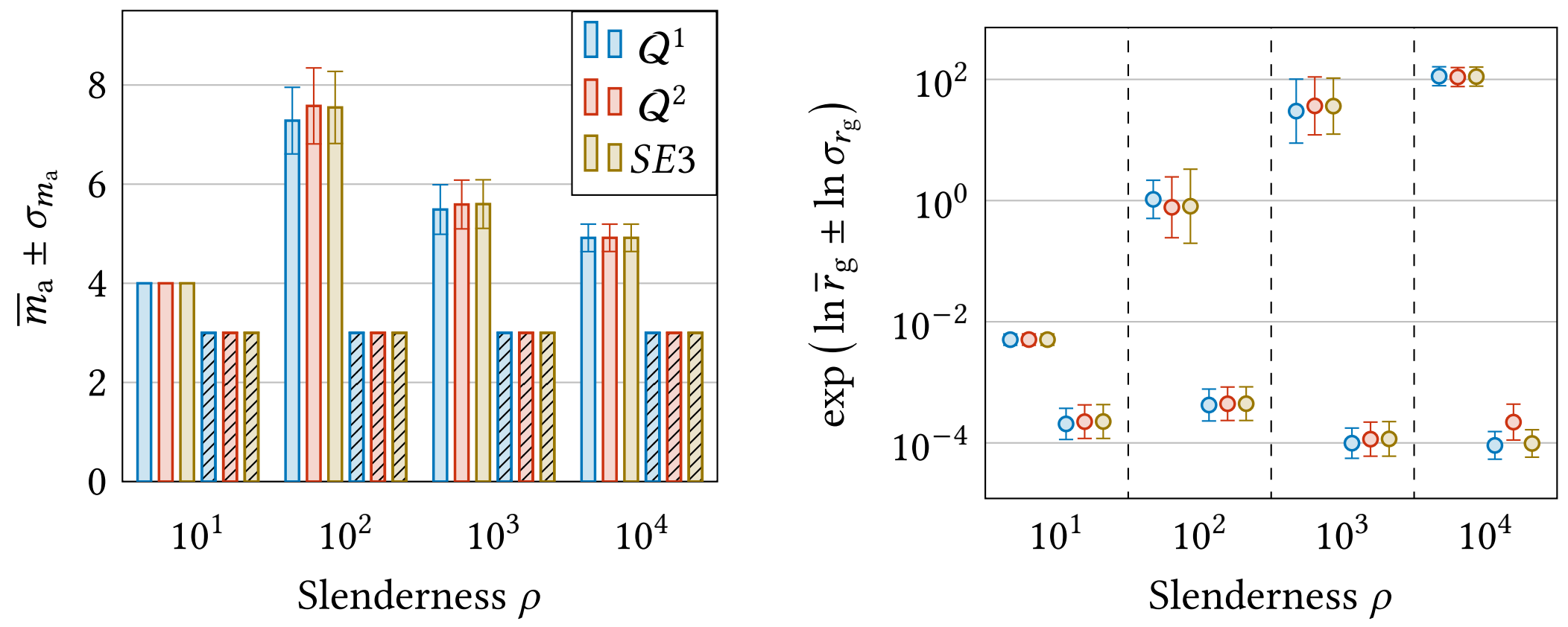

**(a)** Newton iterations **(b)** Local quadratic convergence rate

**Figure 6** Helix experiment: (a) Arithmetic mean $\overline{m}_a$ of Newton iterations and standard deviation $\sigma_{m_a}$ for different slenderness $\rho$ and different formulations. The hatched bars represent mixed formulations. (b) Geometric mean $\overline{r}_g$ of the local quadratic convergence rate and geometric standard deviation $\sigma_{r_g}$ for different slenderness $\rho$ and different formulations. The hatched markers represent mixed formulations.

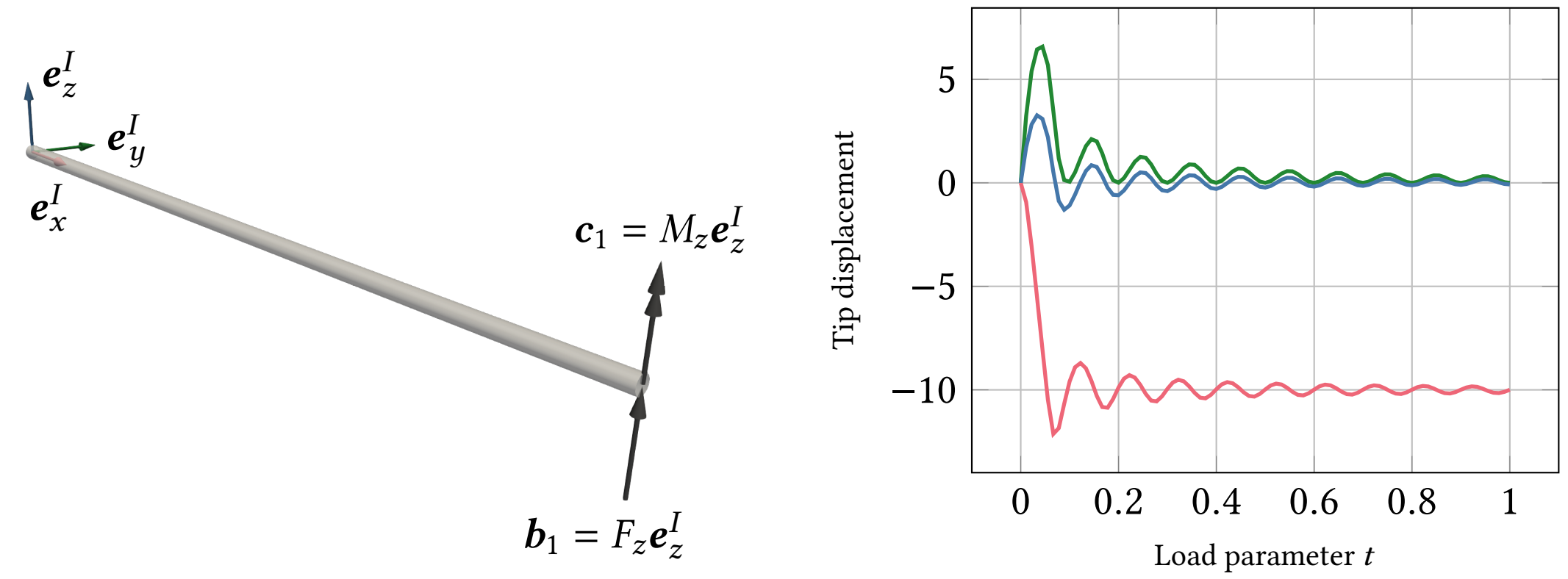

**(a)** Setup of the experiment **(b)** Tip displacement formulation. Because $DB_{\text{full}}$ formulations show locking behaviour, and for mixed formulations the reduced integration method leads to almost the same results as the full integration, we will only compare the $DB_{\text{red}}$ and the $MX_{\text{full}}$ formulations in this experiment.

**Figure 7** Rod bent to a helical form: Overview of the experiment. The tip displacement was obtained using the $Q^2_{MX_{\text{full}}}$ interpolation with $n_{\text{el}} = 30$ elements, i.e., $N = 61$ nodes, with 90 load increments. The components of the tip displacement $(\Delta r_x, \Delta r_y, \Delta r_z) = {}_I\boldsymbol{r}_{OC}(1) - {}_I\boldsymbol{r}^0_{OC}(1)$ are indicated by (——), (——) and (——), respectively.

Table 4 shows the minimum number of load increments to get a converged solution. The number of load increments is again determined by testing powers of two. Similar to the previous experiment, all mixed formulations need much fewer load increments than the displacement-based formulations. For the following analysis, we used 90 load increments for all mixed formulations and 2048 load increments for all displacement-based formulations. It is convenient to introduce the load parameter $t$ as the ratio between the current load increment and the number of load increments.

| $DB_{\text{red}}$ | | | $MX_{\text{full}}$ | | |
|---|---|---|---|---|---|
| $Q^1$ | $Q^2$ | $SE3$ | $Q^1$ | $Q^2$ | $SE3$ |
| 2048 | 2048 | 2048 | 64 | 64 | 64 |

**Table 4** Rod bent to a helical form: Minimum number of load increments needed to reach convergence for the different $DB_{\text{red}}$ and $MX_{\text{full}}$ formulations. The number of load increments was determined by testing only powers of two. All formulations were tested with $N = 61$ nodes.

Figure 7(b) shows the displacement of the rod's free tip with respect to the load parameter $t$. The results are obtained using the $Q^2_{MX_{\text{full}}}$ formulation. Figure 8 shows the deformed configurations of the same formulation every 10 load increments. The horizontal "oscillations" of the free tip are clearly visible (for comparison, see Harsch and Eugster 2023).

The resultant contact forces and moments of the final configuration are shown in Figure 9. The displacement-based formulations of the Quaternion interpolation show the typical fluctuation, leading to a range of almost 10000 for the resultant contact force of the $Q^1_{DB_{\text{red}}}$ formulation. In comparison, the range of the resultant contact force of the mixed formulation is only slightly larger than 50. Additionally, note that the limits of the force axis of the $SE3_{DB_{\text{red}}}$ formulation have been adjusted to match those of the mixed formulations. Results presented in the literature (see Harsch and Eugster 2023; Mäkinen 2007; Ibrahimbegović 1997) are in agreement with the results that we obtained.

## 4.4 Deployment of an elastic ring

In this benchmark experiment, a circular ring with rectangular cross-section is deployed. It has been proposed by Goto et al. (1992) to analyse the buckling phenomenon of a ring when exposed to a moment. It is used as benchmark experiment by various authors (see Da Costa E Silva et al. 2020; Greco et al. 2022; Greco et al. 2024; I. Romero 2004; C. Meier et al. 2014; Smoleński 1999; Borković et al. 2023). The centerline of the ring in the reference configuration has a circular shape with radius $R = 20$. Point $\mathcal{A}$ is located at $\xi = 0$, where the ring is clamped to the origin. To keep the ring closed in the simulation, the other end of the ring is also clamped. At point $\mathcal{C}$, which is located at $\xi = 1/2$, an external moment $M$ along the $\boldsymbol{e}^I_x$ axis is applied, such that the rotated angle around this axis matches a prescribed angle $\theta$. We introduce the basis at point $\mathcal{C}$ by $\boldsymbol{A}_{I\mathcal{C}} = \boldsymbol{A}_{IB}(\xi = 1/2)$.

Furthermore, point $\mathcal{B}$[2] is located at $\xi = 1/4$ and is used to investigate the deformation of the ring. The placement of point $\mathcal{B}$ is introduced as ${}_I\boldsymbol{r}_{O\mathcal{B}} = {}_I\boldsymbol{r}_{OC}(\xi = 1/4)$. The setup of the experiment is shown in Figure 10(a). The cross-section of the ring has height $h = 1$ and width in radial direction $w = 1/3$. With Young's modulus $E = 2.1 \cdot 10^7$ and Poisson's ratio $\nu = 0.3$, the shear modulus is given by $G = E/(2(1+\nu))$. The stiffnesses are given by $k_{\text{e}} = EA$, $k_{\text{s}_y} = k_{\text{s}_z} = GA$, $k_{\text{t}} = G \cdot 9.753 \cdot 10^{-3}$, $k_{\text{b}_y} = EI_y$ and $k_{\text{b}_z} = EI_z$ with $A = hw$, $I_y = wh^3/12$ and $I_z = hw^3/12$. Note that the value for the torsional stiffness $k_{\text{t}}$ is not given by the geometric value of the cross-section, instead a custom value, as done by Greco et al. (2022); Greco et al. (2024), is taken.

The ring is discretised with $n_{\text{el}} = 20$ elements of the $Q^2_{MX_{\text{full}}}$ formulation. Since the number of elements is even, point $\mathcal{C}$ is located at the boundary of adjacent elements so that the requirement on external loads is fulfilled as described after (23). The tolerance for Newton's solver is given by $\epsilon = 10^{-6}$, and the angle of rotation in point $\mathcal{C}$ is increased within 120 increments from $\theta = 0$ to $\theta = 4\pi$. Figure 10(b) shows the required moment $M$ to enforce the rotated angle $\theta$ at point $\mathcal{C}$. The configurations for different rotated angles are shown in Figure 11. The components of the displacement of point $\mathcal{B}$ are shown in Figure 12, where the components are introduced as $(\Delta r_x, \Delta r_y, \Delta r_z) = {}_I\boldsymbol{r}_{O\mathcal{B}} - {}_I\boldsymbol{r}^0_{O\mathcal{B}}$. Due to symmetry, one could also just simulate one half of the ring. This divides also the required moment $M$ by two (*e.g.* Da Costa E Silva et al. 2020; Greco et al. 2024). The results obtained using the herein proposed rod formulation show good agreement with results from literature.

[2] In literature, the points are labeled as $A$, $B$ and $C$. To avoid confusion with the centerline curve ${}_I\boldsymbol{r}_{OC}$ and the transformation matrix describing the cross-section orientation $\boldsymbol{A}_{IB}$, we label them $\mathcal{A}$, $\mathcal{B}$ and $\mathcal{C}$, respectively.

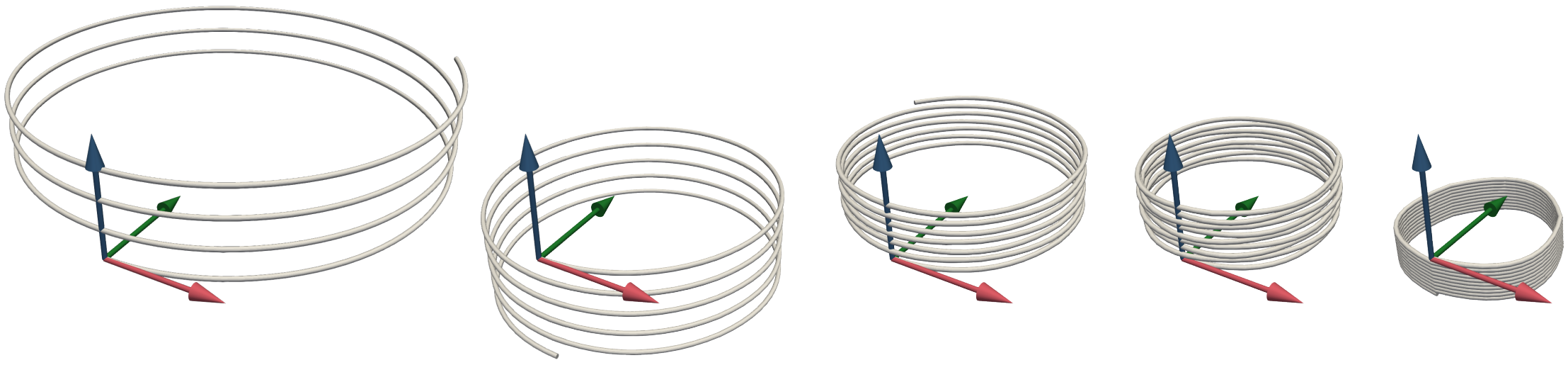

**Figure 8** Rod bent to a helical form: Deformed configurations for different values of the load parameter $t$. The deformed configurations were obtained using the $Q^2_{MX_{\text{full}}}$ interpolation with $n_{\text{el}} = 30$ elements, i.e., $N = 61$ nodes. From the left to the right, the deformed configurations were obtained after 30, 45, 60, 75 and 90 load increments, corresponding to $t \in \{1/3, 1/2, 2/3, 5/6, 1\}$. The shown axes of the $I$-system are scaled by $1/4$.

## 4.5 Constrained cantilever

This experiment shows the capability of the proposed formulation to realise constraints on the deformation. An initially straight rod with length $L = 2\pi$ is considered. The rod is clamped at $\xi = 0$ with $\boldsymbol{A}_{IB}(0) = \boldsymbol{I}$ and ${}_I\boldsymbol{r}_{OC}(0) = \boldsymbol{0}$.

Two load cases are considered, where in the first load case only an external force ${}_I\boldsymbol{b}_1 = (0, -P, 0)$ is applied on the free tip. The scalar force $P$ is given in terms of the bending stiffness around the $\boldsymbol{e}_z^B$ axis, $k_{\mathrm{b}_z} = 2$, and is given by $P = k_{\mathrm{b}_z}(\alpha/L)^2$, where $\alpha^2 \in [0, 10]$ can be interpreted as scaled load parameter. In the second load case, an additional external moment ${}_B\boldsymbol{c}_1 = (0, 0, eP)$ with the lever arm $e = 2.5$ is applied at the free tip. Note, that both load cases result in a deformation only in the ${}_I\boldsymbol{e}_x^I$-${}_I\boldsymbol{e}_y^I$-plane. Furthermore, we consider three different (constrained)

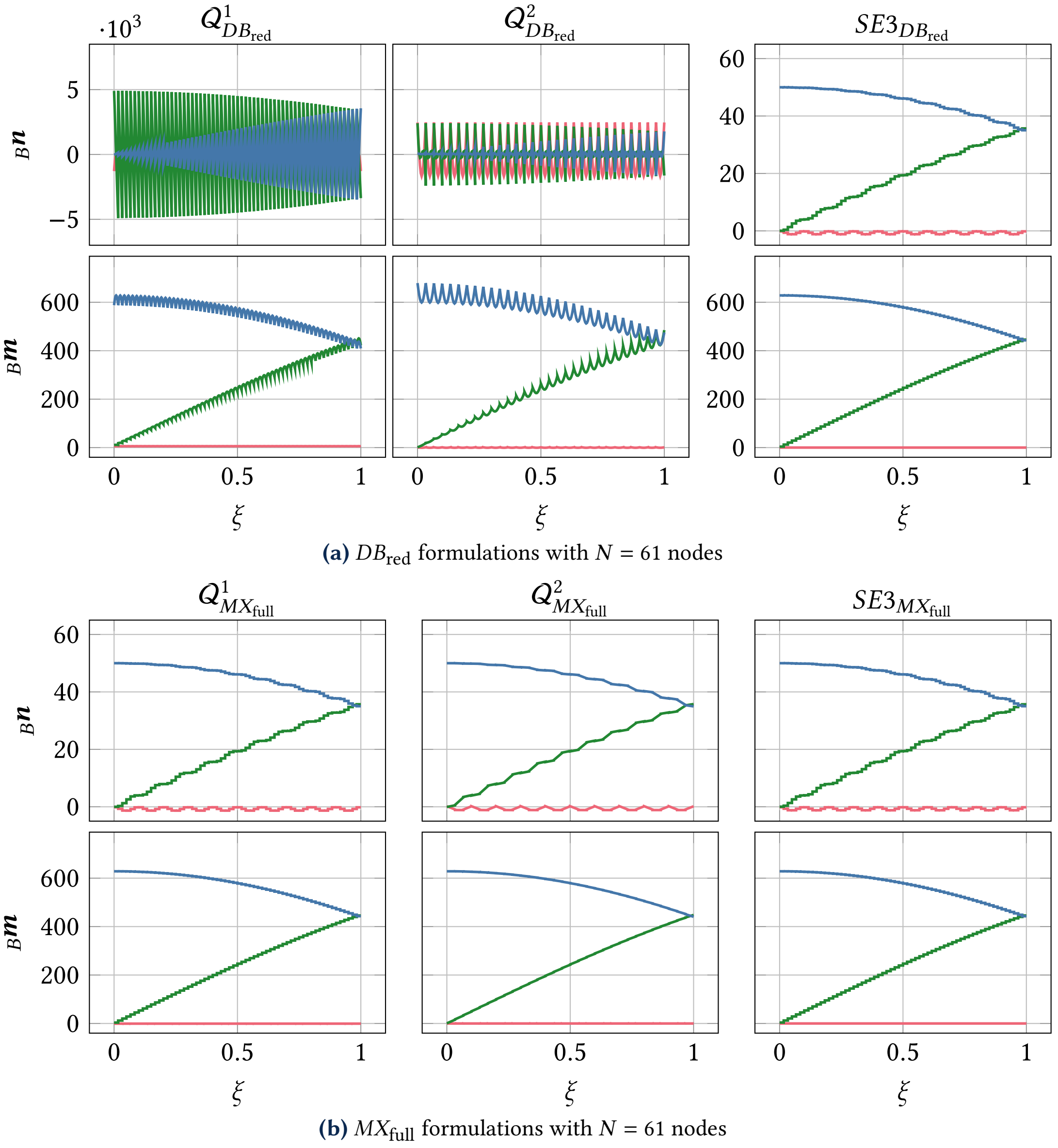

**Figure 9** Rod bent to a helical form: Components of the resultant contact forces ${}_B\boldsymbol{n}$ (top row) and moments ${}_B\boldsymbol{m}$ (bottom row) of the different formulations. The components in $\boldsymbol{e}_x^B$, $\boldsymbol{e}_y^B$ and $\boldsymbol{e}_z^B$ direction are shown as (—), (—) and (—), respectively.

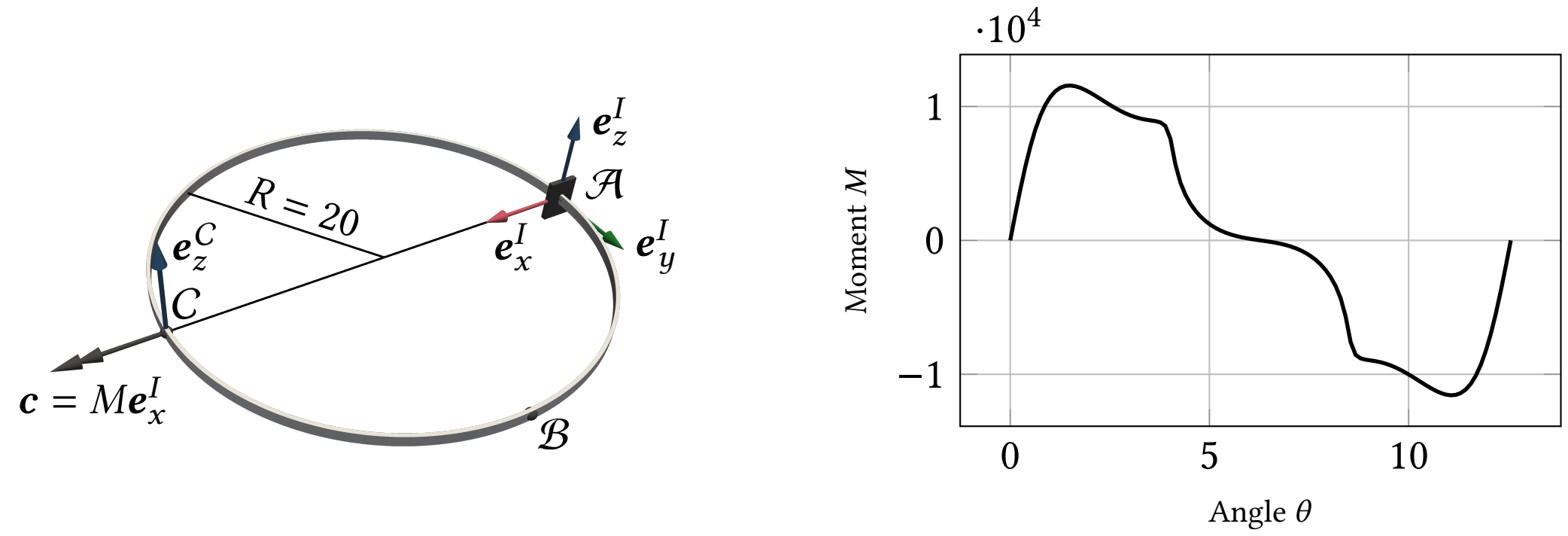

**(a)** Setup of the experiment **(b)** Moment $M$ in dependency of the rotated angle $\theta$

**Figure 10** Deployment of an elastic ring: Overview of the experiment. The experiment was computed using the $Q^2_{MX_{\text{full}}}$ formulation with $n_{\text{el}} = 20$ elements, i.e., $N = 41$ nodes with 120 increments in the prescribed angle $\theta$.

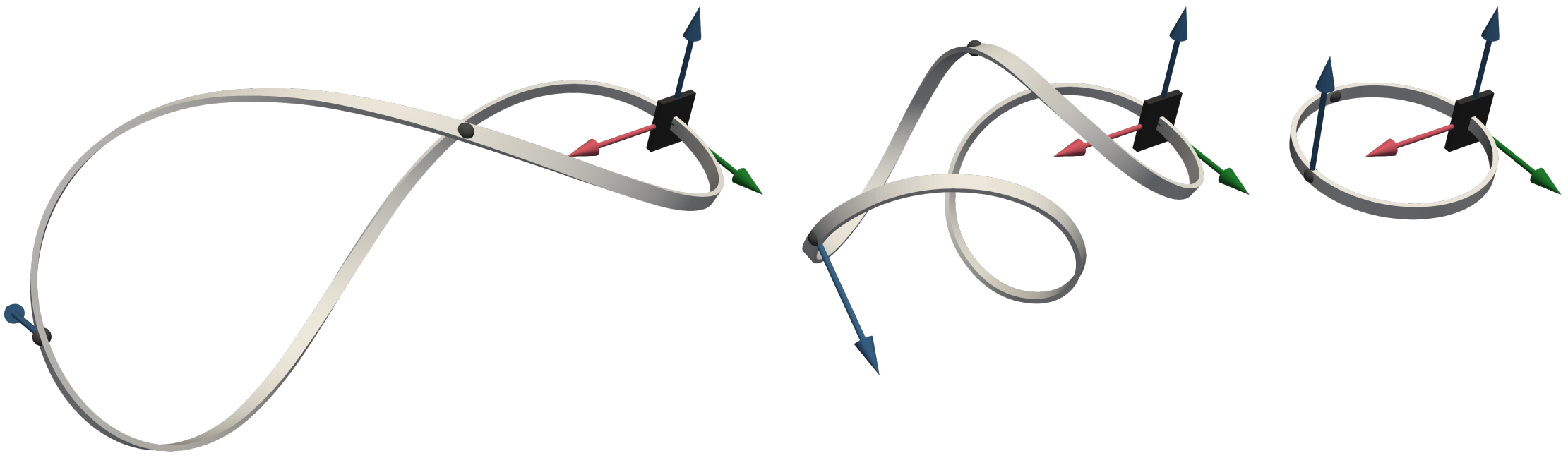

**Figure 11** Deployment of an elastic ring: Deformed configurations for different rotated angles $\theta$. The deformed configurations were obtained using the $Q^2_{MX_{\text{full}}}$ interpolation with $n_{\text{el}} = 20$ elements, *i.e.*, $N = 41$ nodes. From the left to the right, the deformed configurations were obtained after 20, 40 and 60 increments of the angle, corresponding to $\theta \in \{2\pi/3, 4\pi/3, 2\pi\}$.

constitutive laws for the simulation, resulting in an unconstrained rod, a shear-stiff rod and an inextensible shear-stiff rod. Their inverse stiffness matrices are listed in Table 5. As the experiment has a planar setup, the chosen values of the shear stiffness in $e_z^B$ direction, the torsional stiffness and the bending stiffness around $e_y^B$ are not affecting the results. Considering the inextensible shear-stiff rod, also known as Euler's elastica, the analytical solution can also be found using elliptic integrals as shown in Harsch et al. (2021); Frisch-Fay (1962); Bisshopp and Drucker (1945).

The rod was discretised using $n_{\text{el}} = 4$ of the $Q^2_{MX_{\text{full}}}$ formulation. The load was applied within 40 load increments and the tolerance for Newton's solver was $\epsilon = 10^{-12}$.

Figure 13 shows the centerline curves for different load parameters $\alpha \in \{0, 1, 2, 4, 10\}$ as well as the normalised tip placement. The load case with only the force acting is shown in Figure 13(a) and the load case with force and moment acting is shown in Figure 13(b). The results of the inextensible shear-stiff rod are in agreement with the analytical solution, as visible in the middle and the right plots. The results of the unconstrained and shear-stiff rods are in agreement with the results from Harsch et al. (2021).

Finally, we use this experiment to study the condition number of the iteration matrices. For each formulation, we compute in each load step and in each iteration the condition number of the iteration matrix $\partial \boldsymbol{f}/\partial \boldsymbol{x}(\boldsymbol{x})$, containing also the equations for the clamping on the left side.

| Unconstrained | Shear-stiff | Inextensible shear-stiff |
|---|---|---|
| $C_{\gamma}^{-1} = \text{diag}(0.2, 1, 1)$ | $C_{\gamma}^{-1} = \text{diag}(0.2, 0, 0)$ | $C_{\gamma}^{-1} = \text{diag}(0, 0, 0)$ |
| $C_{\kappa}^{-1} = \text{diag}(2, 0.5, 0.5)$ | $C_{\kappa}^{-1} = \text{diag}(2, 0.5, 0.5)$ | $C_{\kappa}^{-1} = \text{diag}(2, 0.5, 0.5)$ |

**Table 5** Constrained cantilever experiment: Inverse elasticity matrix for the different constraints.

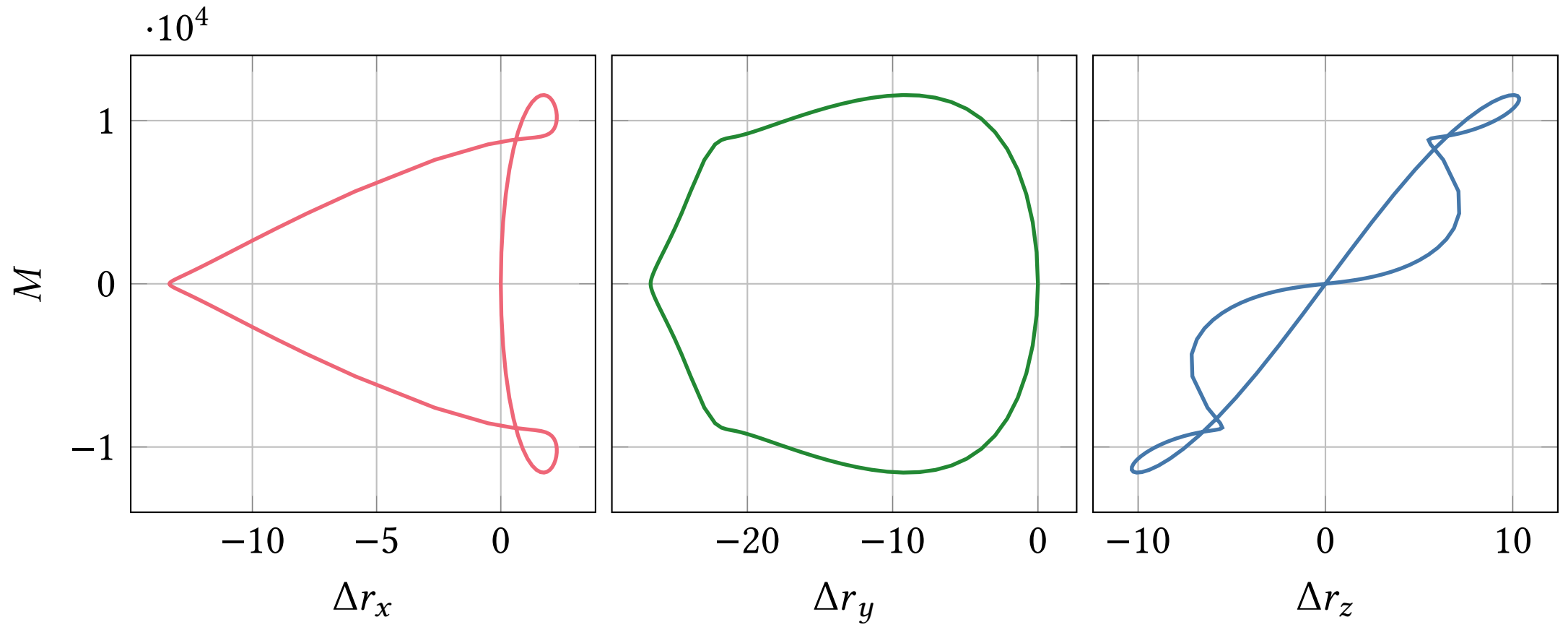

**Figure 12** Deployment of an elastic ring: Displacement of point $B$. The experiment was performed using the $Q^2_{MX_{\text{full}}}$ formulation with $n_{\text{el}} = 20$ elements, *i.e.*, $N = 41$ nodes with 120 increments in the prescribed angle $\theta$. The components of the displacement are $(\Delta r_x, \Delta r_y, \Delta r_z) = {}_I\boldsymbol{r}_{O\mathcal{B}} - {}_I\boldsymbol{r}^0_{O\mathcal{B}}$ and coloured in (——), (——) and (——), respectively.

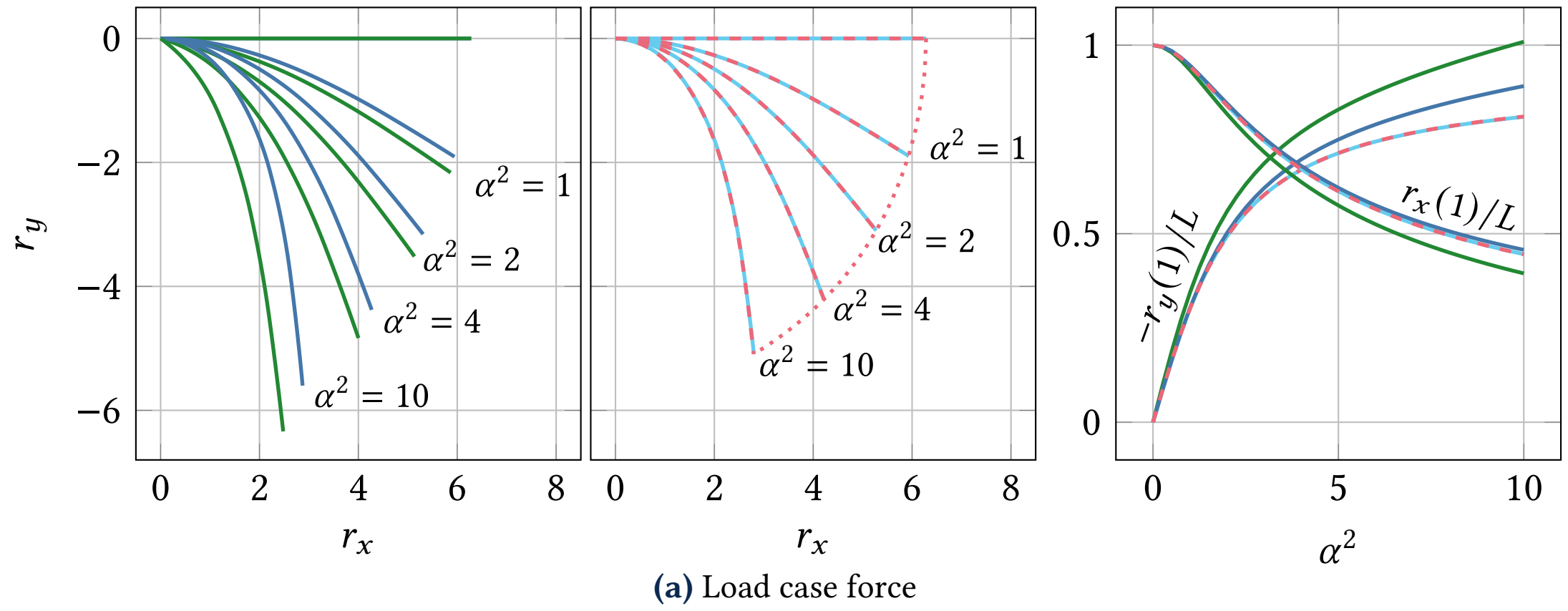

**(a)** Load case force

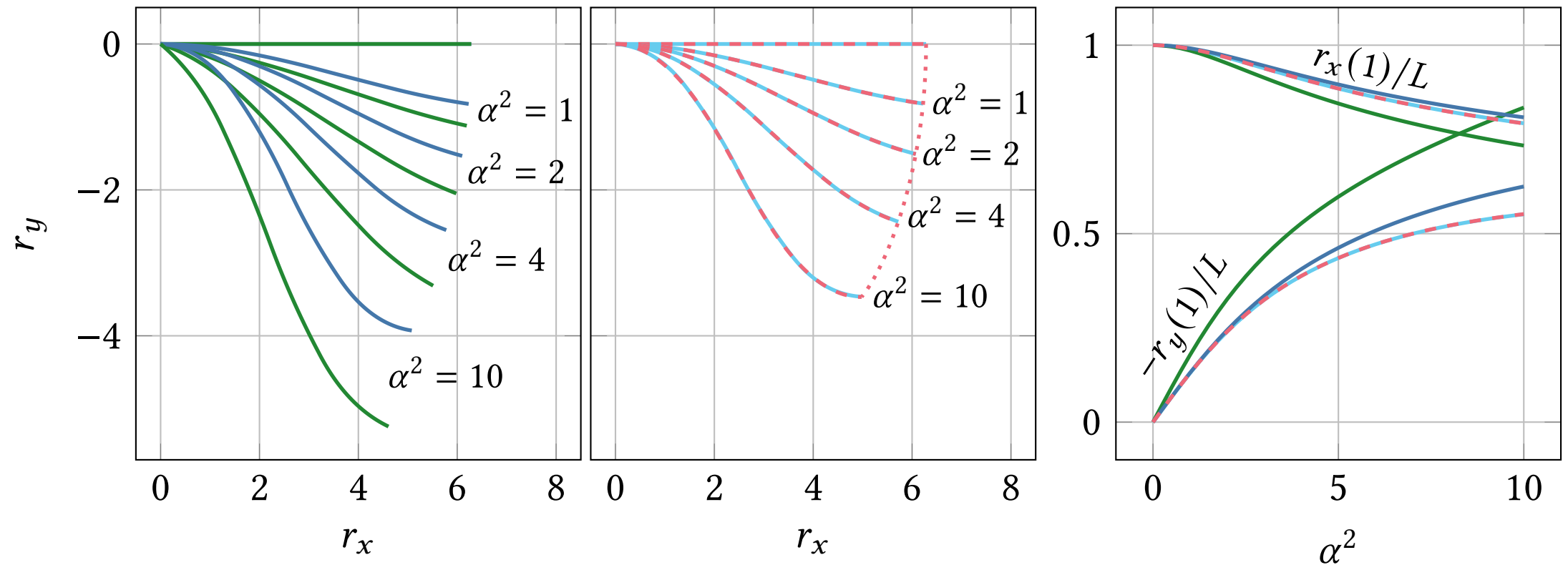

**(b)** Load case force and moment

**Figure 13** Constrained cantilever experiment: (a) Load case only force, (b) load case force and moment. The left plot in each row shows the centerline of the unconstrained (——) Cosserat rod and the centerline of the shear-rigid constrained (——) Cosserat rod. The middle plot of each row shows the centerline of the inextensible and shear-rigid constrained (——) Cosserat rod and the centerline of the analytical solution (- - -). The right plot in each row shows the normalised tip displacements for all rods. The components of the placement function are given by ${}_I\boldsymbol{r}_{OC}(\xi) = (r_x(\xi), r_y(\xi), 0)$.

The minimum, maximum and average values are given in Table 6. To have good comparability between all formulations, also $Q^2_{MX_{\text{red}}}$ and $Q^2_{DB_{\text{red}}}$ formulations are considered. The condition numbers for the $MX_{\text{red}}$ and $MX_{\text{full}}$ formulations are very similar. Similarly, the differences between the different constraints (unconstrained, shear-stiff and inextensible shear-stiff) are small within each formulation. Large differences can be observed between the displacement-based and the mixed formulations. It is important to note that the condition numbers provide valuable insight into the stability of the linear Newton updates, but do not characterise the difficulty of the nonlinear problem itself.

### 4.6 Lateral-torsional buckling of a ribbon

The final experiment validates the proposed formulation in a scenario where many deformation constraints are applied. Following the buckling analysis of Audoly and Neukirch (2021); V. Romero et al. (2021) of a ribbon, we use the presented formulation to replicate this analysis. Using the cross-section geometry of width $a = 0.05$ and thickness $t = 0.002$, together with the material properties of Young's modulus $E = 2.1 \cdot 10^{11}$ and Poisson's ratio $\nu = 0.4$, the characteristic curvature $\varepsilon^\star = [12(1-\nu^2)]^{-1/2} t/a^2$ is computed. We use two different lengths, given by $L_A = 0.0157/\varepsilon^\star$ and $L_B = 0.315/\varepsilon^\star$. The rod is clamped at one side, and initialized such that the length, thickness and width are in $\boldsymbol{e}_x^I$, $\boldsymbol{e}_y^I$ and $\boldsymbol{e}_z^I$ direction, respectively, as shown in Figure 14(b). The constrained deformations of the ribbon are extension, both shear deformations and bending around $\boldsymbol{e}_y^B$. For the remaining deformations, i.e., torsion $\varepsilon_t$ and bending around $\boldsymbol{e}_z^B$, $\varepsilon_b$, three different material laws are used. They are defined by the strain energy density functions

$$\begin{aligned} W_{\text{quad}}(\varepsilon_t, \varepsilon_b) &= \frac{1}{2}\frac{Eat^3}{6(1+\nu)}\varepsilon_t^2 + \frac{1}{2}\frac{Eat^3}{12}\varepsilon_b^2 \,, \\ W_{\text{rib}}(\varepsilon_t, \varepsilon_b) &= \frac{1}{2}\frac{Eat^3}{6(1+\nu)}\varepsilon_t^2 + \frac{1}{2}\frac{Eat^3}{12}\varepsilon_b^2 + \frac{1}{2}\frac{Ea^5 t}{2}\left(\nu\varepsilon_b^2 + \varepsilon_t^2\right)^2 \varphi\left(\frac{\varepsilon_b}{\varepsilon^\star}\right), \\ W_{\text{sdw}}(\varepsilon_t, \varepsilon_b) &= \frac{1}{2}\frac{Eat^3}{12(1-\nu^2)}\varepsilon_{\text{sdw}}^2(\varepsilon_t, \varepsilon_b)\,, \end{aligned} \tag{30}$$

with the auxiliary function $\varphi$ as defined by Audoly and Neukirch (2021)[Equations 2.8 and B.1], and the auxiliary strain

$$\varepsilon_{\text{sdw}}(\varepsilon_t, \varepsilon_b) = \begin{cases} \frac{\varepsilon_t^2+\varepsilon_b^2}{\varepsilon_b}\,, & \varepsilon_b > \varepsilon_s\,, \\ \frac{\varepsilon_t^2+\varepsilon_s^2}{\varepsilon_s} + \frac{1}{2\varepsilon_s}\left(1 - \left(\frac{\varepsilon_t}{\varepsilon_s}\right)^2\right)\left(\varepsilon_b^2 - \varepsilon_s^2\right), & \varepsilon_b \leqslant \varepsilon_s\,. \end{cases} \tag{31}$$

As the Sadowsky potential follows from the ribbon potential in the limit case of large bending strains, *i.e.*, $\varepsilon_b \gg \varepsilon^\star$, the quadratic function as illustrated in the second line of (31) is used in the case of $\varepsilon_b \leqslant \varepsilon_s = 10^{-6}$ to overcome singularities around $\varepsilon_b = 0$. Using a load parameter $\hat{t} \in [0, 1]$, which is increased in 100 steps, the distributed load is given by

$$\begin{aligned} {}_I\boldsymbol{b}(\xi, \hat{t}) &= \frac{Eat^3}{12L^3}\begin{cases} -(5\hat{t})(4\gamma_{\text{crit}}^{\text{rib}})\left({}_I\boldsymbol{e}_z^I + \frac{\pi}{40}(1-5\hat{t})\,{}_I\boldsymbol{e}_y^I\right), & \hat{t} \leqslant 1/5\,, \\ -\gamma\left(\frac{5\hat{t}-1}{4}\right){}_I\boldsymbol{e}_z^I\,, & 1/5 < \hat{t}\,, \end{cases} \\ \gamma(\tau) &= \begin{cases} 4\gamma_{\text{crit}}^{\text{rib}} + \left(\bar{\gamma}^{\square} - 4\gamma_{\text{crit}}^{\text{rib}}\right)\frac{e^{-6\tau}-1}{e^{-6\bar{\tau}}-1}\,, & \tau \leqslant \bar{\tau}\,, \\ \bar{\gamma}^{\square}\frac{\tau-1}{\bar{\tau}-1}\,, & \bar{\tau} < \tau\,, \end{cases} \end{aligned} \tag{32}$$

| | Unconstrained | | | Shear-stiff | | | Inextensible shear-stiff | | |
|---|---|---|---|---|---|---|---|---|---|
| | min | max | average | min | max | average | min | max | average |
| $Q^2_{MX_{\text{full}}}$ | 213.8 | 455.1 | 289.1 | 198.4 | 383.1 | 259.5 | 185.4 | 355.1 | 243.0 |
| | 341.3 | 713.4 | 459.1 | 353.8 | 793.0 | 650.0 | 341.1 | 833.5 | 659.8 |
| $Q^2_{MX_{\text{red}}}$ | 213.2 | 455.2 | 289.2 | 199.2 | 383.1 | 260.0 | 186.2 | 355.1 | 243.5 |
| | 341.2 | 716.4 | 459.3 | 353.5 | 794.9 | 651.1 | 340.9 | 835.6 | 660.9 |
| $Q^2_{DB_{\text{red}}}$ | 1435.2 | 4250.6 | 2456.6 | 1225.6 | 15779.8 | 2657.8 | 1154.2 | 6182.0 | 2243.5 |
| | 4075.3 | 4412.7 | 4316.4 | 2034.7 | 5394.5 | 3519.7 | 3670.8 | 4910.1 | 4082.9 |

**Table 6** Constrained cantilever experiment: Condition numbers of the different formulations. While the first row of each formulation corresponds to the load case force, the second row corresponds to the load case force and moment.

with $\bar{\tau} = 15/16$ and $\bar{\gamma}^{\square} = (1 - 10^{-3})\gamma^{\square}_{\text{crit}}$, where $\gamma^{\square}_{\text{crit}}$ is either $\gamma^{\text{rib}}_{\text{crit}} = \gamma^{\text{quad}}_{\text{crit}} = 18.178/\sqrt{1+\nu}$ or $\gamma^{\text{sdw}}_{\text{crit}} = 21.491/(1-\nu^2)$, depending on the chosen potential from (30). We use this specific load profile to force the ribbon into the buckled configuration and to have a high resolution around the bifurcation point while still using the standard Newton–Raphson solver. For $\hat{t} \leqslant 1/5$, the load increases and includes a component in $\boldsymbol{e}^I_y$ direction, whereas for $\hat{t} > 1/5$ the load decreases and acts purely in (negative) $\boldsymbol{e}^I_z$ direction.

We refer to $\gamma$ as the dimensionless force and normalise it by division by $\gamma^{\text{rib}}_{\text{crit}}$, which is shown in Figure 14 and used in the buckling analysis. As the analysis of Audoly and Neukirch (2021) shows, the solution of is categorised by the three dimensionless parameters $(\nu, L\varepsilon^{\star}, \gamma)$. Therefore, the choice of parameters is to some extent arbitrary.

The rod is discretised with $n_{\text{el}} = 32$ elements and the tolerance for the Newton–Raphson solver is set to $10^{-11}$. When the quadratic potential is used, the full mixed formulation is used. For the ribbon potential and the Sadowsky potential, the displacement-based formulation is used for the unconstrained deformation. The constrained deformations follow the briefly introduced

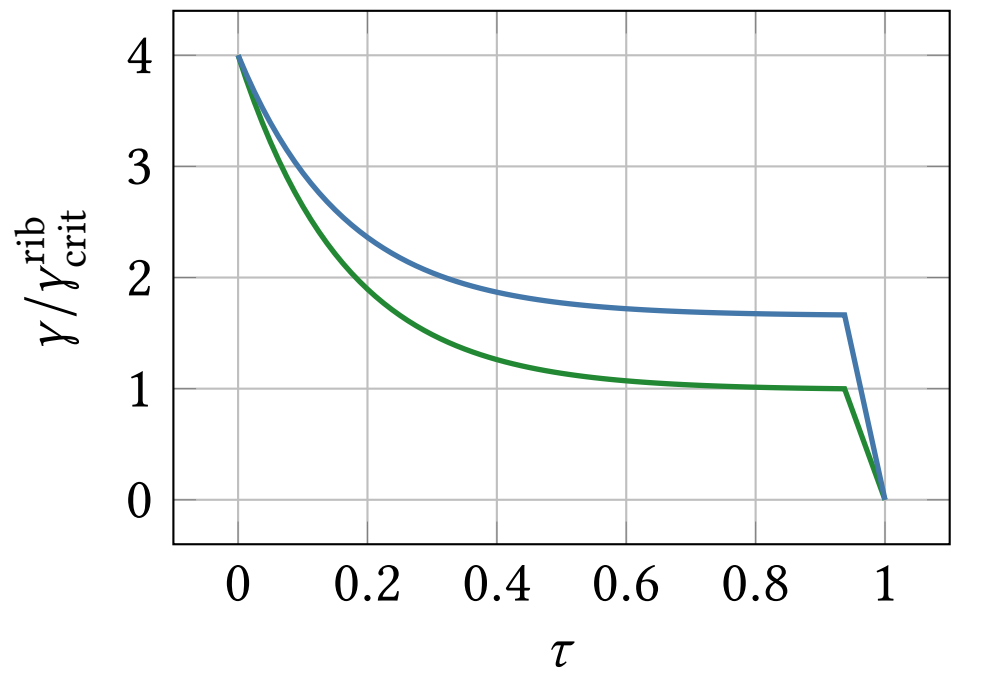

**(a)** Normalised dimensionless force $\gamma(\tau)/\gamma^{\text{rib}}_{\text{crit}}$ over $\tau$ for the different material laws, $W_{\text{quad}}$ & $W_{\text{rib}}$ (——) and $W_{\text{sdw}}$ (——).

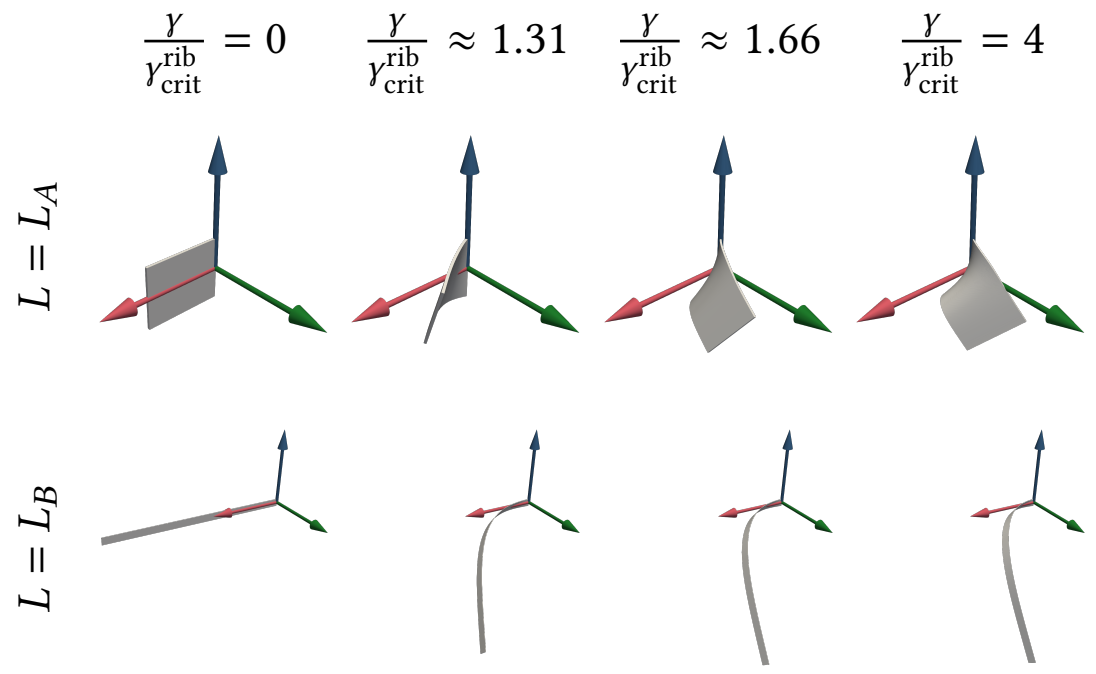

**(b)** Deformed configurations obtained using the $W_{\text{rib}}$ potential.

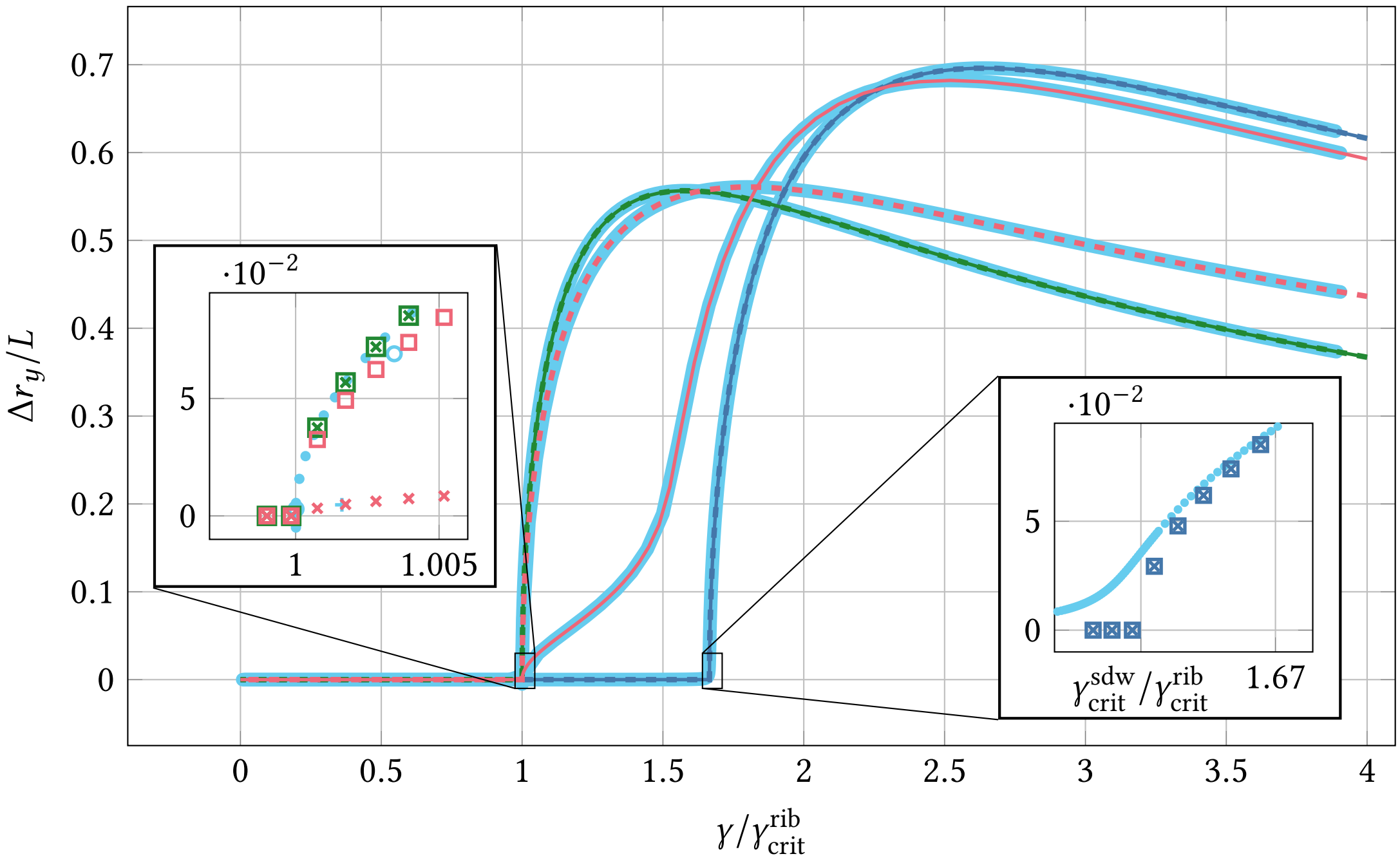

**(c)** Normalised tip displacement in $\boldsymbol{e}^I_y$ direction, $\Delta r_y = {}_I\boldsymbol{r}^{\top}_{OC}(1)\,{}_I\boldsymbol{e}^I_y$ over normalised force. The different material laws $W_{\text{rib}}$, $W_{\text{sdw}}$ and $W_{\text{quad}}$ are indicated by the colours (——), (——) and (——), respectively. The solid lines correspond to the ribbon with length $L_A$, the dashed lines correspond to the ribbon with length $L_B$. In the zoomed sections, length $L_A$ and $L_B$ are indicated by crosses and squares respectively. The solutions from Audoly and Neukirch (2021) are indicated by (——).

**Figure 14** Lateral-torsional buckling experiment: (a) normalised load factor, (b) deformed configurations, (c) normalised tip displacement over normalised load factor.

scheme from Section 3.4. For all simulations, the quadratic quaternion interpolation with reduced integration is used.

Different configurations of the rod with $W_{\text{rib}}$ potential for both lengths and different loads can be seen in Figure 14(b). Figure 14(c) shows the bifurcation diagram as the normalised tip displacement in $\boldsymbol{e}_y^I$ direction over the normalised load. The figure shows the results of the three different material laws and both lengths of our simulation together with the results from Audoly and Neukirch (2021). Our model indicates the same bifurcation points as the one found by Audoly and Neukirch (2021) for all material models. Additionally, the remainder of the normalised displacement diagram is in very good agreement with the one from the literature. For the Sadowsky potential, only the straight configurations, obtained for $\gamma < \gamma_{\text{crit}}^{\text{sdw}}$, are using the approximation of $\varepsilon_{\text{sdw}}$ in the second line of (31).

## 5 Conclusion

In this work, we have proposed a total Lagrangian mixed Petrov–Galerkin finite element formulation for Cosserat rods that addresses fundamental numerical challenges of classical rod formulations. The proposed method employs a Lagrange quaternion interpolation that relies on a minimal set of mathematical concepts while avoiding singularities and ensuring objectivity without additional complexity. In contrast to intrinsically locking-free interpolations based on Lie groups, such as SE(3), our approach avoids the need for trigonometric functions or complex tangent mappings. The formulation is based on a Petrov–Galerkin projection of the internal virtual work and employs the Hellinger–Reissner principle, which introduces independent resultant contact forces and moments and thus leads to a truly locking-free finite element formulation. Furthermore, since the same Petrov–Galerkin projection is used, our formulation can be seamlessly combined with other kinematic interpolations that use the same projection, such as the SE(3) interpolation (as exemplified in this paper), the $\mathbb{R}^{12}$ or the $\mathbb{R}^3 \times \text{SO}(3)$ interpolation outlined in Eugster and Harsch (2023). With the Petrov–Galerkin method, the expressions for the virtual work contributions become much simpler than with a Bubnov–Galerkin approach, since the interpolation of the virtual rotation is not dependent on the generalised coordinates. Additionally, no Lagrange multiplier is required to ensure unit length of the nodal quaternions. A drawback of the Petrov–Galerkin method is that the iteration matrix for the Newton–Raphson is not symmetric, which might need more computation time when iterative solvers are used.

The numerical examples clearly demonstrate the advantages of the proposed mixed formulation over displacement-based approaches. Locking effects, including membrane and shear locking, are completely eliminated, obviating the need for reduced integration or other numerical strategies. Reduced integration mitigates locking in displacement-based formulations in the sense that the centerline and orientation converge to the correct solutions. However, reduced integration does not prevent large fluctuations in deformations within an element and, consequently, in the resultant contact forces and moments. Therefore, the mixed formulation allows for the application of other spatial integration strategies, such as Gauss–Lobatto quadrature, that may be more suitable for certain applications.

Although the mixed formulation requires additional degrees of freedom, it increases computational robustness and efficiency. Compared to the displacement-based formulation, it achieves faster convergence in Newton–Raphson iterations, requires fewer load increments, and reduces the number of iterations per load increment. Moreover, the convergence rate of the Newton–Raphson scheme remains unaffected by an increasing slenderness of the rod. The improved computational robustness of the mixed formulation can also be observed in the intrinsically locking-free $SE(3)$ interpolation. Another nice feature is that the mixed formulation naturally incorporates constrained rod theories, such as Kirchhoff–Love rods, providing a unified framework that encompasses Cosserat and Kirchhoff–Love rod models.

In this study, we examined the robustness of the mixed formulation in static analyses. Subsequent steps will involve analysing the robustness of the mixed formulation in dynamic simulations. To facilitate this analysis, the internal virtual work expression in Harsch et al. (2023) must be replaced by the mixed formulation proposed in this paper. In combination with an appropriate time-integration scheme, we aim to increase robustness and computational efficiency in dynamic simulations to pave the way for real-time simulation applications, which are particularly relevant for control applications, such as model predictive control.

# A Boundary value problem of the Cosserat rod

The boundary value problem of the Cosserat rod can be derived for both formulations of the internal virtual work (4) and (11). Starting with the first one, we perform an integration by parts on (4) and write ${}_I\boldsymbol{n} = \boldsymbol{A}_{IB}\, {}_B\boldsymbol{n}$, which yields

$$\begin{aligned}\delta W^{\text{int}} &= -\left[{}_I\delta\boldsymbol{r}_C^{\top}\, {}_I\boldsymbol{n} + {}_B\delta\boldsymbol{\phi}_{IB}^{\top}\, {}_B\boldsymbol{m}\right]_{\mathcal{J}} \\ &\quad + \int_{\mathcal{J}} \left\{ {}_I\delta\boldsymbol{r}_C^{\top}\, {}_I\boldsymbol{n}_{,\xi} + {}_B\delta\boldsymbol{\phi}_{IB}^{\top}\left[({}_B\boldsymbol{m})_{,\xi} + {}_B\bar{\boldsymbol{\gamma}} \times {}_B\boldsymbol{n} + {}_B\bar{\boldsymbol{\kappa}}_{IB} \times {}_B\boldsymbol{m}\right]\right\} \mathrm{d}\xi\,. \end{aligned} \tag{A.1}$$

Further, using

$$\begin{aligned}({}_B\boldsymbol{m})_{,\xi} &= (\boldsymbol{A}_{IB}^{\top}\, {}_I\boldsymbol{m})_{,\xi} = \boldsymbol{A}_{IB,\xi}^{\top}\, {}_I\boldsymbol{m} + \boldsymbol{A}_{IB}^{\top}\, {}_I\boldsymbol{m}_{,\xi} \\ &= -{}_B\tilde{\bar{\boldsymbol{\kappa}}}_{IB}\boldsymbol{A}_{IB}^{\top}\, {}_I\boldsymbol{m} + \boldsymbol{A}_{IB}^{\top}\, {}_I(\boldsymbol{m}_{,\xi}) = -{}_B\bar{\boldsymbol{\kappa}}_{IB} \times {}_B\boldsymbol{m} + {}_B(\boldsymbol{m}_{,\xi}) \end{aligned} \tag{A.2}$$

simplifies the internal virtual work to

$$\begin{aligned}\delta W^{\text{int}} &= -\left[{}_I\delta\boldsymbol{r}_C^{\top}\, {}_I\boldsymbol{n} + {}_B\delta\boldsymbol{\phi}_{IB}^{\top}\, {}_B\boldsymbol{m}\right]_{\mathcal{J}} \\ &\quad + \int_{\mathcal{J}} \left\{ {}_I\delta\boldsymbol{r}_C^{\top}\, {}_I\boldsymbol{n}_{,\xi} + {}_B\delta\boldsymbol{\phi}_{IB}^{\top}\left[{}_B(\boldsymbol{m}_{,\xi}) + {}_B\bar{\boldsymbol{\gamma}} \times {}_B\boldsymbol{n}\right]\right\} \mathrm{d}\xi\,. \end{aligned} \tag{A.3}$$

Using this internal virtual work and the external virtual work from (12) and applying the principle of virtual work, which states that the totality of virtual work has to vanish for all virtual displacements and virtual rotations, leads to the strong variational formulation of the rod

$$\delta W^{\text{tot}} = \delta W^{\text{int}} + \delta W^{\text{ext}} \overset{!}{=} 0 \quad \forall_I\delta\boldsymbol{r}_C\,,\ \forall_B\delta\boldsymbol{\phi}_{IB}\,, \tag{A.4}$$

from which the boundary value problem can be easily extracted as

$$\begin{aligned} &{}_I\boldsymbol{n} = -{}_I\boldsymbol{b}_0\,, && {}_B\boldsymbol{m} = -{}_B\boldsymbol{c}_0 && \text{at } \xi = \xi_0\,, \\ &{}_I\boldsymbol{n} = {}_I\boldsymbol{b}_1\,, && {}_B\boldsymbol{m} = {}_B\boldsymbol{c}_1 && \text{at } \xi = \xi_1\,, \\ &{}_I\boldsymbol{n}_{,\xi} + J\,{}_I\boldsymbol{b} = \boldsymbol{0}\,, && {}_B(\boldsymbol{m}_{,\xi}) + {}_B\bar{\boldsymbol{\gamma}} \times {}_B\boldsymbol{n} + J\,{}_B\boldsymbol{c} = \boldsymbol{0}\,, && \left.\right\} \text{in } \operatorname{int}(\mathcal{J})\,. \\ &{}_B\boldsymbol{n} = \boldsymbol{C}_{\boldsymbol{\gamma}}\,\boldsymbol{\varepsilon}_{\boldsymbol{\gamma}}\,, && {}_B\boldsymbol{m} = \boldsymbol{C}_{\boldsymbol{\kappa}}\,\boldsymbol{\varepsilon}_{\boldsymbol{\kappa}} && \end{aligned} \tag{A.5}$$

Note that since ${}_B\bar{\boldsymbol{\gamma}} = \boldsymbol{A}_{IB}^{\top}\, {}_I\boldsymbol{r}_{OC,\xi} = {}_B(\boldsymbol{r}_{OC,\xi})$, this boundary value problem is equivalent to the boundary value problem derived in Eugster (2014). Further note that the last line states the constitutive law in primal form.

In the second case, we perform integration by parts on the second line of the internal virtual work from (11) and use the same abbreviations as above, which leads to

$$\begin{aligned}\delta W^{\text{int,HR}} &= -\int_{\mathcal{J}} \left\{\delta({}_B\boldsymbol{n})^{\top}\left[\bar{\boldsymbol{\varepsilon}}_{\boldsymbol{\gamma}} - J\,\boldsymbol{C}_{\boldsymbol{\gamma}}^{-1}\,{}_B\boldsymbol{n}\right] + \delta({}_B\boldsymbol{m})^{\top}\left[\bar{\boldsymbol{\varepsilon}}_{\boldsymbol{\kappa}} - J\,\boldsymbol{C}_{\boldsymbol{\kappa}}^{-1}\,{}_B\boldsymbol{m}\right]\right\}\mathrm{d}\xi \\ &\quad - \left[{}_I\delta\boldsymbol{r}_C^{\top}\, {}_I\boldsymbol{n} + {}_B\delta\boldsymbol{\phi}_{IB}^{\top}\, {}_B\boldsymbol{m}\right]_{\mathcal{J}} + \int_{\mathcal{J}} \left\{ {}_I\delta\boldsymbol{r}_C^{\top}\, {}_I\boldsymbol{n}_{,\xi} + {}_B\delta\boldsymbol{\phi}_{IB}^{\top}\left[{}_B(\boldsymbol{m}_{,\xi}) + {}_B\bar{\boldsymbol{\gamma}} \times {}_B\boldsymbol{n}\right]\right\} \mathrm{d}\xi\,. \end{aligned} \tag{A.6}$$

Adding the external virtual work from (12) leads to the totality of virtual work, which has now not only to vanish for all ${}_I\delta\boldsymbol{r}_C$ and ${}_B\delta\boldsymbol{\phi}_{IB}$, but also for all $\delta({}_B\boldsymbol{n})$ and $\delta({}_B\boldsymbol{m})$, i.e.,

$$\delta W^{\text{tot}} = \delta W^{\text{int,HR}} + \delta W^{\text{ext}} \overset{!}{=} 0 \quad \forall_I\delta\boldsymbol{r}_C\,,\ \forall_B\delta\boldsymbol{\phi}_{IB}\,,\ \forall\delta({}_B\boldsymbol{n})\,,\ \forall\delta({}_B\boldsymbol{m})\,. \tag{A.7}$$

In this case we can also easily extract the boundary value problem as

$$\begin{aligned} &{}_I\boldsymbol{n} = -{}_I\boldsymbol{b}_0\,, && {}_B\boldsymbol{m} = -{}_B\boldsymbol{c}_0 && \text{at } \xi = \xi_0\,, \\ &{}_I\boldsymbol{n} = {}_I\boldsymbol{b}_1\,, && {}_B\boldsymbol{m} = {}_B\boldsymbol{c}_1 && \text{at } \xi = \xi_1\,, \\ &{}_I\boldsymbol{n}_{,\xi} + J\,{}_I\boldsymbol{b} = \boldsymbol{0}\,, && {}_B(\boldsymbol{m}_{,\xi}) + {}_B\bar{\boldsymbol{\gamma}} \times {}_B\boldsymbol{n} + J\,{}_B\boldsymbol{c} = \boldsymbol{0}\,, && \left.\right\} \text{in } \operatorname{int}(\mathcal{J})\,. \\ &\bar{\boldsymbol{\varepsilon}}_{\boldsymbol{\gamma}} = J\,\boldsymbol{C}_{\boldsymbol{\gamma}}^{-1}\,{}_B\boldsymbol{n}\,, && \bar{\boldsymbol{\varepsilon}}_{\boldsymbol{\kappa}} = J\,\boldsymbol{C}_{\boldsymbol{\kappa}}^{-1}\,{}_B\boldsymbol{m} && \end{aligned} \tag{A.8}$$

Note that this boundary value problem is equivalent to the boundary value problem from (A.5), where the constitutive law in the last line is given in dual form with the strain measures from (1).

# B Details of the finite element formulation

## B.1 Lagrange basis functions

For element $e$ with the interval $\mathcal{J}^e = [\xi^e, \xi^{e+1})$, the Lagrange basis function of polynomial degree $q$ and its derivative are given by

$$N_i^{q,e}(\xi) = \prod_{\substack{0\leq j\leq q \\ j\neq i}} \frac{\xi - \xi_j}{\xi_i - \xi_j}\,, \qquad N_{i,\xi}^{q,e}(\xi) = N_i^{q,e}(\xi) \sum_{\substack{k=0 \\ k\neq i}}^{q} \frac{1}{\xi - \xi_k}\,, \tag{B.1}$$

where $\xi_i = \xi^e + i/q\,(\xi^{e+1} - \xi^e)$ with $i \in \{0, \dots q\}$. Divisions of the form $0/0$ in the case of $q = 0$ are treated as 0, such that $\xi_0 = \xi^e$. One obtains then only one constant basis function $N_0^{0,e}(\xi) = 1$ with $N_{0,\xi}^{0,e}(\xi) = 0$, resulting from the empty product being 1 and the empty sum being 0. This case occurs for the interpolation of the resultant contact forces and moments for a linear element ($p = 1$). Figure B.1 shows the Lagrange basis functions for $q \in \{0, 1, 2\}$.

## B.2 Bookkeeping

This section shows the assembly of the tuple of generalised coordinates $\boldsymbol{q}$ and the tuple of resultant contact forces and moments $\boldsymbol{\lambda}_c$ exemplary for $n_{\mathrm{el}} = 3$ and $p = 2$. Here, we use $\boldsymbol{I}_n \in \mathbb{R}^{n\times n}$ to refer to the identity matrix of $\mathbb{R}^n$. With the total number of nodes $N = 7$, the tuple of generalised coordinates is given by

$$\boldsymbol{q} = \left({}_I\boldsymbol{r}_{OC_0}, \boldsymbol{P}_{IB_0}, \dots {}_I\boldsymbol{r}_{OC_k}, \boldsymbol{P}_{IB_k}, \dots {}_I\boldsymbol{r}_{OC_6}, \boldsymbol{P}_{IB_6}\right) \in \mathbb{R}^{49}\,. \tag{B.2}$$

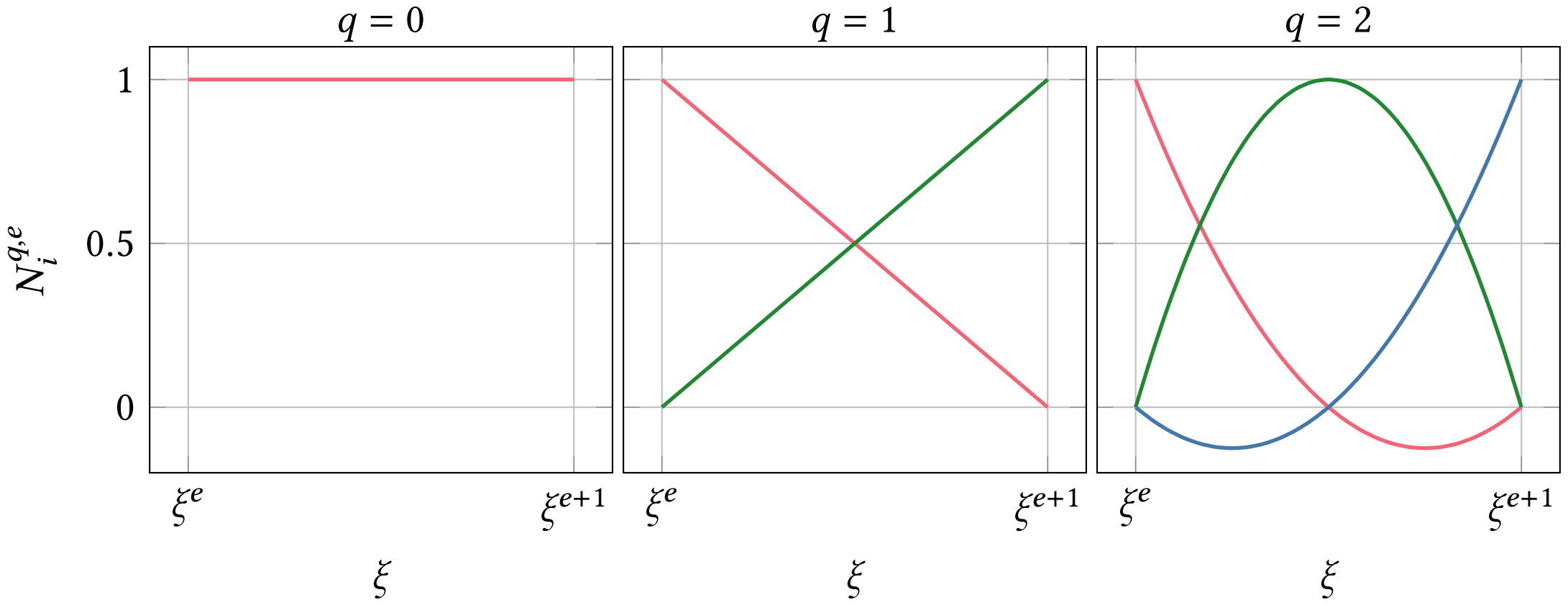

**Figure B.1** Lagrange basis functions $N_i^{q,e}$ for $q = 0$ (left), $q = 1$ (middle) and $q = 2$ (right). The functions are indicated by $N_0^{q,e}$ (red line), $N_1^{q,e}$ (green line) and $N_2^{q,e}$ (blue line).

With the Boolean connectivity matrices

$$\boldsymbol{C}_{\boldsymbol{q},0} = \left(\boldsymbol{I}_{21},\ \boldsymbol{0}_{21\times 28}\right), \quad \boldsymbol{C}_{\boldsymbol{q},1} = \left(\boldsymbol{0}_{21\times 14},\ \boldsymbol{I}_{21},\ \boldsymbol{0}_{21\times 14}\right), \quad \boldsymbol{C}_{\boldsymbol{q},2} = \left(\boldsymbol{0}_{21\times 28},\ \boldsymbol{I}_{21}\right), \tag{B.3}$$

the generalised coordinates of an element are given by

$$\boldsymbol{q}^e = \left({}_I\boldsymbol{r}_{OC_{2e}}, \boldsymbol{P}_{IB_{2e}}, {}_I\boldsymbol{r}_{OC_{2e+1}}, \boldsymbol{P}_{IB_{2e+1}}, {}_I\boldsymbol{r}_{OC_{2e+2}}, \boldsymbol{P}_{IB_{2e+2}}\right) \in \mathbb{R}^{21}\,. \tag{B.4}$$

The matrices of Lagrange shape functions for the element-wise interpolation of the centerline ${}_I\boldsymbol{r}_{OC}$ and the quaternion $\boldsymbol{P}_{IB}$ are then given by

$$\begin{aligned} \boldsymbol{N}_{\boldsymbol{r}}^e &= \begin{pmatrix} N_0^{2,e}\boldsymbol{I}_3 & \boldsymbol{0}_{3\times 4} & N_1^{2,e}\boldsymbol{I}_3 & \boldsymbol{0}_{3\times 4} & N_2^{2,e}\boldsymbol{I}_3 & \boldsymbol{0}_{3\times 4} \end{pmatrix} \in \mathbb{R}^{3\times 21}\,, \\ \boldsymbol{N}_{\boldsymbol{P}}^e &= \begin{pmatrix} \boldsymbol{0}_{4\times 3} & N_0^{2,e}\boldsymbol{I}_4 & \boldsymbol{0}_{4\times 3} & N_1^{2,e}\boldsymbol{I}_4 & \boldsymbol{0}_{4\times 3} & N_2^{2,e}\boldsymbol{I}_4 \end{pmatrix} \in \mathbb{R}^{4\times 21}\,, \end{aligned} \tag{B.5}$$

where the quadratic Lagrange shape functions from Appendix B.1 were used. The tuple of generalised variations $\delta\boldsymbol{s} \in \mathbb{R}^{42}$ and the element tuples $\delta\boldsymbol{s}^e \in \mathbb{R}^{18}$ are ordered identically. Due to their different size, the identity blocks and zero blocks in $\boldsymbol{C}_{\delta s,e} \in \mathbb{R}^{18\times 42}$, $\boldsymbol{N}^e_{\delta\boldsymbol{r}} \in \mathbb{R}^{3\times 18}$ and $\boldsymbol{N}^e_{\delta\boldsymbol{\phi}} \in \mathbb{R}^{3\times 18}$ must be chosen accordingly.

The tuple of resultant contact forces and moments is given by

$$\boldsymbol{\lambda}_c = \left(\boldsymbol{\lambda}_c^0, \boldsymbol{\lambda}_c^1, \boldsymbol{\lambda}_c^2\right) \in \mathbb{R}^{36} \tag{B.6}$$

with the element tuples

$$\boldsymbol{\lambda}_c^e = \left({}_B\boldsymbol{n}_0^e, {}_B\boldsymbol{m}_0^e, {}_B\boldsymbol{n}_1^e, {}_B\boldsymbol{m}_1^e\right) \in \mathbb{R}^{12}\,. \tag{B.7}$$

The Boolean connectivity matrices for the tuples of resultant contact forces and moments are therefore given by

$$\boldsymbol{C}_{\boldsymbol{\lambda}_c,0} = \left(\boldsymbol{I}_{12},\ \boldsymbol{0}_{12\times 24}\right), \quad \boldsymbol{C}_{\boldsymbol{\lambda}_c,1} = \left(\boldsymbol{0}_{12\times 12},\ \boldsymbol{I}_{12},\ \boldsymbol{0}_{12\times 12}\right), \quad \boldsymbol{C}_{\boldsymbol{\lambda}_c,2} = \left(\boldsymbol{0}_{12\times 24},\ \boldsymbol{I}_{12}\right). \tag{B.8}$$

The matrices of Lagrange shape functions for the element-wise interpolation of the resultant contact force ${}_B\boldsymbol{n}$ and moment ${}_B\boldsymbol{m}$ are given by

$$\begin{aligned} \boldsymbol{N}^e_{\boldsymbol{n}} &= \begin{pmatrix} N_0^{1,e}\boldsymbol{I}_3 & \boldsymbol{0}_{3\times 3} & N_1^{1,e}\boldsymbol{I}_3 & \boldsymbol{0}_{3\times 3} \end{pmatrix} \in \mathbb{R}^{3\times 12}\,, \\ \boldsymbol{N}^e_{\boldsymbol{P}} &= \begin{pmatrix} \boldsymbol{0}_{3\times 3} & N_0^{1,e}\boldsymbol{I}_3 & \boldsymbol{0}_{3\times 3} & N_1^{1,e}\boldsymbol{I}_3 \end{pmatrix} \in \mathbb{R}^{3\times 12}\,, \end{aligned} \tag{B.9}$$

where the linear Lagrange shape functions from Appendix B.1 were used. The ordering of the tuples of generalised coordinates, generalised variations and resultant contact forces and moments can be effortlessly extended to different numbers of elements and different polynomial degrees. Accordingly, the sizes of the Boolean connectivity matrices and the matrices of shape function changes, but their structure is similarly.

## B.3 Initialisation of the rod

Let the initial configuration of the rod be described by ${}_I\boldsymbol{r}^*_{OC}(\xi) \in \mathbb{R}^3$ and $\boldsymbol{A}^*_{IB}(\xi) \in SO(3)$. The nodal positions and quaternions at the $N$ nodes $\xi_k = k/(N-1)$ with $k \in \{0, \dots N-1\}$ are then given such that

$${}_I\boldsymbol{r}_{OC_k} = {}_I\boldsymbol{r}^*_{OC}(\xi_k)\,, \quad \boldsymbol{A}(\boldsymbol{P}_{IB_k}) = \boldsymbol{A}^*_{IB}(\xi_k)\,, \quad \|\boldsymbol{P}_{IB_k}\| = 1\,, \quad \boldsymbol{P}^\top_{IB_k}\boldsymbol{P}_{IB_{k+1}} > 0\,. \tag{B.10}$$

The extraction of the quaternion from a given transformation matrix can be done using Spurrier's algorithm (see Spurrier 1978). The last two conditions ensure at first the fulfillment of the quaternion constraint $\boldsymbol{g}_S(\boldsymbol{q}) = \boldsymbol{0}_N$, and secondly that two consecutive quaternions are in the same hemisphere, since there exist two unit quaternions leading to the same transformation matrix, i.e., $\boldsymbol{A}(\boldsymbol{P}) = \boldsymbol{A}(-\boldsymbol{P})$.

# C Tangent operator for quaternion

To give the derivation of the tangent operator for non-unit quaternions in a comprehensible way, we introduce at first two properties related to the cross product. Afterwards, we use these properties to precompute quantities that will be used in the last step, where the tangent operator is finally derived.

In the following, let $\boldsymbol{a}, \boldsymbol{b}, \boldsymbol{c}, \boldsymbol{v} \in \mathbb{R}^3$ be arbitrary triples. Starting with Grassmann's identity $\boldsymbol{a} \times (\boldsymbol{b} \times \boldsymbol{c}) = (\boldsymbol{a}^\top\boldsymbol{c})\boldsymbol{b} - (\boldsymbol{a}^\top\boldsymbol{b})\boldsymbol{c}$, we can concatenate its result with another cross product from the left side, which results in

$$\tilde{\boldsymbol{a}}\tilde{\boldsymbol{b}}\tilde{\boldsymbol{c}}\boldsymbol{v} = \boldsymbol{a} \times (\boldsymbol{b} \times (\boldsymbol{c} \times \boldsymbol{v})) = \boldsymbol{a} \times (\boldsymbol{c}(\boldsymbol{b}^\top\boldsymbol{v}) - \boldsymbol{v}(\boldsymbol{b}^\top\boldsymbol{c})) = \boldsymbol{a} \times \boldsymbol{c}\,(\boldsymbol{b}^\top\boldsymbol{v}) - (\boldsymbol{b}^\top\boldsymbol{c})\boldsymbol{a} \times \boldsymbol{v} \tag{C.1}$$

As this has to hold for all $\boldsymbol{v} \in \mathbb{R}^3$, we can conclude

$$\tilde{\boldsymbol{a}}\tilde{\boldsymbol{b}}\tilde{\boldsymbol{c}} = (\boldsymbol{a} \times \boldsymbol{c})\boldsymbol{b}^\top - (\boldsymbol{b}^\top\boldsymbol{c})\tilde{\boldsymbol{a}}\,. \tag{C.2}$$

The second property is derived by the difference

$$(\tilde{\boldsymbol{a}}\tilde{\boldsymbol{b}} - \tilde{\boldsymbol{b}}\tilde{\boldsymbol{a}})\boldsymbol{v} = \boldsymbol{a} \times (\boldsymbol{b} \times \boldsymbol{v}) - \boldsymbol{b} \times (\boldsymbol{a} \times \boldsymbol{v}) = \left\{\boldsymbol{b}(\boldsymbol{a}^\top \boldsymbol{v}) - \boldsymbol{v}(\boldsymbol{a}^\top \boldsymbol{b})\right\} - \left\{\boldsymbol{a}(\boldsymbol{b}^\top \boldsymbol{v}) - \boldsymbol{v}(\boldsymbol{b}^\top \boldsymbol{a})\right\} . \tag{C.3}$$

Since the latter terms in the braces cancel each other out, we are left with only the first terms. They can be simplified by

$$\boldsymbol{b}(\boldsymbol{v}^\top \boldsymbol{a}) - \boldsymbol{a}(\boldsymbol{v}^\top \boldsymbol{b}) = \boldsymbol{v} \times (\boldsymbol{b} \times \boldsymbol{a}) = (\boldsymbol{a} \times \boldsymbol{b}) \times \boldsymbol{v} = (\widetilde{\tilde{\boldsymbol{a}}\boldsymbol{b}})\boldsymbol{v} , \tag{C.4}$$

where we have used Grassmann's identity in the first step and the anticommutativity of the cross product in the second step. Again, as this has to hold for all $\boldsymbol{v} \in \mathbb{R}^3$, we can use the $\tilde{\bullet}$ notation and conclude

$$\tilde{\boldsymbol{a}}\tilde{\boldsymbol{b}} - \tilde{\boldsymbol{b}}\tilde{\boldsymbol{a}} = (\widetilde{\tilde{\boldsymbol{a}}\boldsymbol{b}}) . \tag{C.5}$$

Let $\boldsymbol{P} = (p_0, \boldsymbol{p}) \in \mathbb{R}^4$ with $p_0 \in \mathbb{R}$ and $\boldsymbol{p} \in \mathbb{R}^3$ be a non-unit quaternion. Following Egeland and Gravdahl (2003, Equation 6.166) with normalisation, the orthonormal transformation matrix is defined by

$$\boldsymbol{A}_{IB} = \boldsymbol{I} + \frac{2}{\|\boldsymbol{P}\|^2}(p_0\tilde{\boldsymbol{p}} + \tilde{\boldsymbol{p}}\tilde{\boldsymbol{p}}) . \tag{C.6}$$

At first, some later required quantities are precomputed. Using the vectorial part of the quaternion in (C.2) results in

$$\tilde{\boldsymbol{p}}\tilde{\boldsymbol{p}}\tilde{\boldsymbol{p}} = -(\boldsymbol{p}^\top \boldsymbol{p})\tilde{\boldsymbol{p}} = p_0^2\tilde{\boldsymbol{p}} - \|\boldsymbol{P}\|^2\tilde{\boldsymbol{p}} , \tag{C.7}$$

where we have used $\|\boldsymbol{P}\|^2 = p_0^2 + \boldsymbol{p}^\top \boldsymbol{p}$ in the second step. This is used to compute

$$\boldsymbol{A}_{IB}\,\tilde{\boldsymbol{p}} = \tilde{\boldsymbol{p}} + \frac{2}{\|\boldsymbol{P}\|^2}(p_0\tilde{\boldsymbol{p}}\tilde{\boldsymbol{p}} + \tilde{\boldsymbol{p}}\tilde{\boldsymbol{p}}\tilde{\boldsymbol{p}}) = -\tilde{\boldsymbol{p}} + \frac{2p_0}{\|\boldsymbol{P}\|^2}(p_0\tilde{\boldsymbol{p}} + \tilde{\boldsymbol{p}}\tilde{\boldsymbol{p}}) = -\tilde{\boldsymbol{p}} + p_0(\boldsymbol{A}_{IB} - \boldsymbol{I}) . \tag{C.8}$$

Using the trivial identity $\boldsymbol{A}_{IB}^\top = \boldsymbol{I} - (\boldsymbol{I} - \boldsymbol{A}_{IB}^\top)$ yields

$$\boldsymbol{A}_{IB}^\top\,\tilde{\boldsymbol{p}}' = \tilde{\boldsymbol{p}}' - (\boldsymbol{I} - \boldsymbol{A}_{IB}^\top)\,\tilde{\boldsymbol{p}}' . \tag{C.9}$$

Post-multiplying this with $\tilde{\boldsymbol{p}}$, using $\tilde{\boldsymbol{p}}'\tilde{\boldsymbol{p}} = \boldsymbol{p}(\boldsymbol{p}')^\top - \boldsymbol{I}\boldsymbol{p}^\top\boldsymbol{p}'$ and further $(\boldsymbol{I} - \boldsymbol{A}_{IB}^\top)\boldsymbol{p} = \boldsymbol{0}$ yields

$$\boldsymbol{A}_{IB}^\top\,\tilde{\boldsymbol{p}}'\tilde{\boldsymbol{p}} = \tilde{\boldsymbol{p}}'\tilde{\boldsymbol{p}} + \boldsymbol{p}^\top\boldsymbol{p}'(\boldsymbol{I} - \boldsymbol{A}_{IB}^\top) . \tag{C.10}$$

Almost identical to (C.8), we obtain $\boldsymbol{A}_{IB}^\top\,\tilde{\boldsymbol{p}} = -\tilde{\boldsymbol{p}} - p_0(\boldsymbol{A}_{IB}^\top - \boldsymbol{I})$. This is used to get the last precomputed quantity as

$$\boldsymbol{A}_{IB}^\top\,\tilde{\boldsymbol{p}}\,\tilde{\boldsymbol{p}}' = -\tilde{\boldsymbol{p}}\,\tilde{\boldsymbol{p}}' - p_0(\boldsymbol{A}_{IB}^\top - \boldsymbol{I})\,\tilde{\boldsymbol{p}}' . \tag{C.11}$$

With these quantities, let us now have a look at the derivative of the quaternion mapping. It is given by

$$\boldsymbol{A}_{IB}' = -\boldsymbol{P}^\top\boldsymbol{P}'\frac{2}{\|\boldsymbol{P}\|^2}\left(\frac{2}{\|\boldsymbol{P}\|^2}(p_0\tilde{\boldsymbol{p}} + \tilde{\boldsymbol{p}}\tilde{\boldsymbol{p}}) + \boldsymbol{I} - \boldsymbol{I}\right) + \frac{2}{\|\boldsymbol{P}\|^2}(p_0'\tilde{\boldsymbol{p}} + p_0\tilde{\boldsymbol{p}}' + \tilde{\boldsymbol{p}}'\tilde{\boldsymbol{p}} + \tilde{\boldsymbol{p}}\tilde{\boldsymbol{p}}') , \tag{C.12}$$

where we inserted $+\boldsymbol{I} - \boldsymbol{I}$ in the first bracket to identify the whole bracket as $\boldsymbol{A}_{IB} - \boldsymbol{I}$. Further expanding $\boldsymbol{P}^\top\boldsymbol{P}' = p_0p_0' + \boldsymbol{p}^\top\boldsymbol{p}'$ and grouping by $p_0'$ and $\boldsymbol{p}'$ yields

$$\boldsymbol{A}_{IB}' = \frac{2}{\|\boldsymbol{P}\|^2}\left[p_0'\big(\tilde{\boldsymbol{p}} - p_0(\boldsymbol{A}_{IB} - \boldsymbol{I})\big) + \big(p_0\tilde{\boldsymbol{p}}' + \tilde{\boldsymbol{p}}'\tilde{\boldsymbol{p}} + \tilde{\boldsymbol{p}}\tilde{\boldsymbol{p}}' - \boldsymbol{p}^\top\boldsymbol{p}'(\boldsymbol{A}_{IB} - \boldsymbol{I})\big)\right] . \tag{C.13}$$

The left round bracket can be identified as $-\boldsymbol{A}_{IB}\,\tilde{\boldsymbol{p}}$, as shown in (C.8). The product of $\boldsymbol{A}_{IB}^\top$ with the right round bracket can be easily evaluated with the quantities from (C.9) to (C.11) and is given by

$$\begin{aligned} &\boldsymbol{A}_{IB}^\top\big(p_0\tilde{\boldsymbol{p}}' + \tilde{\boldsymbol{p}}'\tilde{\boldsymbol{p}} + \tilde{\boldsymbol{p}}\tilde{\boldsymbol{p}}' - \boldsymbol{p}^\top\boldsymbol{p}'(\boldsymbol{A}_{IB} - \boldsymbol{I})\big) \\ &= p_0\tilde{\boldsymbol{p}}' - p_0(\boldsymbol{I} - \boldsymbol{A}_{IB}^\top)\tilde{\boldsymbol{p}}' + \tilde{\boldsymbol{p}}'\tilde{\boldsymbol{p}} + \boldsymbol{p}^\top\boldsymbol{p}'(\boldsymbol{I} - \boldsymbol{A}_{IB}^\top) \\ &\quad - \tilde{\boldsymbol{p}}\tilde{\boldsymbol{p}}' - p_0(\boldsymbol{A}_{IB}^\top - \boldsymbol{I})\tilde{\boldsymbol{p}}' - \boldsymbol{p}^\top\boldsymbol{p}'(\boldsymbol{I} - \boldsymbol{A}_{IB}^\top) \\ &\qquad = p_0\tilde{\boldsymbol{p}}' - (\tilde{\boldsymbol{p}}\tilde{\boldsymbol{p}}' - \tilde{\boldsymbol{p}}'\tilde{\boldsymbol{p}}) = p_0\tilde{\boldsymbol{p}}' - (\widetilde{\tilde{\boldsymbol{p}}\boldsymbol{p}'}) , \end{aligned} \tag{C.14}$$

where a lot of terms cancel in the second expression, and the identity (C.5) was used in the last step. Finally, the skew symmetric matrix of the scaled curvature is obtained by

$$ {}_B\tilde{\bar{\boldsymbol{\kappa}}}_{IB} = \boldsymbol{A}_{IB}^{\top}\boldsymbol{A}_{IB}' = \frac{2}{\|\boldsymbol{P}\|^2}\big[-p_0'\tilde{\boldsymbol{p}} + p_0\tilde{\boldsymbol{p}}' - (\widetilde{\tilde{\boldsymbol{p}}\boldsymbol{p}'})\big]\,. \tag{C.15} $$

From that, one can extract the vector representing the scaled curvature as

$$ {}_B\bar{\boldsymbol{\kappa}}_{IB} = \frac{2}{\|\boldsymbol{P}\|^2}\big[-p_0'\boldsymbol{p} + p_0\boldsymbol{p}' - \tilde{\boldsymbol{p}}\boldsymbol{p}'\big] = \boldsymbol{T}(\boldsymbol{P})\boldsymbol{P}'\,, \tag{C.16} $$

where the tangent operator is given by

$$ \boldsymbol{T}(\boldsymbol{P}) = \frac{2}{\|\boldsymbol{P}\|^2}\begin{pmatrix}-\boldsymbol{p} & p_0\boldsymbol{I} - \tilde{\boldsymbol{p}}\end{pmatrix} \in \mathbb{R}^{3\times 4}\,. \tag{C.17} $$

This result is identical to the result given by Rucker (2018) and can also be found in Wasmer and Betsch (2024).

# D Numerical convergence analysis

## D.1 Spatial convergence

To quantify the accuracy of a formulation, the spatial convergence is used. For that, the root-mean-square error of relative twists $e_{\boldsymbol{\theta}}^k$ is computed. The superscript $k$ refers to the number of points $0 \leqslant \xi_i = i/(k-1) \leqslant 1$, where the relative twist $\Delta\boldsymbol{\theta}(\xi_i)$ is computed. The relative twist contains a unified error measure for positions and cross-section orientations. For a given Euclidean transformation matrix

$$ \boldsymbol{H}_{\mathcal{I}\mathcal{B}} = \begin{pmatrix}\boldsymbol{A}_{IB} & {}_I\boldsymbol{r}_{OC}\\ \boldsymbol{0}_{1\times 3} & 1\end{pmatrix}\,, \tag{D.1} $$

the corresponding twist can be obtained by $\boldsymbol{\theta} = \mathrm{Log}_{SE(3)}(\boldsymbol{H}_{\mathcal{I}\mathcal{B}}) \in \mathbb{R}^6$, where the definition of the $SE(3)$-logarithm map can be found in Eugster and Harsch (2023, Equation 32).

The error measure is given by

$$ e_{\boldsymbol{\theta}}^k = \sqrt{1/k\sum_{i=0}^{k-1}\Delta\boldsymbol{\theta}(\xi_i)^{\top}\Delta\boldsymbol{\theta}(\xi_i)}\,,\quad \Delta\boldsymbol{\theta}(\xi_i) = \mathrm{Log}_{SE(3)}\big(\boldsymbol{H}_{\mathcal{I}\mathcal{B}}(\xi_i)^{-1}\boldsymbol{H}_{\mathcal{I}\mathcal{B}}^*(\xi_i)\big) \in \mathbb{R}^6\,, \tag{D.2} $$

where $\boldsymbol{H}_{\mathcal{I}\mathcal{B}}^*(\xi_i)$ is the Euclidean transformation matrix of the reference solution and $\boldsymbol{H}_{\mathcal{I}\mathcal{B}}(\xi_i)$ is the Euclidean transformation of the rod to be compared. We will analyse how the root-mean-square error of relative twists $e_{\boldsymbol{\theta}}^k$ changes, when the number of nodes $N$, respectively the number of elements $n_{\mathrm{el}}$, is increased.

## D.2 Solver convergence

To quantify the computational performance of the different formulations, the number of iterations and the interim results of Newton's method are analysed. Let the $k$-th Newton-iterate of all unknowns be given by $\boldsymbol{x}_{[k]} \in \mathbb{R}^n$ and let $m$ be the number of required iterations to fulfil the termination criteria $\sqrt{\boldsymbol{f}(\boldsymbol{x})^{\top}\boldsymbol{f}(\boldsymbol{x})} < \varepsilon\sqrt{n}$. Based on Albersmeyer and Diehl (2010), we approximate the local quadratic convergence rate by

$$ r \approx \frac{\|\boldsymbol{x}_{[m-1]} - \boldsymbol{x}_{[m]}\|}{\|\boldsymbol{x}_{[m-2]} - \boldsymbol{x}_{[m]}\|^2}\,. \tag{D.3} $$

As multiple load increments are used, we can compute the arithmetic mean $\overline{m}_{\mathrm{a}}$ and the standard deviation $\sigma_{m_{\mathrm{a}}}$ of the number of Newton iterations during all load steps. Since the local quadratic convergence rate is always positive, it is convenient to visualise it with a logarithmic scale. We can therefore use the geometric mean of $\overline{r}_{\mathrm{g}}$ and the geometric standard deviation $\sigma_{r_{\mathrm{g}}}$ of the convergence rates during all load steps. Note that the natural logarithm of the geometric mean is equivalent to the arithmetic mean of the natural logarithm of each local quadratic convergence rate. Similarly, the natural logarithm of the geometric standard deviation is equivalent to the standard deviation of the natural logarithm of each local quadratic convergence rate. Therefore, no variance in the local quadratic convergence rates corresponds to a geometric standard deviation of 1.

**Authors' contributions** MH: Conceptualisation, Methodology, Software, Formal Analysis, Visualisation Investigation, Writing - original draft, Writing - review & editing. DC: Methodology, Investigation, Software, Writing - original draft. JB: Methodology, Software, Writing - review & editing. ICG: Formal Analysis, Supervision, Writing - review & editing. LG: Methodology, Supervision, Writing - review & editing. SRE: Conceptualisation, Methodology, Supervision, Writing - review & editing.

**Competing interests** The authors declare that they have no competing interests.

**Journal's Note** JTCAM remains neutral with regard to the content of the publication and institutional affiliations.